\documentclass[english,11pt,a4paper]{amsart}
\usepackage[colorlinks=true,linkcolor=blue,citecolor=red,urlcolor=blue]{hyperref}
 \usepackage[left=1.05in,right=1.05in,bottom=1in,top=1in]{geometry}
 \usepackage{dsfont}

\usepackage{tikz,amsthm,amsmath,amstext,amssymb,amscd,epsfig,euscript, mathrsfs, dsfont,pspicture,multicol, mathtools,graphpap,graphics,graphicx,times,enumerate,subfig,sidecap,wrapfig,color}
 \usepackage{url}
\makeatletter
\def\@tocline#1#2#3#4#5#6#7{\relax
  \ifnum #1>\c@tocdepth 
  \else
    \par \addpenalty\@secpenalty\addvspace{#2}%
    \begingroup \hyphenpenalty\@M
    \@ifempty{#4}{%
      \@tempdima\csname r@tocindent\number#1\endcsname\relax
    }{%
      \@tempdima#4\relax
    }%
    \parindent\z@ \leftskip#3\relax \advance\leftskip\@tempdima\relax
    \rightskip\@pnumwidth plus4em \parfillskip-\@pnumwidth
    #5\leavevmode\hskip-\@tempdima
      \ifcase #1
       \or\or \hskip 1em \or \hskip 2em \else \hskip 3em \fi%
      #6\nobreak\relax
    \dotfill\hbox to\@pnumwidth{\@tocpagenum{#7}}\par
    \nobreak
    \endgroup
  \fi}
\makeatother

\def \rr {{\mathbb R}}
\def \cc {{\mathbb C}}
\def \nn {{\mathbb N}}

\DeclareGraphicsExtensions{.png}
\def\z{\zeta}
\def\ep{\epsilon}
\def\u{\upsilon}
\def\g{\gamma}
\def\l{\lambda}
\def\k{\kappa}
\def\e{\eta}

\def\s{\sigma}
\def\m{\mu}

\def\D{\Delta}
\def\G{\Gamma}

\def\O{\Omega}

\def\S{\Sigma}

\def \cc {{\mathbb{C}}}

\newcommand{\bela}{\begin{equation} \label}
\newcommand{\eeq}{\end{equation}}
\newcommand{\ba}{\begin{array}}
\newcommand{\ea}{\end{array}}

\newtheorem{definition}{Definition}[section]
\newtheorem{theorem}{Theorem}[section]
\newtheorem{proposition}{Proposition}[section]
\newtheorem{lemma}{Lemma}[section]
\newtheorem{corollary}{Corollary}[section]
\newtheorem{remark}{Remark}[section]
\newtheorem{notation}{Notation}[section]

\catcode`\@=12

\begin{document}


\title[Dirac operators with $\delta$-Shell Interactions]{Spectral Properties of the Dirac Operator coupled with $\delta$-Shell Interactions}

\author[Badreddine Benhellal]{Badreddine Benhellal}
\address{Departamento de Matem\'aticas, Universidad del Pa\' is Vasco, Barrio Sarriena s/n 48940 Leioa, SPAIN\\ and
 Universit\' e de Bordeaux, IMB, UMR 5251, 33405 Talence Cedex, FRANCE.}
\email{benhellal.badreddine@ehu.eus and badreddine.benhellal@u-bordeaux.fr}


\begin{abstract}
Let $\O\subset\rr^3$ be an open set. We study the spectral properties of the free Dirac operator  $ \mathcal{H} :=- i \alpha \cdot\nabla + m\beta$ coupled with the  singular potential  $V_\k=(\ep I_4 +\mu\beta + \eta(\alpha\cdot N))\delta_{\partial\O}$, where $\k=(\ep,\m,\eta)\in\rr^3$. The open set $\O$ can be either a $\mathcal{C}^2$-bounded domain or a locally deformed half-space. In both cases, self-adjointness is proved and several spectral properties are given. In particular, we give a complete description of the essential spectrum of $ \mathcal{H}+V_\k $ in the case of a  locally deformed half-space, for the so-called critical combinations of coupling constants.  Finally,  we introduce a new model of Dirac operators with $\delta$-interactions  and deal with its spectral properties.  More precisely, we study the coupling $\mathcal{H}_{\z,\u}=\mathcal{H}+  \left( -i\z\alpha_1\alpha_2\alpha_3+ i\u\beta\left(\alpha\cdot N\right)\right)\delta_{\partial\O}$, with $\z,\u\in\rr$. In particular, we show that  $\mathcal{H}_{ 0,\pm2}$ is essentially self-adjoint and generates confinement. 
\end{abstract}

\maketitle
\tableofcontents

\section{Introduction}

 In this paper,  we investigate  in $\rr^3$ the self-adjointness character and the spectral properties of the free  Dirac operator $\mathcal{H}$ coupled with  combinations of the following singular potentials (see section \ref{sec1} for notations):
 \begin{align*} 
V_{\ep}=\ep I_4\delta_{\partial\O},\quad  V_{\mu}=\mu\beta\delta_{\partial\O},  \quad  V_\e=\eta(\alpha\cdot N)\delta_{\partial\O},\quad  \ep,\m,\e\in\rr,\\
V_{\z}=  -i\z\alpha_1\alpha_2\alpha_3\delta_{\partial\O}, \quad V_{\u}=  i\u\beta \left(\alpha\cdot N\right)\delta_{\partial\O},  \quad \z,\u\in\rr,
\end{align*} 
 here  $\partial\O$ is the boundary of an open set $\O$ of $\rr^3$, $\mathit{N}$ is the unit normal vector field at $\partial\O$ which points outwards of $\O$ and  the $\delta$-potential is the Dirac distribution supported on $\partial\O$.  In relativistic quantum mechanics,  the Dirac Hamiltonian $\mathcal{H}+ \mathbb{V}$,  where $\mathbb{V}$ is a combination of the above potentials,  describes the dynamics of the massive relativistic particles of spin-$1/2$  in the external potential $\mathbb{V}$. From this physical point of view, the singular interactions given by the coupling constants $\ep,\m,\e$ and $\u$  are called respectively electrostatic, Lorentz scalar, magnetic and anomalous magnetic potential; as far as we know the singular interaction given by the coupling constant $\z$ is introduced for the first time in this article, and we will see that it can be considered as an interactions of electrostatic type. Finally, the surface $\partial\O$ supporting the interactions is called a shell. 

Recently,  Dirac operators  with $\delta$-shell interactions have been studied extensively. Namely,  the coupling $ \mathcal{H}$ with the electrostatic and the Lorentz scalar $\delta$-shell interactions which we denote by $\mathcal{H}_{\ep,\m}$; we refer to the survey \cite{OBP} for a review on the topic. To our knowledge, the spectral study of the Dirac operator  $\mathcal{H}_{\ep,\m}$  goes back to the papers \cite{DES} and \cite{DA}, where the authors studied the spherical case (i.e $\partial\O$ is a sphere). Moreover, in \cite{DES} the authors point out that under the assumption $\ep^2-\m^2=-4$, the shell becomes impenetrable. Physically, this means that at the time $t = 0$, if the particle in consideration (an electron for example) is in the region $\O$ (respectively in $\rr^3\setminus\overline{\O}$), then during the evolution in time, it cannot cross the surface $\partial\O$ to join the region $\rr^3\setminus\overline{\O}$ (respectively $\O$ ) for all $t> 0$. Mathematically, this means that the Dirac operator in consideration  decouples into a direct sum of two Dirac operators acting respectively on $ \O$ and  $\rr^3\setminus\overline{\O}$,  with appropriate boundary conditions. In particular, when $\ep=0$, this phenomenon  has been known to physicists since the $1970$'s (cf.\cite{CJJTW} and \cite{John} for example); and its mathematical model described by the Dirac operator with MIT boundary conditions has been the subject of several mathematical papers (we refer to the recent paper \cite{BHM} as well as the references cited there). 
All these physical motivations made the mathematical study of Dirac operators with $\delta$-shell interactions a relevant  subject. However, unlike the non-relativistic counterpart (i.e. Schr\"{o}dinger operators with $\delta$-shell interactions) the study of relativistic $\delta$-interactions has known a long period of silence. Indeed, apart from \cite{SV} where the authors studied the scattering theory and the non-relativistic limit of $\mathcal{H}_{\ep,\m}$ (in the spherical case), the spectral study of $\mathcal{H}_{\ep,\m}$ has been forgotten for two decades.  Since then, it has been relaunched in  \cite{AMV1}, where the authors developed a new technique to characterize the self-adjointness of the free Dirac operator coupled with a measure-valued potential. As a particular case, they dealt with the pure electrostatic $\delta$-shell interactions (i.e $\m=0$)  supported on the boundary of a bounded regular domain, and  they proved that the perturbed operator is self-adjoint for all $\ep\neq\pm2$. The same authors continued the spectral study of  the electrostatic
case; for instance, the existence of point spectrum and related problems; see \cite{AMV2} and \cite{AMV3}. In \cite{AMV2}  it is shown that $\mathcal{H}_{\ep,\m}$  still generates confinement under the condition $\ep^2-\m^2=-4$. 

A different approach based on the concept of quasi-boundary triples and their Weyl functions has been  used in \cite{BEHL1}  to study the Dirac operators with electrostatic $\delta$-shell interactions. In there, the authors prove the self-adjointness for all $\ep\neq\pm2$, and investigate several spectral properties, adding the scattering theory and asymptotic properties of the model. In all the above papers, the case $\ep=\pm2$ (known as the critical interaction strengths) has not been considered. This gap has been covered in \cite{OBV} and \cite{BH} with different approaches. Indeed, in this particular case, it turns out that the Dirac operator with electrostatic $\delta$-shell interactions is essentially self-adjoint, and functions in the domain of the closure are less regular comparing to the non critical case. Moreover,  in \cite{BH}  it is  shown that if $\partial\O$ contains a flat part,  then the point $0$ belongs to the essential spectrum of $\overline{\mathcal{H}_{\pm2,0}}$. Similar phenomenon appears when we study the Dirac operator $\mathcal{H}_{\ep,\m}$. In fact, in this case, the critical combinations of coupling constants are  $\ep^2-\m^2=4$; see \cite{BEHL2} for example. The self-adjointness in this critical case  was proved for the two dimensional analogue of $\mathcal{H}_{\ep,\m}$ in \cite{BHOP}, where the authors considered $\delta$-interactions supported on  a smooth closed curve. Furthermore, by making use of complex analysis and periodic pseudo-differential operators techniques, they showed  that
\begin{align}\label{QQ1}
\mathrm{Sp}_{\mathrm{ess}}(\overline{\mathcal{H}_{\ep,\m}})=\big(-\infty,-m\big]\cup\left\{-\frac{m\m}{\ep}\right\}\cup \big[m,+\infty\big).
\end{align}
 Of course, such techniques are no longer available in the three dimensional case. Nevertheless, at this stage,  one may ask the following question:
\begin{itemize}
  \item[(Q1)] In the three dimensional setting, when $\ep^2-\m^2=4$, does \eqref{QQ1} hold true?
\end{itemize}

The main objective of the current manuscript is to study the spectral properties and the phenomenon of the confinement  of the following couplings:
\begin{align*}
\mathcal{H}_{\k}&:=\mathcal{H}+ (\ep I_4 +\mu\beta+\eta(\alpha\cdot N))\delta_{\partial\O},  \quad \k:=(\ep,\m,\e)\in\rr^3,\\
\mathcal{H}_{\z,\u}&:=\mathcal{H} + \left( -i\z\alpha_1\alpha_2\alpha_3+ i\u\beta\left(\alpha\cdot N\right)\right)\delta_{\partial\O},\quad (\z,\u)\in\rr^2.
\end{align*} 

  Let us present the context we are considering and summarize the main results of our work.  We shall  assume that the open set $\O$  satisfies one of the following hypotheses:
\begin{itemize}
  \item[(1)]   $\O$ is a $\mathit{C}^2$-bounded domain.
  \item[(2)]  $\O:= \O_\nu:=  \lbrace (x,t)\in \rr^2\times\rr:  t> \nu \phi(x)\rbrace$, where $\nu \in\rr$ and $\phi : \rr^2\rightarrow \rr$ is a $\mathit{C}^2$-smooth, compactly supported function.
  \end{itemize}
 Following the strategy of \cite{AMV1}, we define the Dirac operators  $\mathcal{H}_{\bullet}$, $\bullet =\k$ or $(\z,\u)$,  on the domain 
 \begin{align*}
\mathrm{dom}(\mathcal{H}_{\bullet})=\left\{ u+\Phi[g]: u\in\mathit{H}^1(\rr^3)^4, g\in\mathit{L}^2(\partial\O)^4, u_{|\partial\O}=-\Lambda_{+}[g]\right\},
\end{align*}
where $\Phi$ is an appropriate fundamental solution of the unperturbed operator $\mathcal{H}$, and  $\Lambda_{\pm}$ are bounded linear operators  acting on $\mathit{L}^2(\partial\O)^4$   (see Notation \ref{Lambda}). We mention that the operator  $\Lambda_{\pm}$ also appears in several works when the quasi-boundary triples theory is used to study the Dirac operator $\mathcal{H}_{\ep,\mu}$, see  \cite[Lemma 5.4]{BH} and   \cite[Proposition 4.3]{BHOP} for example. We point out that  the consideration of the second assumption is motivated by \cite{EL}, where the Schr\"{o}dinger operator with $\delta$-shell interaction was considered. 

 As a first step of the current paper, we study the self-adjointness character of  $\mathcal{H}_{\k}$, when $\O$ satisfies the assumption $(1)$ or $(2)$. We begin by  proving that  $\mathcal{H}_{\k}$ is self-adjoint when $\ep^2-\m^2-\e^2\neq4$ (i.e in the non-critical case), and we show that $\mathrm{dom}(\mathcal{H}_{\k})\subset\mathit{H}^1(\rr^3\setminus\partial\O)^4$, which means that functions in $\mathrm{dom}(\mathcal{H}_{\k})$ have a Sobolev regularity; cf. Theorem \ref{main1}.  To prove this result we develop a strategy very close to \cite{OBV}, it is based essentially on the fact that the anticommutators of Cauchy operator $ \mathit{C}_{\partial\O}$ (see \eqref{Cauchyop} for the definition) with $\beta$ or with $(\alpha\cdot N)$ have a regularizing effect. Indeed, as it was observed in several works  (see \cite{AMV2} for example), the operators $\Lambda_{\mp}\Lambda_{\pm}$ involve the above anticommutators and it turns out that in the non-critical case, the regularization effect of these anticommutator pushes $\Lambda_{+}$ to regularize the functions in $\mathrm{dom}(\mathcal{H}_{\k})$ to have the $\mathit{H}^1$-Sobolev regularity.  When $\ep^2-\m^2-\e^2=4$, which is actually the critical case, we show that $\mathcal{H}_{\k}$ is essentially self-adjoint (i.e $\overline{\mathcal{H}_{\k}}$ is self-adjoint). In addition, we point out the relation between the self-adjointness of $\mathcal{H}_{\k}$  and the operator $\Lambda_{+}$, which is essentially the main idea behind the concept of quasi boundary triples theory (see Subsection \ref{sub3.2}).

 As a second step, we turn to the spectral study of $\mathcal{H}_{\k}$,  we focus on the case where $\O$  satisfies the second assumption and we show several spectral properties of  $\mathcal{H}_{\k}$. Namely, using Fourier analysis and compactness arguments,  we compute precisely the essential spectrum of  $\mathcal{H}_{\k}$ in the non-critical case for $\eta=0$; cf Theorem \ref{Prop non}. More precisely, under certain conditions on the sign of $\ep,\m$ and $(\ep^3-\m^2)$, it turns out that the continuous spectrum might emerge in the gap $(-m,m)$,  which leads to the phenomenon of the essential spectrum
instability. In particular, we have $\mathrm{Sp}_{\mathrm{ess}}(\mathcal{H}_{\k})=\rr$  when $\m=-2$ and $\ep=\eta=0$; see Theorem \ref{Prop non} for more details. 
 
  In the critical case, we give a complete characterization of the essential spectrum of $\overline{\mathcal{H}_{\k}}$  when  $\O$ satisfies the second assumption. More precisely, we prove in Theorem \ref{cas deforme}  that 
$$\mathrm{Sp}_{\mathrm{ess}}(\overline{\mathcal{H}_{\k}})=\big(-\infty,-m\big]\cup\left\{-\frac{m\m}{\ep}\right\}\cup \big[m,+\infty\big),\quad \ep^2-\m^2-\eta^2=4,$$
which answers positively to the question $(\mathrm{Q1})$, hence generalizing the result of  \cite{BHOP} to this kind of surfaces. The proof is based on the use of  compactness and localization arguments. We remark that even after adding the perturbation by the potential $V_{\e}$,  the point which appears in the gap remains the same (see the discussion after Theorem \ref{cas deforme} for more details).  All these results  will be proven using an adapted Birman-Schwinger principle, a Krein-type resolvent formula and compactness arguments. Nevertheless, the situation is more delicate, in particular the use of compactness arguments. 

The last part of this paper is devoted to the spectral study of the Dirac operator $\mathcal{H}_{\z,\u}$.  We mention that while preparing this manuscript, it turns out that the authors of the paper \cite{CLMT} considered the two-dimensional analog of  $\mathcal{H}_{0,\u}$, and our results intersect on this point (see Section \ref{sec6} for more details). Assuming that $\O$ satisfies the assumption $(1)$, one of the most important properties that we show for this operator is that, in the critical case $\z^2+\u^2=4$,  $\mathcal{H}_{ \z,\u}$ is essentially self-adjoint. Furthermore, in one hand $\mathcal{H}_{ \pm\z,0}$ coincides with  the Dirac operator coupled with the electrostatic $\delta$-interactions of strength $-\z$. In another hand,  $\mathcal{H}_{ 0,\pm\u}$ decouples in a direct sum of two Dirac operators acting respectively on $ \O$ and  $\rr^3\setminus\overline{\O}$, with boundary conditions in $\mathit{H}^{-1/2}(\partial\O)$. Thus,  $\overline{\mathcal{H}_{ 0,\u}}$ generates confinement for $\u=\pm2$, and hence $\partial\O$ becomes impenetrable. Moreover, the inner part of $\overline{\mathcal{H}_{ 0,\pm\u}}$ which acts on $ \O$ coincide with so-called Dirac operator with \textit{Zig-zag} boundary condition, see Section \ref{sec6}. 

We mention that several statements on the self-adjointness  and spectral properties of the above operators have been extended in \cite{BB} in the non-critical cases for compact surfaces with low regularity. 

\textbf{Organisation of the paper.} The structure of the paper is as follows. In the second section, we set up the necessary notations and recall the relevant material from \cite{AMV1}. In Section \ref{sec3}, we study the self-adjointness of  $\mathcal{H}_{\k}$, when $\partial\O$ satisfies the first and the second assumption, the main results being Theorem \ref{main1}. Section \ref{sec4} is devoted to the spectral study of $\mathcal{H}_{\k}$.  We focus namely on the case where $\O$ is a locally deformed half space and we give a complete description of the essential spectrum of $\mathcal{H}_{\k}$, for the non-critical and critical combinations of coupling constants in Theorem \ref{Prop non} and Theorem \ref{cas deforme}, respectively.   Finally,  in Section \ref{sec6},  we study the spectral properties of the Dirac operator $\mathcal{H}_{\z,\u}$, for all possible combinations of interaction strengths. The main results in this section are Theorem   \ref{main6} and Theorem \ref{main7}.

\section{Notations and Preliminaries} \label{sec1}
\setcounter{equation}{0}

We consider a surface $\Sigma\subset\rr^3$ dividing the space into two regions $\O_{\pm}$. More precisely, we assume that $\S$ satisfies one of the hypotheses:
\begin{itemize}
  \item[(H1)]  $\S=\partial\Omega_+$ with $\Omega_+$ a $\mathit{C}^2$-bounded domain.
  \item[(H2)]  $\S:= \S_\nu:=  \lbrace (x_1,x_2,x_3)\in \rr^3:  x_3= \nu \phi(x_1,x_2)\rbrace$, where $\nu \in\rr_{+}$ and $\phi : \rr^2\rightarrow \rr$ is a $\mathit{C}^2$-smooth, compactly supported function. We denote by  $F$ the flat  part of $\S_\nu$ i.e.
  \begin{align}\label{plat}
   F:= \{ x=(x_1,x_2,\nu \phi(x_1,x_2))\in \S_\nu: (x_1,x_2)\notin \mathrm{supp}(\phi) \}.
   \end{align}
\end{itemize} 
We parameterize $\S_\nu$ by the mapping 
\begin{equation}\label{param}
   \left\{
\begin{aligned}
\tau: \rr^2 &\longrightarrow \rr^3\\
\overline{x}&\longmapsto (\overline{x},\nu \phi(\overline{x}))
\end{aligned}
  \right.
\end{equation}
For $x=(\overline{x},\nu \phi(\overline{x}))\in\S_\nu$, we express the surface mesure on $\S_\nu$ via the formula $\mathrm{dS}(x)=J_\nu(\overline{x})\mathrm{d}\overline{x}$, where  $J_\nu$ is the Jacobian given by 
\begin{align}\label{leJ}
   J_\nu(\overline{x})= \sqrt{1+ \nu^2| \nabla \phi(\overline{x})|^2}.
\end{align}
Throughout the paper,  we shall  work on the Hilbert space $\mathit{L}^2(\rr^{d})^4$ (respectivelly, $\mathit{L}^2(\O_\pm)^4$) with respect to the Lebesgue measure, and we will make use of the orthogonal decomposition $\mathit{L}^2(\rr^3)^4= \mathit{L}^2(\O_{+})^4\oplus\mathit{L}^2(\O_{-})^4$. That is, for $\varphi\in\mathit{L}^2(\rr^3)^4$ we write $\varphi=(\varphi_+,\varphi_-)$,  where $\varphi_\pm=\varphi\mathds{1}_{\O_\pm}=:\varphi\downharpoonright_{\Omega_\pm}$. $\mathcal{D}(\O_\pm)^4$ denotes the usual space of indefinitely
differentiable functions with compact support, and $\mathcal{D}^{\prime}(\O_\pm)^4$ is the  space of distributions defined as the dual space of $\mathcal{D}(\O_\pm)^4$. We define the unitary Fourier-Plancherel operator $ \mathcal{F}:\mathit{L}^2(\rr^d)^4\rightarrow \mathit{L}^2(\rr^{d})^4$ as follows
\begin{align}
 \mathcal{F}[u](\xi)=\frac{1}{(2\pi)^{d/2}}\int_{\rr^d}e^{-ix\cdot\xi}u(x)\mathrm{d}x,\quad \forall\xi \in\rr^d, 
\end{align}
and by $ \mathcal{F}^{-1}$ we denote the inverse  Fourier-Plancherel operator $ \mathcal{F}^{-1}:\mathit{L}^2(\rr^{d})^4\rightarrow \mathit{L}^2(\rr^{d})^4$, given by
\begin{align}
 \mathcal{F}^{-1}[u](x)=\frac{1}{(2\pi)^{d/2}}\int_{\rr^d}e^{i\xi\cdot x}u(\xi)\mathrm{d}\xi,\quad \forall x \in\rr^{d}.
\end{align}
Given $\overline{x}\in\rr^{d-1}$, by  $\mathcal{F}_{\overline{x}}$   we abbreviate the partial Fourier-Plancherel operator on the variable  $\overline{x}$.  Given $s\in[-1,1]$, we denote by  $\mathit{H}^s(\rr^{d})^4$ the Sobolev space of  order $s$, defined as
\begin{align}
 \mathit{H}^s(\rr^{d})^4:=\{ u\in\mathit{L}^2(\rr^{d})^4: \int_{\rr^d}(1+|\xi|^2)^{s} \left|\mathcal{F}[u](\xi)\right|^2\mathrm{d}\xi<\infty\}.
\end{align}
The Sobolev space $\mathit{H}^1(\O_\pm)^4$ is defined as follows:
\begin{align*}
\mathit{H}^1(\O_\pm)^4=\{ \varphi\in\mathit{L}^2(\O_\pm)^4: \text{ there exists } \tilde{\varphi}\in\mathit{H}^1(\rr^{3})^4 \text{ such that }  \tilde{\varphi}|_{\O_\pm} =\varphi\},
\end{align*}
and we  use notation 
$$\mathit{H}^1(\rr^3\setminus\S)^4=\mathit{H}^1(\Omega_+)^4\oplus\mathit{H}^1(\Omega_-)^4.$$ 
By $\mathit{L}^2(\S, \mathrm{d}S)^4:=\mathit{L}^2(\S)^4$ we denote the usual $\mathit{L}^2$-space over $\S$.  Given  $s\in[0,1]$, if $\S$ satisfies $(\mathrm{H2})$, we then define the Sobolev spaces $\mathit{H}^s(\S)^4$  in terms of the Sobolev spaces over $\rr^2$ as usual. That is given $g\in \mathit{L}^2(\S)^4$, we define $g_{\phi}(\overline{x}) =g(\overline{x},\nu\phi(\overline{x}))$, for $\overline{x}\in\rr^{2}$. Then 
\begin{align}
 \mathit{H}^s(\S)^4:=\{ g\in\mathit{L}^2(\S)^4:g_{\phi}\in \mathit{H}^s(\rr^2)^4\} , \text{ for all } s\in[0,1],
\end{align}
and then define $\mathit{H}^{-s}(\S)^4$ to be the completion of $\mathit{L}^2(\S)^4$ with following norm:
\begin{align}
\| g\|_{\mathit{H}^{-s}(\S)^4}:= \|g_{\phi} J_\nu \|_{\mathit{H}^{-s}(\rr^2)^4}, \text{ for all } s\in[0,1].
\end{align}
Recall that $\mathit{H}^{-s}(\S)^4$  is a realization  of the dual space of $\mathit{H}^{s}(\S)^4$; see \cite{Mc} for example. Now,  if $\S$ satisfies $(\mathrm{H1})$, we then define the Sobolev spaces  $\mathit{H}^s(\S)^4$ using local
coordinates representation on the surface $\S$; see \cite{Mc}. We will use the notation $\langle\, , \, \rangle_{\mathit{H}^{-1/2},\mathit{H}^{1/2}}$ for the duality pairing between $\mathit{H}^{-1/2}(\S)^4$ and $\mathit{H}^{1/2}(\S)^4$.
By $t_{\S}:\mathit{H}^1(\O_\pm)^4\rightarrow \mathit{H}^{1/2}(\S)^4$ we denote the classical trace operator, and by   $E_{\O_\pm}:  \mathit{H}^{1/2}(\S)^4\rightarrow  \mathit{H}^1(\O_\pm)^4$ the extension operator, i.e $t_{\S}E_{\O_\pm}$ is the identity operator.  For a function $u\in\mathit{H}^1(\rr^3)^4$, with a slight abuse of terminology we will refer to $t_{\S}u:=u\downharpoonright_{\S}\in\mathit{H}^{1/2}(\S)^4$ as the restriction of $u$ on $\S$, and by $E:  \mathit{H}^{1/2}(\S)^4\rightarrow  \mathit{H}^1(\rr^3)^4$ the continuous right inverse.

Let $x\in\S$ and $a>0$, denote the nontangential approach regions of opening $a$ at the point $x$ by
\begin{align}\label{ntregion}
\Gamma_a^{\O_\pm}(x)=\{ y\in\O_\pm: |x-y|<(1+a)\mathrm{dist}(y,\S)\}.
\end{align}
We fix $a >0$ large enough such that $x\in\overline{\Gamma_a^{\O_\pm}(x)}$ for all $x\in\S$. If $x\in\S$, then 
\begin{align}\label{ntTrace}
U_\pm(x):= \lim\limits_{\Gamma_{a}^{\O_\pm}(x) \ni y\xrightarrow[]{}x}U(y),
\end{align}
is the nontangential limit of $U$ with respect to $\O_\pm$ at $x$.  If  $a>0$ is fixed, we shall write $\Gamma^{\O_\pm}(x)$  instead of $\Gamma_{a}^{\O_\pm}(x)$.

  Let $\alpha= (\alpha_1,\alpha_2,\alpha_3)$ and $\beta$ be the $4\times4$ Hermitian and unitary matrices given by
\begin{align}\label{les matrice}
\alpha_k=\begin{pmatrix}
0 & \sigma_k\\
\sigma_k & 0
\end{pmatrix}\quad \text{ for } k=1,2,3
\quad
\beta=\begin{pmatrix}
\mathit{I}_2 & 0\\
0 & -\mathit{I}_2
\end{pmatrix},
\end{align}  
where $\sigma = (\sigma_1,\sigma_2,\s_3)$ are the Pauli matrices defined by
\begin{align}\label{Pauli}
\sigma_1=\begin{pmatrix}
0 & 1\\
1 & 0
\end{pmatrix},\quad \sigma_2=
\begin{pmatrix}
0 & -i\\
i & 0
\end{pmatrix} ,
\quad
\sigma_3=\begin{pmatrix}
1 & 0\\
0 & -1
\end{pmatrix}.
\end{align} 
We denote by $\mathit{N}$ and  $\delta_\S$  the unit normal vector field at $\S$ which points outwards of $\O_+$ and the Dirac distribution supported on $\S$, respectively. Given $m>0$, we consider the Dirac operator 
\begin{align}
\mathcal{H}_{\k} =  \mathcal{H} + V_{\k}=- i \alpha \cdot\nabla + m\beta+(\ep I_4 +\mu\beta + \eta(\alpha\cdot N))\delta_{\S},\quad  \k:=(\ep,\m,\e)\in\rr^3,
\end{align}  
in the Hilbert space $\mathit{L}^2(\rr^3)^4$. We recall that $(\mathcal{H},\mathit{H}^1(\rr^3)^4)$ is self-adjoint (see \cite[ subsection~1.4]{Tha}) and its spectrum is given by
\begin{align*}
  \mathrm{Sp}(\mathcal{H}) =\mathrm{Sp}_{\mathrm{ess}}(\mathcal{H})=(-\infty,-m]\cup [m,+\infty),
\end{align*}
here and in the rest of the paper, for a closed operator T,  its  resolvent set,  spectrum, essential spectrum, point and discrete spectrum are denoted by  $\rho(T), \mathrm{Sp}(T),\mathrm{Sp}_{\mathrm{ess}}(T),\mathrm{Sp}_{\mathrm{pp}}(T)$ and $\mathrm{Sp}_{\mathrm{disc}}(T)$, respectively.

The rest of this section will be devoted to give a first definition of the Hamiltonian $\mathcal{H}_{\k}$ . For this and for the convenience of the reader, we recall the relevant material from \cite{AMV1} (without detailed proofs), thus making our exposition self-contained. 


\subsection{Integral operators associated to the Dirac operator}
Here we list some well known results about integral operators associated to the fundamental solution of the Dirac operator. Given $z\in\cc\setminus\left((-\infty,-m]\cup[m,\infty)\right)$  with the convention that $\mathrm{Im}\sqrt{z^2-m^2}>0$, we recall that the fundamental solution of $(\mathcal{H}-z)$ is given by
\begin{align}
\label{}
  \phi^z(x)=\frac{e^{i\sqrt{z^2-m^2}|x|}}{4\pi|x|}\left(z +m\beta+( 1-i\sqrt{z^2-m^2}|x|)i\alpha\cdot\frac{x}{|x|^2}\right), \quad \text{for all } x\in\rr^3\setminus\{0\},
\end{align}
  see for example \cite[Section 1.E]{Tha}. Next, we define the following operators
\begin{align}\label{Phi}
\begin{split}
\Phi^z: \mathit{L}^2(\S)^4&  \longrightarrow \mathit{L}^2(\rr^3)^4\\
g&\longmapsto  \Phi^{z}[g](x) =\int_\S \phi^{z}(x-y)g(y)\mathrm{dS}(y), \quad \text{for all } x\in\rr^3\setminus\S,
\end{split}
\end{align}
then, $\Phi^z: \mathit{L}^2(\S)^4  \longrightarrow \mathit{L}^2(\rr^3)^4$ is a bounded operator. Furthermore,  $(\mathcal{H}-z)\Phi^{z}[g]=0$ holds in $\mathcal{D}^{\prime}(\Omega_\pm)^4$, for all $g \in\mathit{L}^2(\S)^4$. Moreover, thanks to  \cite[Proposition 4.2]{BH}, we know that $\Phi^z$ gives  rise to a bounded operator from $\mathit{H}^{1/2}(\S)^4 $ onto $\mathit{H}^1(\rr^3\setminus\S)^4$.  It is worth mentioning that the arguments leading to \cite[Proposition 4.2]{BH} remain valid when $\S$ satisfies $(\mathrm{H2})$, and do not depend on the fact that $\Omega_+$ is bounded. see also Remark \ref{une remarque} for an alternative proof.

  Given $x\in\S$ and $g \in\mathit{L}^2(\S)^4$, we set
\begin{align}\label{Cauchyop}
  \mathit{C}^{z}_{\S}[g](x)=  \lim\limits_{\rho\searrow 0}\int_{|x-y|>\rho}\phi^{z}(x-y)g(y)\mathrm{dS}(y)\quad \text{ and } \quad \mathit{C}^{z}_{\pm}[g](x)=  \lim\limits_{\Gamma^{\O_\pm}(x) \ni y\xrightarrow[]{}x}\Phi^{z}[g](y).
\end{align}
Then, we have the following lemma.

\begin{lemma}\label{lemme 2.1} Let $\mathit{C}^{z}_{\S}$ and $\mathit{C}^{z}_{\pm}$ be as above. Then $ \mathit{C}^{z}_{\S}[g](x)$ and $\mathit{C}^{z}_{\pm}[g](x)$ exist for $\mathrm{dS}$-a.e. $x\in\S$, and $\mathit{C}^{z}_{\S}, \mathit{C}^{z}_{\pm}: \mathit{L}^2(\S)^4\rightarrow \mathit{L}^2(\S)^4$ are linear bounded operators. Furthermore, the following hold:
\begin{itemize}
  \item[(i)] $\mathit{C}^{z}_{\pm}= \mp\frac{i}{2}(\alpha\cdot\mathit{N}) + \mathit{C}^{z}_{\S} $ ,(Plemelj-Sokhotski jump formula).
  \item[(ii)]   $(\mathit{C}^{z}_{\S} (\alpha\cdot\mathit{N}) )^2 =-\frac{1}{4}\mathit{I}_4$. In particular,  $\Vert \mathit{C}^{z}_{\S}\Vert\geqslant \frac{1}{2}$.
\end{itemize}
\end{lemma}
\textbf{Proof.}   If $\S$ satisfies $(\mathrm{H1})$, then the proof is  analogous to the one of  \cite[Lemma 2.2]{AMV2},   where the authors use essentially the Green's theorem and the following well known result on the trace of derivatives of a single-layer potential. Indeed, given $g\in \mathit{L}^2(\S)$, then for $\mathrm{dS}$-a.e. $x\in\S$, we have
\begin{align}\label{jusmprelation}
   \lim\limits_{\Gamma^{\O_\pm}(x)  \ni y\xrightarrow[]{}x}\int\frac{y-w}{4\pi |y-w|^3}g(w)\mathrm{dS}(w)= \mp\frac{1}{2}g(x)\mathit{N}(x) +  \lim\limits_{\rho\searrow 0}\int_{|x-w|>\rho}\frac{x-w}{4\pi |x-w|^3}g(w)\mathrm{dS}(w).
\end{align}
Note that this result is also true if  $\S$ satisfies $(\mathrm{H2})$, see \cite[Theorem 5.4.7]{Me} for example. Thus, one can adapt the proof of \cite[Lemma 3.3]{AMV1} in this case and get the claimed results, we omit the details.\qed
\newline

\begin{remark}\label{remark2.1} Note that in the same setting, Lemma \ref{lemme 2.1} still holds true if for example $\S$ is a compact Lipschitz surface or the graph of a Lipschitz function $\phi : \rr^2\rightarrow \rr$; see  \cite{AGHM} and \cite[Remark 3.14]{AMV1}. Moreover, since $(\Phi^{\bar{z}})^{\ast}=(\mathcal{H}-z)^{-1}\downharpoonright_{\S}$,  by duality and interpolation  arguments, it follows that $\Phi^z$ gives rise to a bounded operator from $\mathit{L}^{2}(\S)^4 $ onto $\mathit{H}^{1/2}(\rr^3\setminus\S)^4$; cf. \cite[Subsection 3.3]{BHM}.  Hence, the non-tangential limit in Lemma \ref{lemme 2.1}$(\mathrm{i})$ coincides with the trace operator for all data in $\mathit{H}^{1/2}(\S)^4$.  
\end{remark}
\begin{corollary} Let $z\in\cc\setminus\left((-\infty,-m]\cup[m,\infty)\right)$. Then, the operator $\mathit{C}^{z}_{\S}$ is bounded from  $\mathit{H}^{1/2}(\S)^4$ onto itself. Moreover, it holds that $(\mathit{C}^{z}_{\S})^{\ast}[g]= \mathit{C}^{\overline{z}}_{\S}[g]$, for all  $g\in\mathit{L}^2(\S)^4$. In particular, $\mathit{C}^{z}_{\S}$ is  self-adjoint in $\mathit{L}^2(\S)^4$, for all $z\in(-m,m)$.
\end{corollary}
\textbf{Proof.} Given $g\in \mathit{H}^{1/2}(\S)^4$.  Since $\Phi^z[g]\in \mathit{H}^{1}(\rr^3\setminus\S)^4$,  it follows that  $\mathit{C}^{z}_{\pm}[g]\in \mathit{H}^{1/2}(\S)^4$. Thus,  from Lemma \ref{lemme 2.1} $(\mathrm{i})$  we deduce that $2\mathit{C}^{z}_{\S} [g]= (\mathit{C}^{z}_{+}+\mathit{C}^{z}_{-})[g]\in \mathit{H}^{1/2}(\S)^4$. This proves the first statement. The second statement is a direct consequence of the fact that $\overline{\phi^{z}(y-x)}=\phi^{\overline{z}}(x-y)$.  \qed 
\newline

\begin{notation}\label{Lambda} For $\k=(\ep,\m,\e)\in\rr^3$, we set  
\begin{align}\label{sgndef}
 \mathrm{sgn}(\k):= \ep^2- \m^2-\e^2\neq0.
 \end{align}
  We define the operators $\Lambda^{z}_{\pm}$ as follows:
\begin{align}
\Lambda^{z}_{\pm}=\frac{1}{\mathrm{sgn}(\k)}(\ep I_4 \mp(\mu\beta+ \eta(\alpha\cdot N)) )\pm\mathit{C}^{z}_{\S},\quad\forall z\in\cc\setminus\left((-\infty,-m]\cup[m,\infty)\right).
\end{align}
Since $(\alpha\cdot N)$ is $\mathit{C}^{1}$-smooth and symmetric, it easily follows  that $\Lambda^{z}_{\pm}$ are bounded (and self-adjoint for $z\in(-m,m)$) from  $\mathit{L}^{2}(\S)^4$ onto itself, and bounded from $\mathit{H}^{1/2}(\S)^4$ onto itself.
\end{notation}
In the sequel, we shall write $\Phi$, $ \mathit{C}_{\S} $,  $\mathit{C}_{\pm}$ and $\Lambda^{}_{\pm}$ instead of $\Phi^0$, $ \mathit{C}^{0}_{\S} $,  $ \mathit{C}^{0}_{\pm}$ and $\Lambda^{0}_{\pm}$, respectively. Now we are in position to give the first definition of the Dirac Hamiltonian with $\delta$-interactions supported on $\S$, the main object of the present paper.

\begin{definition}\label{def1} Let $\k=(\ep,\m,\e)\in\rr^3$ be such that $\mathrm{sgn}(\k)\neq0 $. The Dirac operator coupled with a combination of electrostatic, Lorentz scalar and normal vector field $\delta$-shell interactions of strength $\ep$, $\m$ and $\e$ respectively, is the operator $\mathcal{H}_{\k}=  \mathcal{H} + V_{\k}$, acting in $ \mathit{L}^2(\rr^3)^4$ and defined on the domain
 \begin{align}\label{dom}
\mathrm{dom}(\mathcal{H}_{\k})=\left\{ u+\Phi[g]: u\in\mathit{H}^1(\rr^3)^4, g\in\mathit{L}^2(\S)^4, t_{\S}u=-\Lambda_{+}[g]\right\},
\end{align}
where 
\begin{align}
\label{}
V_\k(\varphi)=\frac{1}{2}(\ep I_4 +\mu\beta+ \eta(\alpha\cdot N)))(\varphi_+ +\varphi_-)\delta_{\S},
\end{align}
with $\varphi_\pm= t_\S u + \mathit{C}^{}_\pm [g]$. Hence,  $\mathcal{H}_{\k}$ acts in the sense of distributions as $ \mathcal{H}_{\k}(\varphi)= \mathcal{H}u$, for all $\varphi=u+\Phi^{}[g]\in\mathrm{dom}(\mathcal{H}_{\k})$. 
\end{definition}
\setcounter{equation}{0}
\section{Self-adjointness of $ \mathcal{H}_{\k}$  }\label{sec3}
In this section, we study the self-adjointness of the Dirac operator $ \mathcal{H}_{\k}$.  In our setting, it turns out that the special value $\mathrm{sgn}(\k)=4$ plays a critical role in the analysis of the spectral properties of  $ \mathcal{H}_{\k}$.  Before stating the main result of this part, some notations and auxiliary results are needed. 

Denote by $\mathit{H}(\alpha,\O_\pm)$ the Sobolev space associated to the Dirac operator on $\O_\pm$, defined by
 \begin{align}\label{Sobolev Dirac}
\mathit{H}(\alpha,\O_\pm)=\{ \varphi\in\mathit{L}^2(\O_\pm)^4: (\alpha\cdot\nabla)\varphi\in\mathit{L}^2(\O_\pm)^4 \},
\end{align}
then $\mathit{H}(\alpha,\O_\pm)$ is a Hilbert space with respect to the following scalar product (see \cite[Section 2.3]{OBV}) 
\begin{align*}
\langle \varphi, \psi \rangle_{\mathit{H}(\alpha,\O_\pm)}=\langle \varphi, \psi\rangle_{\mathit{L}^{2}(\O_\pm)^4}+\langle (-i\alpha\cdot\nabla)\varphi, (-i\alpha\cdot\nabla)\psi\rangle_{\mathit{L}^{2}(\O_\pm)^4},\quad \varphi,\psi\in\mathit{H}(\alpha,\O_\pm).
\end{align*}

The following proposition gathers some properties of the Sobolev space $\mathit{H}(\alpha,\O_\pm)$ and the operators  $\Phi^z$ and $\mathit{C}^z_{\S}$.
 \begin{proposition}\label{extension}(\cite{OBV},\cite{BH}) Let $\Phi^z$, $\mathit{C}^{z}_{\pm}$ and $\mathit{C}^z_{\S}$ be as in Lemma \ref{lemme 2.1}. Then, the following hold:
\begin{itemize}
 \item[(i)] The trace operator $t_{\S}$ (which until now was defined on  $\mathit{H}^{1}(\Omega_\pm)^4$) has a unique extension to a bounded linear operator from $\mathit{H}(\alpha,\O_\pm)$ to  $\mathit{H}^{-1/2}(\S)^4$, and we have 
 \begin{align}\label{Geen formula}
\langle (-i\alpha\cdot\nabla)\varphi, \psi \rangle_{\mathit{L}^{2}(\O_\pm)^4}-\langle \varphi,  (-i\alpha\cdot\nabla)\psi\rangle_{\mathit{L}^{2}(\O_\pm)^4}=\pm\langle (-i\alpha\cdot\mathit{N})t_{\S}\varphi, t_{\S}\psi\rangle_{\mathit{H}^{-1/2},\mathit{H}^{1/2}},
\end{align}
 for all $\varphi\in\mathit{H}(\alpha,\O_\pm)$ and $\psi\in\mathit{H}^{1}(\Omega_\pm)^4$.
  \item[(ii)] The operator $\Phi^{z}$ admits a continuous extension from $\mathit{H}^{-1/2}(\S)^4$ to $\mathit{H}(\alpha,\O_+)\oplus \mathit{H}(\alpha,\O_-)$, which we still denote by $\Phi^{z}$.
  \item[(iii)] The operator $\mathit{C}^{z}_{\S}$ admits a continuous extension  $\tilde{\mathit{C}^{z}_{\S}}:\mathit{H}^{-1/2}(\S)^4\rightarrow\mathit{H}^{-1/2}(\S)^4$. Moreover, we have
  \begin{align}\label{dual}
  \begin{split}
  t_{\S}(\Phi^z[h]\downharpoonright_{\Omega_\pm})=:\tilde{\mathit{C}^{z}_{\pm}}[h] = (\mp\frac{i}{2}(\alpha\cdot\mathit{N}) +\tilde{\mathit{C}^{z}_{\S}})[h],\\
        \langle \tilde{\mathit{C}^{z}_{\S}}[h],g\rangle_{\mathit{H}^{-1/2},\mathit{H}^{1/2}}=\langle h,\mathit{C}^{\overline{z}}_{\S}[g]\rangle_{\mathit{H}^{-1/2},\mathit{H}^{1/2}},
        \end{split}
    \end{align}
    for any $g\in\mathit{H}^{1/2}(\S)^4$ and  $ h\in\mathit{H}^{-1/2}(\S)^4$.
      \item[(iv)] If $\varphi\in\mathit{H}(\alpha,\O_\pm)$  and $t_{\S}\varphi\in\mathit{H}^{1/2}(\S)^4$, then $ \varphi\in\mathit{H}^{1}(\Omega_\pm)^4$.  Moreover, for all $\varphi_\pm\in\mathit{H}(\alpha,\O_\pm)$ one has $ \left(1/2 \mp i\tilde{\mathit{C}^z_{\S}}(\alpha\cdot\mathit{N})\right)t_{\S} \varphi_\pm\in\mathit{H}^{1/2}(\S)^4$.
\end{itemize}
\end{proposition}
\textbf{Proof.} First, assume that  $\S$ satisfies $(\mathrm{H1})$. Then the statements of assertion $(\mathrm{i})$ follow by  \cite[Proposition 2.1]{OBV} and  \cite[Corollary 2.15]{OBV} respectively. Assertion $(\mathrm{ii})$ can be proved as much the same way as in \cite[Theorem 2.2]{OBV}, see also \cite[Proposition 4.4]{BH}. Since $(\mathit{C}^{z}_{\S})^{\ast}= \mathit{C}^{\overline{z}}_{\S}$, and  $\mathit{C}^{z}_{\S}$ is  bounded from   $\mathit{H}^{1/2}(\S)^4$ onto itself, by duality we get the first statement of $(\mathrm{iii})$, and  \eqref{dual} follows by density arguments,  for a detailed proof we refer to \cite[Proposition 3.5]{BHOP} and \cite[Proposition 4.4 $(\mathrm{ii})$]{BH}.  Assertion $(\mathrm{iv})$ follows by \cite[Proposition 2.7 and Proposition 2.16]{OBV} when $z=0$, and the case $z\neq0$ follows in an analogous way.  Finally, we note that the  arguments leading to the above results  depend only on the continuity of $t_{\S}$ on  $\mathit{H}^{1}(\Omega_\pm)^4$ and the density arguments which remain valid when $\S$ satisfies $(\mathrm{H2})$, thus assertions $(\mathrm{i})-(\mathrm{iv})$ are still hold true in that case.  
\qed
\newline
\begin{remark}\label{une remarque} Another way to prove the boundedness of  $\Phi^z$  from $\mathit{H}^{1/2}(\S)^4 $ into $\mathit{H}^1(\rr^3\setminus\S)^4$ is as follows.  Firstly, use  density arguments to extend continuously $\mathit{C}^{z}_{\S}$   from  $\mathit{H}^{-1/2}(\S)^4$ into itself, and then by duality from  $\mathit{H}^{1/2}(\S)^4$ into itself. Secondly, extend the trace formula \eqref{dual} by density. Finally, as $\mathit{N}$ is $\mathit{C}^1$-smooth, we then get the claimed result by Proposition \ref{extension}-$(\mathrm{iv})$.
\end{remark}

 In the following, we shall denote by $\tilde{\Lambda}^{z}_{\pm}$ the continuous extension of $\Lambda^{z}_{\pm}$ defined from $\mathit{H}^{-1/2}(\S)^4$ onto itself.  Now,  we can state the first main theorem of the paper, the remainder of this part will be devoted to the proof of this result.
 
 \begin{theorem}\label{main1} Let $ \mathcal{H}_{\k}$ be as in the definition \ref{def1}. Then, the following statements hold true:
\begin{itemize}
 \item[(i)] If  $\mathrm{sgn}(\k)\neq 4$, then $ \mathcal{H}_{\k}$ is self-adjoint and we have 
  \begin{align}
  \mathrm{dom}(\mathcal{H}_{\k})=\left\{ u+\Phi[g]: u\in\mathit{H}^1(\rr^3)^4, g\in\mathit{H}^{1/2}(\S)^4, t_{\S}u=-\Lambda_+[g]\right\}.
  \end{align}
  \item[(ii)]  If  $\mathrm{sgn}(\k)= 4$, then $\mathcal{H}_{\k}$ is essentially self-adjoint and we have 
  \begin{align}
  \mathrm{dom}(\overline{\mathcal{H}_{\k}})=\left\{ u+\Phi[g]: u\in\mathit{H}^1(\rr^3)^4, g\in\mathit{H}^{-1/2}(\S)^4, t_{\S}u=-\tilde{\Lambda}_+[g]\right\}.
  \end{align}
 \end{itemize}
\end{theorem}
     
\begin{proposition}\label{closable} Let  $ \mathcal{H}_{\k}$ be as in the definition \ref{def1}. Then, $ \mathcal{H}_{\k}$ is closable.
\end{proposition}
\textbf{Proof.}  As any symmetric operator on a Hilbert space with dense domain of definition always admits a closure, to prove the proposition it suffices to show the following:
\begin{itemize}
  \item[(i)]  $\mathrm{dom}(\mathcal{H}_{\k})$ is dense in $\mathit{L}^{2}(\rr^3)^4$.
  \item[(ii)]  $\mathcal{H}_{\k}$ is symmetric on $\mathrm{dom}(\mathcal{H}_{\k})$.
\end{itemize}
   First, observe that $\mathit{C}^{\infty}_{0}(\rr^3\setminus\S)^4\subset\mathrm{dom}(\mathcal{H}_{\k})\subset\mathit{L}^{2}(\rr^3)^4$. Thus $(\mathrm{i})$ follows from this and the fact that $\mathit{C}^{\infty}_{0}(\rr^3\setminus\S)^4$ is a dense subspace of $\mathit{L}^{2}(\rr^3)^4$. Now we prove  $(\mathrm{ii})$, let $\varphi,\psi\in\mathrm{dom}(\mathcal{H}_{\k})$ with $\varphi=u+\Phi[g]$ and $\psi= v+\Phi[h]$. Then, we have
\begin{align*}\label{}
  \langle \mathcal{H}_{\k}\varphi,\psi\rangle_{\mathit{L}^2(\rr^3)^4}- \langle \varphi,\mathcal{H}_{\k}\psi \rangle_{\mathit{L}^2(\rr^3)^4}&=  \langle \mathcal{H}u,v+\Phi[h]\rangle_{\mathit{L}^2(\rr^3)^4}- \langle u+\Phi[g],\mathcal{H}v \rangle_{\mathit{L}^2(\rr^3)^4}\\ 
& =  \langle \mathcal{H}u,\Phi[h]\rangle_{\mathit{L}^2(\rr^3)^4}- \langle \Phi[g],\mathcal{H}v \rangle_{\mathit{L}^2(\rr^3)^4}\\
    &=\langle t_{\S}u, h \rangle_{\mathit{L}^2(\S)^4}- \langle g,t_{\S}v \rangle_{\mathit{L}^{2}(\S)^4}.
\end{align*}
Using the conditions $ t_{\S}u=-\Lambda_+[g]$ and $ t_{\S}v=-\Lambda_+[h]$,  and that $\Lambda_+$ is self-adjoint, we obtain
\begin{align}\label{relationself}
  \langle \mathcal{H}_{\k}\varphi,\psi\rangle_{\mathit{L}^2(\rr^3)^4}- \langle \varphi,\mathcal{H}_{\k}\psi \rangle_{\mathit{L}^2(\rr^3)^4} =\langle -\Lambda_+[g],h\rangle_{\mathit{L}^2(\S)^4}- \langle g,-\Lambda_+[h] \rangle_{\mathit{L}^{2}(\S)^4}=0.
\end{align}
Thus, $\mathcal{H}_{\k}$ is symmetric on $\mathrm{dom}(\mathcal{H}_{\k})$ and densely defined in  $\mathit{L}^{2}(\rr^3)^4$. This finishes the proof. \qed
\newline

   The following proposition gives a description of the domain of the adjoint operator $\mathcal{H}^{\ast}_{\k}$. 
   \begin{proposition}\label{adjoint1} Let $ \mathcal{H}_{\k}$ be as in the definition \ref{def1}. Then we have
   \begin{align}\label{adjoint}
  \mathrm{dom}(\mathcal{H}^{\ast}_{\k})=\left\{ u+\Phi[g]: u\in\mathit{H}^1(\rr^3)^4, g\in\mathit{H}^{-1/2}(\S)^4, t_{\S}u=-\tilde{\Lambda}_{+}[g]\right\}.
  \end{align}
   \end{proposition}
   \textbf{Proof.} Let $D$ be the set on the right-hand of \eqref{adjoint}. First, we prove the inclusion $D\subset  \mathrm{dom}(\mathcal{H}^{\ast}_{\k})$.  Given $\varphi:=v+\Phi[h]\in D$ and $\psi= u+\Phi[g]\in\mathrm{dom}(\mathcal{H}_{\k})$, then 
\begin{align*}\label{}
  \langle\varphi, \mathcal{H}_{\k}\psi\rangle_{\mathit{L}^2(\rr^3)^4}&=  \langle  \mathcal{H}v,u\rangle_{\mathit{L}^2(\rr^3)^4}+  \langle\Phi[h], \mathcal{H}u\rangle_{\mathit{L}^2(\rr^3)^4}=  \langle  \mathcal{H}v,u\rangle_{\mathit{L}^2(\rr^3)^4}+  \langle h,t_{\S}u\rangle_{\mathit{H}^{-1/2},\mathit{H}^{1/2}}\\
   &=  \langle  \mathcal{H}v,u\rangle_{\mathit{L}^2(\rr^3)^4}+  \langle h,-\Lambda_+[g]\rangle_{\mathit{H}^{-1/2},\mathit{H}^{1/2}}=  \langle  \mathcal{H}v,u\rangle_{\mathit{L}^2(\rr^3)^4}+  \langle t_{\S}v ,g\rangle_{\mathit{H}^{-1/2},\mathit{H}^{1/2}}\\
&=  \langle  \mathcal{H}v,\psi\rangle_{\mathit{L}^2(\rr^3)^4}.
\end{align*}
  Which yields  $\varphi\subset\mathrm{dom}(\mathcal{H}_{\k}^{\ast})$ and thus $D\subset\mathrm{dom}(\mathcal{H}_{\k}^{\ast})$. 
  
   Now we prove the inclusion $\mathrm{dom}(\mathcal{H}_{\k}^{\ast})\subset D$. Fix $\varphi \in\mathrm{dom}(\mathcal{H}_{\k}^{\ast})$, we first show that there exist functions $v\in\mathit{H}^1(\rr^3)^4$ and $ h\in\mathit{H}^{-1/2}(\S)^4$  uniquely determined by $\varphi$ such that $\varphi=v+\Phi[h]$. For that, let $\psi=(\psi_+,\psi_-)\in\mathcal{D}(\O_+)^4\oplus\mathcal{D}(\O_-)^4$, then by definition there is $U=(U_+,U_-)\in\mathit{L}^2(\rr^3)^4$  such that 
 \begin{align*}
   \langle\mathcal{H}\varphi, \psi\rangle_{\mathcal{D}^{\prime}(\rr^3)^4,\mathcal{D}(\rr^3)^4}&= \langle  \varphi,\mathcal{H}\psi\rangle_{\mathcal{D}^{\prime}(\rr^3)^4,\mathcal{D}(\rr^3)^4}= \langle \varphi_+, \mathcal{H}\psi_+\rangle_{\mathit{L}^{2}(\O_+)^4}+\langle \varphi_-, \mathcal{H}\psi_-\rangle_{\mathit{L}^{2}(\O_-)^4}\\
   &=\langle U_+, \psi_+\rangle_{\mathit{L}^{2}(\O_+)^4}+\langle U_-, \psi_-\rangle_{\mathit{L}^{2}(\O_-)^4}= \langle U, \psi\rangle_{\mathit{L}^{2}(\rr^3)^4}
\end{align*}
 Thus we obtain $ \mathcal{H} \varphi_\pm= U_\pm$ in $\mathcal{D}^{\prime}(\O_\pm)^4$ and then in $\mathit{L}^2(\O_\pm)^4$. From this we conclude that $\varphi\in\mathit{H}(\alpha,\O_+)\oplus \mathit{H}(\alpha,\O_-)$. Set 
\begin{align}\label{determinifunct}
 h=i(\alpha\cdot\mathit{N})(t_{\S}\varphi_+-t_{\S}\varphi_-)\, \text{ and }\, v=\varphi-\Phi[h].
\end{align}
As $t_{\S}\varphi_\pm\in\mathit{H}^{-1/2}(\S)^4$ holds by Proposition \ref{extension}, it follows that $h\in\mathit{H}^{-1/2}(\S)^4$ and $v\in\mathit{H}(\alpha,\O_+)\oplus \mathit{H}(\alpha,\O_-)$. Moreover,   a simple computation yields that
\begin{align*}
 t_{\S}(v\downharpoonright_{\Omega_\pm})=\left( \frac{1}{2} -i\tilde{\mathit{C}_{\S}}(\alpha\cdot\mathit{N})\right) t_{\S}\varphi_+ + \left( \frac{1}{2} +i\tilde{\mathit{C}_{\S}}(\alpha\cdot\mathit{N})\right) t_{\S}\varphi_-.
\end{align*}
Thanks to Proposition \ref{extension}-$(\mathrm{iv})$,  we know that $t_{\S}v\in \mathit{H}^{1/2}(\S)^4$ and  we have $v\in\mathit{H}^1(\rr^3)^4$, which justifies the decomposition $\varphi=v+\Phi[h]$.  Since $\varphi\in\mathrm{dom}(\mathcal{H}_{\k}^{\ast})\cap\mathit{H}(\alpha,\O_+)\oplus \mathit{H}(\alpha,\O_-)$, from \eqref{Geen formula} and \eqref{dual} it follows that 
 \begin{align}\label{equdua}
 \begin{split}
0&=\langle (-i\alpha\cdot\mathit{N})t_{\S}\varphi_+, t_{\S}\psi_+\rangle_{\mathit{H}^{-1/2},\mathit{H}^{1/2}}-\langle (-i\alpha\cdot\mathit{N})t_{\S}\varphi_-, t_{\S}\psi_-\rangle_{\mathit{H}^{-1/2},\mathit{H}^{1/2}}\\
&= \langle t_{\S}v, g\rangle_{\mathit{L}^{2}(\S)^4}-\langle h, t_{\S}u\rangle_{\mathit{H}^{-1/2},\mathit{H}^{1/2}},
\end{split}
\end{align}
for all $\psi=u+\Phi[g]\in\mathrm{dom}(\mathcal{H}_{\k})\cap\mathit{H}^{1}(\O_+)^4\oplus \mathit{H}^{1}(\O_-)^4$. 

Let $g\in\mathit{H}^{1/2}(\S)^4$ and set $u=E(-\Lambda_+[g])\in\mathit{H}^1(\rr^3)^4$, where $E$ is the extension operator. Then $u+\Phi[g]\in\mathrm{dom}(\mathcal{H}_{\k})$ and  by \eqref{equdua} we obtain
 \begin{align}\label{lastequa}
 \langle t_{\S}v, g\rangle_{\mathit{L}^{2}(\S)^4}-\langle h, t_{\S}u\rangle_{\mathit{H}^{-1/2},\mathit{H}^{1/2}}=\langle t_{\S}v+\tilde{\Lambda}_+[h], g\rangle_{\mathit{H}^{-1/2},\mathit{H}^{1/2}}=0,
\end{align}
where the condition $ t_{\S}u=-\Lambda_+[g]$ was used in the last step. As \eqref{lastequa} holds for all $g\in\mathit{H}^{1/2}(\S)^4$ we conclude that $ t_{\S}v=-\tilde{\Lambda}_{+}[h]$ holds in $\mathit{H}^{-1/2}(\S)^4$, and then in $\mathit{H}^{1/2}(\S)^4$, which proves the inclusion  $\mathrm{dom}(\mathcal{H}_{\k}^{\ast})\subset D$ and completes the proof of the proposition.\qed
\newline

Given $z\in\cc\setminus\left((-\infty,-m]\cup[m,\infty)\right)$, it is known that the fundamental solution of $(\Delta +m^2-z^2)I_4$  is given by 
\begin{align}
\label{}
 \psi^z(x)=\frac{e^{i\sqrt{z^2-m^2}|x|}}{4\pi|x|}\mathit{I}_4,\quad\text{ for } x\in\rr^3.
\end{align}
Moreover, the trace of the single-layer associated to $(\Delta +m^2-z^2)I_4$,  denoted by  $S^z$,  has the integral  representation 
\begin{align}\label{SL}
S^{z}[g](x)=\int_{\S} \psi^z(x-y)g(y)\mathrm{dS}(y), \quad \text{ for all }  x\in\S \text{ and } g\in\mathit{L}^{2}(\S)^4.
\end{align}
 If $z=0$, we simply write $S:=S^0$.  
   
  The next result contains  the main tools to prove the self-adjointness of the Dirac operator $\mathcal{H}_{\k}$. Recall that  $\lbrace A, B\rbrace= AB +BA$ is the usual anticommutator bracket.
\begin{lemma}\label{commutator}Given $a\in(-m,m)$, then the following hold:
\begin{itemize}
  \item[(i)] The anticommutator  $\lbrace \beta, \mathit{C}^{a}_{\S}\rbrace$ extends to a  bounded operator from $\mathit{H}^{-1/2}(\S)^4$ onto $\mathit{H}^{1/2}(\S)^4$. In particular, if $\S$ satisfies $\mathrm{(H1)}$, then $\lbrace \beta, \mathit{C}^{a}_{\S}\rbrace$ is a compact operator in $\mathit{L}^2(\S)^4$.
  \item[(ii)] The anticommutator  $\lbrace \alpha\cdot \mathit{N}, \mathit{C}^{a}_{\S}\rbrace$ extends to a  bounded operator from $\mathit{H}^{-1/2}(\S)^4$ to $\mathit{H}^{1/2}(\S)^4$. In particular,   if $\S$ satisfies $\mathrm{(H1)}$, then $\lbrace \alpha\cdot \mathit{N}, \mathit{C}^{a}_{\S}\rbrace$ is a compact operator in $\mathit{L}^2(\S)^4$.
   \item[(iii)] If $\S$ satisfies $\mathrm{(H2)}$, then  $\lbrace \alpha\cdot \mathit{N}, \mathit{C}^{}_{\S}\rbrace$ is a compact operator in $\mathit{L}^2(\S)^4$.
\end{itemize}
\end{lemma}

\textbf{Proof.} We are going to prove item $(\mathrm{i})$.  For this, observe that 
 \begin{align}\label{psi}
\frac{1}{2(m^2-a^2)}(m \mathit{I}_4-a\beta) \{ \beta,\mathit{C}^{a}_{\S}\}[g](x)=S^{a}[g](x).
\end{align}
 Hence,  the first statement of $(\mathrm{i})$ follows by  \cite[Theorem 6.11]{Mc} (see also \cite{Me}) for example. Furthermore, if $\S$ satisfies $(\mathrm{H1})$, then using that  the embedding $\mathit{H}^{1/2}(\S)^4\hookrightarrow\mathit{L}^{2}(\S)^4$ is compact, we then get that $ \{ \beta,\mathit{C}^{a}_{\S}\}$ is a compact operator in $\mathit{L}^{2}(\S)^4$. This finishes the proof of  $(\mathrm{i})$.

Now we prove  item $(\mathrm{ii})$. Let $x\in\S$ and $y\in\rr^3$,  a straightforward computation using the anticommutation relations of the Dirac matrices  yields that
\begin{align}\label{clifford}
  (\alpha\cdot \mathit{N}(x))(\alpha\cdot y)= -(\alpha\cdot y) (\alpha\cdot \mathit{N}(x)) +2(\mathit{N}(x)\cdot y)\mathit{I}_4.
\end{align}
Using \eqref{clifford} it follows that
\begin{align*}
(\alpha\cdot \mathit{N}(x))\phi^{a}(y) =& -\phi^{a}(y)(\alpha\cdot \mathit{N}(x)) -\frac{e^{-\sqrt{m^2-a^2}|y|}}{2i\pi |y|^3}(1+m|y|)(\mathit{N}(x)\cdot y)\mathit{I}_4+ 2a(\alpha\cdot \mathit{N}(x))\psi^a(y).
\end{align*}
Note that there are constants $C_1$ and $C_2$ such that, for all $x,y\in\S$, it holds that 
\begin{align*}
|\mathit{N}(x)-\mathit{N}(y)|\leqslant C_1 |x-y| \quad\text{and }\quad  |\mathit{N}(x)\cdot(x-y)|\leqslant C_2|x-y|^2,
\end{align*}
 see \cite[Lemma 3.15]{Fo} for example. Using this, for $g\in\mathit{L}^2(\S)^4$, we get that 
  \begin{align}\label{dec}
  \begin{split}
    \lbrace \alpha\cdot \mathit{N}, \mathit{C}^{a}_{\S}\rbrace[g](x)&= \int_{\S}K_a(x,y)g(y)dS(y)+ 2a(\alpha\cdot \mathit{N}(x))S^{a}[g](x)\\
    &:= T_{a,1}[g](x) + T_{a,2}[g](x),
    \end{split}
    \end{align}
where the kernel $K_a$ is given by 
\begin{align*}
 K_a(x,y)&= \phi^{a}(x-y)(\alpha\cdot(\mathit{N}(y)-\mathit{N}(x)) -\frac{e^{-\sqrt{m^2-a^2}|x-y|}}{2i\pi |x-y|^3}(1+\sqrt{m^2-a^2}|x-y|)(\mathit{N}(x)\cdot(x-y))\mathit{I}_4.
 \end{align*}
 Since $\S$ is $\mathit{C}^2$-smooth, from $(\mathrm{i})$ it follows  immediately that $T_{a,2}$ is bounded from $\mathit{H}^{-1/2}(\S)^4$ to $\mathit{H}^{1/2}(\S)^4$. Hence, it remains to prove that $T_{a,1}$ is bounded  from $\mathit{H}^{-1/2}(\S)^4$ to $\mathit{H}^{1/2}(\S)^4$. Actually,  if $\S$ satisfies $(\mathrm{H1})$, then the result follows with the same arguments as \cite[ Proposition 2.8]{OBV}, where the authors proved the statement for $a=0$. 
 
 Now, assume that  $\S$ satisfies $(\mathrm{H2})$ and recall that $F$ denotes  the flat part of $\S=\S_\nu$. Remark that  
 \begin{align*}
 \mathit{N}(y)-\mathit{N}(x)=0=\mathit{N}(x)\cdot(x-y)\quad \text{ if } x,y\in F.
 \end{align*} 
Therefore  the kernel  $K_a(x,y)$ vanishes for all $x,y\in F$ and it holds that $|K_a(x,y)|\leqslant C |x-y|^{-1} $. 
 Let $\chi:\S\rightarrow \rr$ be a $\mathit{C}^{\infty}$-smooth  and compactly  supported function on $\S$ and such that  $\S_\nu\setminus F\subsetneq \mathrm{supp}(\chi)$. Using that $K_a(x,y)$ vanishes for all $x,y\in F$,  we then obtain  
  \begin{align}\label{chi decomposition}
 T_{a,1}= \chi T_{a,1}\chi  + \chi T_{a,1} (1-\chi) +(1-\chi) T_{a,1}\chi:= T_1+T_2+T_3.
 \end{align}
 Again, one can extend  $T_1$  to a bounded operators from $\mathit{H}^{-1/2}(\S)^4$ to $\mathit{H}^{1/2}(\S)^4$ as much the same way as in the case of assumption $(\mathrm{H1})$.  
Now, we show that  $T_2$ is bounded from $\mathit{H}^{-1/2}(\S)^4$ into $\mathit{H}^{1/2}(\S)^4$, the proof for $ T_3$ is similar.    Note  that  $T_2$ is not singular  and there is $\delta>0$ such that 
 \begin{align*}
   T_{2}[g](x)= \int_{\substack{ |x-y|>\delta\\y\in\S\setminus \mathrm{supp}(\chi)}}\chi(x) K_a(x,y)(1-\chi(y))g(y)dS(y),  \quad g\in\mathit{L}^2(\S)^4. 
    \end{align*}
Let $V\subset\rr^2$ and  $\phi:V\rightarrow \rr^2$ be a linear affine  function which parametrizes $\mathrm{supp}(\chi)$. Clearly, the mapping $V\ni\overline{x}\rightarrow K_a(\phi(\overline{x}),y)$  is $\mathit{C}^{1}$-smooth for all $y\in \S$, and  the mapping $F\ni\overline{x}\rightarrow K_a(\phi(\overline{x}),y)$  is $\mathit{C}^{\infty}$-smooth for all $\overline{x}\in V$. Using this, it follows that $(T_2[g])\circ\phi$ is differentiable on $V$ and we have 
\begin{align*}
   \partial_{\overline{x}}(T_{2}[g])(\phi (\overline{x}))=&\int_{\substack{ |\phi(\overline{x})-y|>\delta\\ y\in\S\setminus \mathrm{supp}(\chi)}}  (\partial_{\overline{x}} \chi(\phi(\overline{x})) K_a(x,y)(1-\chi(y))g(y)dS(y)\\
   &+   \int_{\substack{ |\phi(\overline{x})-y|>\delta\\ y\in\S\setminus \mathrm{supp}(\chi)}}  \chi(\phi(\overline{x})) \partial_{\overline{x}}K(\phi (\overline{x}),y)(1-\chi(y))g(y)dS(y)\\
   &:=(T_{2,1}[g])(\phi (\overline{x}))+  (T_{2,2}[g])(\phi (\overline{x})). 
    \end{align*}
Using the boundedness of $(\partial_{\overline{x}} \chi(\phi(\cdot))$, it is easy to see that $T_{2,1}$ is bounded from $\mathit{L}^2(\S)^4$ into itself. Now, define the kernels
\begin{align*}
K_1(\phi (\overline{x}),y)&=    |\phi (\overline{x})-y|^{\frac{-1}{2}}(1-\chi(y))\text{ and }K_2(\phi (\overline{x}),y)&=  \chi(\phi(\overline{x})) |\phi (\overline{x})-y|^{\frac{1}{2}}\partial_{\overline{x}}K_a(\phi (\overline{x}),y),
 \end{align*}
 and note that $ |\partial_{\overline{x}}K_a(\phi (\overline{x}),y)|\leqslant C|\phi (\overline{x})-y|^{-2}$. Therefore we obtain 
 \begin{align*}
 \sup_{\overline{x}\in V}\int_\S |K_{1}(\phi (\overline{x}),y)|^2dS(y)<\infty \quad \text{and}\quad \sup_{y\in\S}\int_V |K_{2}(\phi (\overline{x}),y)|^2d\overline{x} <\infty.
  \end{align*}
 Hence, using the Schur test (see \cite[Lemma 0.32]{Te} for example), it follows that 
   \begin{align*}
\left|\left| T_{2,2} \right|\right|^{2}_{\mathit{L}^2(\S)^4} \leqslant \left( \sup_{\overline{x}\in V}\int_\S |K_{1}(\phi (\overline{x}),y)|^2dS(y) \right)\left( \sup_{y\in\S}\int_V |K_{2}(\phi (\overline{x}),y)|^2d\overline{x}\right)<\infty.
 \end{align*} 
Thus $T_{2,2}$ is bounded from $\mathit{L}^2(\S)^4$ into itself, we then  conclude  that  $T_{2}$ is bounded from $\mathit{L}^2(\S)^4$ to $\mathit{H}^{1}(\S)^4$ and  by duality and interpolation arguments one can extend it continuously to a bounded operator from $\mathit{H}^{-1/2}(\S)^4$ to $\mathit{H}^{1/2}(\S)^4$. Therefore   $\lbrace \alpha\cdot \mathit{N}, \mathit{C}^{a}_{\S}\rbrace$ extends to a  bounded operator from $\mathit{H}^{-1/2}(\S)^4$ to $\mathit{H}^{1/2}(\S)^4$. The second statement is a direct consequence of the Sobolev injection, and this completes the proof of $(\mathrm{ii})$.

Now we turn to the proof of $(\mathrm{iii})$.  Assume that  $\S$ satisfies $(\mathrm{H2})$, then from $(\mathrm{ii})$ we know that $ \lbrace \alpha\cdot \mathit{N}, \mathit{C}_{\S}\rbrace$  coincides with $T_{a,1}$  for $a=0$, and it is bounded from $\mathit{H}^{-1/2}(\S)^4$ to $\mathit{H}^{1/2}(\S)^4$.   Hence, $ \lbrace \alpha\cdot \mathit{N}, \mathit{C}_{\S}\rbrace$ is compact on $\mathit{L}^2(\S)^4$ by the decomposition \eqref{chi decomposition} and the compactness of the  Sobolev embedding $\chi\mathit{H}^{1/2}(\S)^4\hookrightarrow\mathit{L}^{2}(\S)^4$. This finishes the proof of the lemma.\qed
 \newline
\begin{remark}\label{ext remark} Actually the above result is not surprising since the integral kernels associated to the anticommutators $\lbrace \alpha\cdot \mathit{N}, \mathit{C}^{a}_{\S}\rbrace$ and $\lbrace \beta, \mathit{C}^{a}_{\S}\rbrace$ behave locally like $|x-y|^{-1}$, when $|x-y|$ tends to zero. Therefore, the operators in consideration are bounded from $\mathit{L}^2(\S)^4$ to $\mathit{H}^1(\S)^4$ because $\S$ is $\mathit{C}^2$-smooth. Moreover, the same result holds true for all $a\in\cc\setminus\left((-\infty,-m]\cup[m,\infty)\right)$. 
\end{remark}
We are now in position to prove Theorem \ref{main1}.

  \textbf{ Proof of Theorem \ref{main1}}   $(\mathrm{i})$ Let $\k\in\rr^3$ and recall that   $\mathrm{sgn}(\k)$ is defined by  \eqref{sgndef}. Assume that $\mathrm{sgn}(\k)\neq 4$. From the definition of $\tilde{\Lambda}_{\pm}^{a}$,  a simple computation using Lemma \ref{lemme 2.1}$(\mathrm{ii})$ gives
\begin{align}
\begin{split}
\tilde{\Lambda}_{\pm}^{a}\tilde{\Lambda}_{\mp}^{a}&=\frac{1}{\mathrm{sgn}(\k)}-(\tilde{\mathit{C}}^{a}_{\S})^2+\frac{\mu}{\mathrm{sng}(\k)}\{\beta,\tilde{\mathit{C}}^{a}_{\S}\} + \frac{\e}{\mathrm{sgn}(\k)}\lbrace \alpha\cdot \mathit{N}, \tilde{\mathit{C}}^{a}_{\S}\rbrace\\
&=\frac{1}{\mathrm{sgn}(\k)}-\frac{1}{4}-\mathit{C}^{a}_{\S} (\alpha\cdot \mathit{N})\lbrace \alpha\cdot \mathit{N}, \tilde{\mathit{C}}^{a}_{\S}\rbrace +\frac{\mu}{\mathrm{sgn}(\k)}\{\beta,\tilde{\mathit{C}}^{a}_{\S}\}+\frac{\e}{\mathrm{sgn}(\k)} \lbrace \alpha\cdot \mathit{N}, \tilde{\mathit{C}}^{a}_{\S}\rbrace\label{multi} .
\end{split}
\end{align}
  Let $g\in\mathit{H}^{-1/2}(\S)^4$ be such that  $\tilde{\Lambda}_{+}[g]\in\mathit{H}^{1/2}(\S)^4$. Then, it follows from \eqref{multi} that 
\begin{align*}
g=\frac{4(\mathrm{sgn}(\k))}{4-\mathrm{sgn}(\k)}\left( \Lambda_{-}^{}\tilde{\Lambda}_{+} +\mathit{C}^{a}_{\S} (\alpha\cdot \mathit{N})\lbrace \alpha\cdot \mathit{N}, \tilde{\mathit{C}}^{a}_{\S}\rbrace -\frac{\mu}{\mathrm{sgn}(\k)}\{\beta,\tilde{\mathit{C}}^{a}_{\S}\} -\frac{\e}{\mathrm{sgn}(\k)} \lbrace \alpha\cdot \mathit{N}, \tilde{\mathit{C}}^{a}_{\S}\rbrace\right)[g].
\end{align*}
Therefore, using Lemma \ref{commutator} we obtain  that $g\in\mathit{H}^{1/2}(\S)^4$. Hence,  given any $\varphi = u +\Phi[g]\in  \mathrm{dom}(\mathcal{H}^{\ast}_{\k})$, since  $ g\in\mathit{H}^{-1/2}(\S)^4$ and  $t_{\S}u=\tilde{\Lambda}_{+}[g]\in\mathit{H}^{1/2}(\S)^4$, we deduce that $ g\in\mathit{H}^{1/2}(\S)^4$. Thus, $ \mathrm{dom}(\mathcal{H}^{\ast}_{\k})=\mathrm{dom}(\mathcal{H}_{\k})$ and it holds that  
  \begin{align}
  \mathrm{dom}(\mathcal{H}_{\k})=\left\{ u+\Phi[g]: u\in\mathit{H}^1(\rr^3)^4, g\in\mathit{H}^{1/2}(\S)^4, t_{\S}u=-\Lambda_+[g]\right\}.
  \end{align}
   This finishes the proof of $(\mathrm{i})$.
   
   $(\mathrm{ii})$ Fix $\k$ such that $\mathrm{sgn}(\k)= 4$. Since $\mathcal{H}_{\k}$ is closable by Proposition \ref{closable}, it follows that  $\overline{\mathcal{H}_{\k}}\subset \mathcal{H}^{\ast}_{\k}$. Let us prove the other inclusion, for this given $\varphi=u+\Phi[g]\in  \mathrm{dom}(\mathcal{H}^{\ast}_{\k}) $ and let $(h_j)_{j\in\mathbb{N}}\subset \mathit{H}^{1/2}(\S)^4$ be a sequence of functions that converges to $g$ in $\mathit{H}^{-1/2}(\S)^4$. Set 
   \begin{align}\label{suite}
   g_j:= g +\frac{2}{\ep} \tilde{\Lambda}_{-}[h_j - g],\quad\forall j\in\mathbb{N}.
   \end{align}
   Then $(g_j)_{j\in\mathbb{N}}, (\Lambda_{+}[g_j])_{j\in\mathbb{N}}\subset \mathit{H}^{1/2}(\S)^4$, and it holds that 
   \begin{align}\label{cv}
g_j \xrightarrow[j\to\infty]{} g \text{ in } \mathit{H}^{-1/2}(\S)^4, \quad \Lambda_{+}[g_j] \xrightarrow[j\to\infty]{} \tilde{\Lambda}_{+}[g] ,  \text{ in } \mathit{H}^{1/2}(\S)^4.
\end{align}
Indeed, remark that $  \tilde{\Lambda}_{+}+ \tilde{\Lambda}_{-}=\ep/2$, thus  one can write $g_j$ as follows 
\begin{align*}
g_j=  \frac{2}{\ep}( \tilde{\Lambda}_{+}[g]+ \tilde{\Lambda}_{-}[h_j]).
\end{align*}
 Using this, \eqref{cv} follows easily since $\tilde{\Lambda}_{\pm}\tilde{\Lambda}_{\mp}$ are bounded from $\mathit{H}^{-1/2}(\S)^4$ to $\mathit{H}^{1/2}(\S)^4$ by Lemma \ref{commutator} and \eqref{multi}. Now, for $j\in\mathbb{N}$, we define $\varphi_j:= u_j +\Phi[g_j]$ where 
 \begin{align}
 u_j=u-v_j  \quad \text{and }\,v_j=E\left(\frac{2}{\ep}\tilde{\Lambda}_{+}\tilde{\Lambda}_{-}[h_j-g]\right) .
  \end{align}
 Clearly, we have $u_j\in\mathit{H}^1(\rr^3)^4 $ and $t_\S u_j= -\Lambda_+[g_j]\in \mathit{H}^{1/2}(\S)^4$, hence $(\varphi_j)_{j\in\mathbb{N}}\subset \mathrm{dom}(\mathcal{H}_{\k})$. Moreover, since $(h_j)_{j\in\mathbb{N}}$ (respectively $(g_j)_{j\in\mathbb{N}}$)  converges to $g$ in $\mathit{H}^{-1/2}(\S)^4$ as $j\longrightarrow \infty$, using the continuity of $\tilde{\Lambda}_{\pm}\tilde{\Lambda}_{\mp}$ it follows that  $(\varphi_j,\mathcal{H}_{\k}\varphi_j) \xrightarrow[j\to\infty]{}(\varphi, \mathcal{H}^{\ast}_{\k}\varphi)$ in $\mathit{L}^2(\rr^3)^4$. Therefore $ \mathcal{H}^{\ast}_{\k}\subset\overline{\mathcal{H}_{\k}}$ and the Theorem is proved. \qed
\newline 
\begin{remark}\label{appelrem} It is worthwhile to mention that, in view of \eqref{determinifunct},  the functions $u$ and $g$ in $\varphi= u+ \Phi[g]\in\mathrm{dom}(\overline{\mathcal{H}_{\k}})$ are uniquely determined by $\varphi$. Moreover, by  Proposition \ref{extension}-$(\mathrm{iv})$ we have that  $(\Phi^z-\Phi)[g]\in\mathit{H}^1(\rr^3)^4$. Consequently, for any $z\in\rho(\overline{\mathcal{H}_{\k}})\cap\rho(\mathcal{H})$ and $\varphi= u+ \Phi[g]\in\mathrm{dom}(\overline{\mathcal{H}_{\k}})$, there exist uniquely determined functions $v\in\mathit{H}^1(\rr^3)^4$ and $g\in\mathit{H}^{1/2}(\S)^4$ (respectively $g\in\mathit{H}^{-1/2}(\S)^4$ when $\mathrm{sgn}(\kappa)=4$) such that  $\varphi= v+ \Phi^z[g]$ and $(\overline{\mathcal{H}_{\k}}-z)\varphi=(\mathcal{H}-z)v$ (just write $\varphi= u -(\Phi^z-\Phi)[g]+\Phi^z[g]$). 
\end{remark}

  In the following, we explain how to define the Dirac operator $\overline{\mathcal{H}_{\k}}$ via transmission condition.   Let $\varphi= u+ \Phi[g]\in\mathrm{dom}(\overline{\mathcal{H}_{\k}})$ and set $\varphi_\pm := \varphi\downharpoonright_{\Omega_\pm}$. It is clear that  $\varphi_\pm, (\alpha\cdot\nabla) \varphi_\pm\in\mathit{L}^2(\Omega_{\pm})^4$. Now, we define $\delta_{\S}\varphi$ as the distribution
\begin{align}
\langle \delta_{\S}\varphi,\psi\rangle_{\mathcal{D}^{\prime}(\rr^3)^4,\mathcal{D}(\rr^3)^4}:= \frac{1}{2}\int_{\S}\langle t_{\S}\varphi_+ +t_{\S}\varphi_-,\psi\rangle_{\cc^4}\mathrm{dS}(x),\quad \text{ for all } \psi\in\mathcal{D}(\rr^3)^4.
\end{align}
Therefore, a simple computation in the sense of distributions yields
\begin{align*}
 (\mathcal{H}+(\ep I_4+\mu\beta + \eta(\alpha\cdot N))&\delta_{\S})\varphi= (-i\alpha\cdot\nabla +m\beta)\varphi +\frac{1}{2}(\ep I_4+\mu\beta+ \eta(\alpha\cdot N))(t_{\S}\varphi_+ +t_{\S}\varphi_-)\delta_{\S},\\
=& (-i\alpha\cdot\nabla +m\beta)\varphi_+\oplus(-i\alpha\cdot\nabla +m\beta)\varphi_- +i\alpha\cdot\mathit{N}(t_{\S}\varphi_+ -t_{\S}\varphi_-)\delta_{\S}\\
&+\frac{1}{2}(\ep I_4+\mu\beta+ \eta(\alpha\cdot N))(t_{\S}\varphi_+ +t_{\S}\varphi_-)\delta_{\S} .
\end{align*}
Using the Plemelj-Sokhotski formula (see Lemma \ref{lemme 2.1}), a computation shows that
\begin{align}\label{trans}
\frac{1}{2}(\ep I_4+\mu\beta+ \eta(\alpha\cdot N))(t_{\S}\varphi_+ +t_{\S}\varphi_-)\delta_{\S} +i\alpha\cdot\mathit{N}(t_{\S}\varphi_+ -t_{\S}\varphi_-)\delta_{\S}=0,
\end{align}
holds in $\mathit{H}^{-1/2}(\S)^4$. Since $(-i\alpha\cdot\nabla +m\beta)\varphi_+\oplus(-i\alpha\cdot\nabla +m\beta)\varphi_- \in \mathit{L}^2(\rr^3)^4$, given $\varphi= (\varphi_+,\varphi_-)\in\mathit{L}^2(\rr^3)^4$ such that $ (\alpha\cdot\nabla) \varphi_\pm\in\mathit{L}^2(\Omega_{\pm})^4$ and satisfying \eqref{trans}, it holds that $\mathcal{H}^{}_{\k}\varphi \in \mathit{L}^2(\rr^3)^4$. In particular, this leads to the following definition:
\begin{definition}\label{def2} Let $\k=(\ep,\mu,\e)\in\rr^3$ be such that $\mathrm{sgn}(\k)\neq0$ and $m>0$. The self-adjoint Dirac operator coupled with a combination of electrostatic, Lorentz scalar and normal vector field $\delta$-shell interactions of strength $\ep$, $\m$ and $\e$ respectively,  is the operator $\overline{\mathcal{H}_{\k}}$ defined on the domain
 \begin{align}\label{dom1}
 \begin{split}
\mathrm{dom}(\overline{\mathcal{H}_{\k}})=\big\{ \varphi=(\varphi_+,\varphi_-)&\in\mathit{L}^2(\Omega_+)^4\oplus\mathit{L}^2(\Omega_-)^4 :\\
& (\alpha\cdot\nabla) \varphi_\pm\in\mathit{L}^2(\Omega_{\pm})^4 \text{ and }\eqref{trans} \text{ holds in } \mathit{H}^{-1/2}(\S)^4 \big\},
\end{split}
\end{align}
 and it acts in the sense of distributions as $ \overline{\mathcal{H}_{\k}}(\varphi)= (\mathcal{H}\varphi_+)\oplus(\mathcal{H}\varphi_-)$, for all $\varphi\in\mathrm{dom}(\overline{\mathcal{H}_{\k}})$. 
\end{definition}
\begin{remark}\label{remarkregu} Assume that $\mathrm{sgn}(\k)\neq 0,4$. Since the operator $\Phi$  is bounded  from $\mathit{H}^{1/2}(\S)^4$ to $\mathit{H}^{1}(\rr^3\setminus\S)^4$, it holds that $\varphi_\pm := \varphi|_{\Omega_\pm}\in\mathit{H}^1(\Omega_\pm)^4$. Moreover, following the same arguments as above, we conclude that the transmission condition  \eqref{trans} holds actually in $\mathit{H}^{1/2}(\S)^4$. Therefore, it follows that
 \begin{align}\label{dom2}
\mathrm{dom}(\mathcal{H}_{\k})=\left\{ \varphi=(\varphi_+,\varphi_-)\in\mathit{H}^1(\Omega_+)^4\oplus\mathit{H}^1(\Omega_-)^4 :  \eqref{trans} \text{ holds in } \mathit{H}^{1/2}(\S)^4 \right \}.
\end{align}
\end{remark}

  Let us make some comments on the technique developed here. Note that the condition on $\S$ of being $\mathit{C}^2$-smooth is minimal to prove the self-adjointness of $\overline{\mathcal{H}_{\k}}$ when $\mathrm{sgn}(\k)=4$. Indeed, the main ingredient that we have used is the continuity  of $\Lambda_{\pm}\Lambda_{\mp}$  from $\mathit{H}^{-1/2}(\S)^4$ to $\mathit{H}^{1/2}(\S)^4$, or equivalently, the continuity of the anticommutators $\{\beta,\mathit{C}_{\S}\}$ and $\lbrace \alpha\cdot \mathit{N}, \mathit{C}_{\S}\rbrace$. Since $\{\beta,\mathit{C}_{\S}\}$ involves the trace of the single-layer potential, we can always extend it to a bounded operator from $\mathit{H}^{-1/2}(\S)$ to $\mathit{H}^{1/2}(\S)$, even if $\S$ is Lipschitz. However,  $\lbrace \alpha\cdot \mathit{N}, \mathit{C}_{\S}\rbrace$ involves the principal value of the double-layer potential, its adjoint and  the commutators $[\mathit{N}_k, R_j]$, where $R_j$ are the Riesz transforms (see \cite[Lemma 4.1]{BB}), and it is well known that the $\mathit{C}^2$ regularity is minimal to extend continuously  these operators from $\mathit{H}^{-1/2}(\S)$ to $\mathit{H}^{1/2}(\S)$. However, if $\mathrm{sgn}(\k)\neq0,4$ and $\O_+$ is a bounded $\mathit{C}^{1,\g}$-smooth domain, for some $\g\in(1/2,1)$, then one can manage to prove the self-adjointness of $\mathcal{H}_{\k}$ using the technique developed in this part,  see  \cite[Section 5]{BB} for more details.
\subsection{On the Dirac Operator with Electrostatic and Lorentz scalar $\delta$-Shell interactions }\label{sub3.1}
 We discuss in this part the self-adjointness of the Dirac operator  $\mathcal{H}_{\k}$ in the case $\e=0$, and we denote it by  $\mathcal{H}_{\ep,\m}$. This operator is well known as  the  Dirac operator with electrostatic and Lorentz scalar $\delta$-shell interactions, cf. \cite{AMV2},\cite{BEHL2},\cite{BHOP}.
  If $|\ep | \neq|\m|$, from Theorem \ref{main1} we get immediately the following result.
  \begin{proposition}\label{DesL} Given $\ep,\m\in\rr\setminus \{ 0 \}$ such that $| \ep | \neq|\m|$, and  define the operators $\Lambda_\pm$ as follows
  \begin{align}
\Lambda^{}_{\pm}=\frac{1}{\ep^2-\mu^2}(\ep I_4 \mp\mu\beta)\pm\mathit{C}^{}_{\S}.
\end{align}
Then, the following hold:
\begin{itemize}
 \item[(i)] If  $\ep^2-\mu^2\neq4$, then $ \mathcal{H}_{\ep,\mu}$ is self-adjoint and we have 
  \begin{align}
  \mathrm{dom}(\mathcal{H}_{\ep,\mu})=\left\{ u+\Phi[g]: u\in\mathit{H}^1(\rr^3)^4, g\in\mathit{H}^{1/2}(\S)^4, t_{\S}u=-\Lambda_+[g]\right\}.
  \end{align}
  \item[(ii)]  If  $\ep^2-\mu^2=4$, then $\mathcal{H}_{\ep,\mu}$ is essentially self-adjoint and we have 
  \begin{align*}
  \mathrm{dom}(\overline{\mathcal{H}_{\ep,\mu}})= \left\{ u+\Phi[g]: u\in\mathit{H}^1(\rr^3)^4, g\in\mathit{H}^{-1/2}(\S)^4, t_{\S}u=-\tilde{\Lambda}_+[g]\right\}.
  \end{align*}
 \end{itemize}
\end{proposition}

  Now we turn to the special case $\ep=\pm\m$.  Set $P_\pm= (\mathit{I}_4\pm\beta)/2$, then $\mathcal{H}_{\ep,\pm\ep}$ is given formally by  
\begin{align}
\mathcal{H}_{\ep,\pm\ep} =  \mathcal{H} + P_\pm V_{\ep,\pm\ep}=- i \alpha \cdot\nabla + m\beta+2\ep  P_\pm\delta_{\S}.
\end{align} 
Define 
\begin{align}\label{lambda2}
\Lambda_{+}=P_\pm\left(1/2\ep+\mathit{C}_{\S}\right)P_\pm \quad \text{and } \Lambda_{-}=P_\pm\left(1/2\ep -\mathit{C}_{\S}\right)P_\pm, \quad\forall \ep \neq0.
\end{align}
Clearly , $\Lambda_{\pm}$ are bounded and self-adjoint from $P_\pm\mathit{L}^{2}(\S)^4$ onto itself (respectively from $P_\pm\mathit{H}^{1/2}(\S)^4$ onto itself). In order to define rigorously $\mathcal{H}_{\ep,\pm\ep}$  as in Definition \ref{def1}, that is  $\mathcal{H}_{\ep,\pm\ep}\varphi=\mathcal{H}u$ in the sense of distributions for $ \varphi = u+\Phi[g]$, with $u\in\mathit{H}^1(\rr^3)^4$ and $ g\in\mathit{H}^{1/2}(\S)^4$. We shall  take $ g\in P_\pm\mathit{H}^{1/2}(\S)^4$ and assume the condition  $  P_\pm t_{\S}u=-P_\pm \Lambda_{+}[g]$.  Indeed, if we set 
\begin{align}
\mathrm{dom}(\mathcal{H}_{\ep,\pm\ep})=\left\{ u+\Phi[g]: u\in\mathit{H}^1(\rr^3)^4, g\in P_\pm\mathit{L}^{2}(\S)^4 \text{ and }P_\pm t_{\S}u=-P_\pm \Lambda_{+}[g]\right\},
\end{align}
Then, in a similar way as in Proposition \ref{closable} and Proposition \ref{adjoint1}, one can check  that $(\mathcal{H}_{\ep,\pm\ep}, \mathrm{dom}(\mathcal{H}_{\ep,\pm\ep}))$ is closable and its adjoint is defined on the domain
   \begin{align}
  \mathrm{dom}(\mathcal{H}^{\ast}_{\ep,\pm\ep})=\left\{ u+\Phi[g]: u\in\mathit{H}^1(\rr^3)^4, g\in P_\pm\mathit{H}^{-1/2}(\S)^4,P_\pm t_{\S}u=-P_\pm\tilde{\Lambda}_{+}[g]\right\},
  \end{align}
where $\tilde{\Lambda}_\pm$ denotes the bounded extension of $\Lambda_\pm $ from $P_\pm\mathit{H}^{-1/2}(\S)^4$ onto itself, and we get the analogous  of Theorem \ref{main1}  in this case which reads as follows:
  \begin{proposition}\label{lametmu} Assume that  $\ep\neq0$, then $(\mathcal{H}_{\ep,\pm\ep}, \mathrm{dom}(\mathcal{H}_{\ep,\pm\ep}))$ is self-adjoint and we have 
    \begin{align}
  \mathrm{dom}(\mathcal{H}_{\ep,\pm\ep})=\left\{ u+\Phi[g]: u\in\mathit{H}^1(\rr^3)^4, g\in P_\pm\mathit{H}^{1/2}(\S)^4, P_\pm t_{\S}u=-P_\pm\Lambda_+[g]\right\}.
  \end{align}
  \end{proposition} 
  \textbf{Proof.}  Fix $\ep\neq0$ and let $\tilde{\Lambda}_\pm$ be as above. Using the relations $P_\pm \alpha_j=P_\mp\alpha_j$ and $P_\pm\beta=\beta P_\pm$, a simple computation yields
 \begin{align}\label{remote}
\tilde{\Lambda}_{-}\tilde{\Lambda}_{+}&= \frac{1}{4\ep^2}P_\pm -\frac{1}{2\ep}P_\pm\tilde{\mathit{C}_{\S}}P_\pm\tilde{\mathit{C}_{\S}}P_\pm= \frac{1}{4\ep^2}P_\pm -\frac{m^2}{2\ep}(S)^2P_\pm,
\end{align}
 where $S$ is given by \eqref{SL}. Recall that $\tilde{\Lambda}_{-}\tilde{\Lambda}_{+}$ and  $S P_\pm$ are bounded from $P_\pm\mathit{H}^{-1/2}(\S)^4$ into $P_\pm\mathit{H}^{1/2}(\S)^4$. Thus, from \eqref{remote} it follows that, if $g\in P_\pm\mathit{H}^{-1/2}(\S)^4$ and  $\tilde{\Lambda}_{+}[g]\in P_\pm\mathit{H}^{1/2}(\S)^4$, then  $g\in P_\pm\mathit{H}^{1/2}(\S)^4$. Which yields that  $\mathrm{dom}(\mathcal{H}_{\ep,\pm\ep})= \mathrm{dom}(\mathcal{H}^{\ast}_{\ep,\pm\ep})$ and the proposition is proved.\qed

\subsection{The operators $\Lambda^a_\pm$}\label{sub3.2}
Let $a\in(-m,m)$ and let $\Lambda^a_\pm$ be as in the Notation \ref{Lambda}. From the proof of Theorem \ref{main1}, it is evident that the study of the self-adjointness character of $\mathcal{H}_{\k}$ is related to the spectral properties of $\Lambda_+$. The goal of this part is to establish the connection between $\mathcal{H}_{\k}$ and $\Lambda_+$. For this, we introduce the Laplace-Beltrami operators $\Delta_\S$ on $\S$ and we define the operator $L:=(c-\Delta_\S)\mathit{I}_4$ (here we assume that $c$ is big enough if $\S$ satisfies $(\mathrm{H2})$, so that $c$ is not in the spectrum of $\Delta_\S$). It is well known that $L^{\pm 1/4}$ is a bijective operator from $\mathit{H}^{\pm1/2}(\S)^4$ onto $\mathit{L}^2(\S)^4$. Hence, one can  write the domain of  $\mathcal{H}_{\k}$ as follows:
 \begin{align}\label{dom22}
\mathrm{dom}(\mathcal{H}_{\k})=\left\{ u+\Phi L^{1/4}[g]: u\in\mathit{H}^1(\rr^3)^4, g\in\mathit{H}^{1/2}(\S)^4 \text{ and }L^{1/4} t_{\S}u=-L^{1/4}\Lambda_{+}L^{1/4}[g]\right\},
\end{align}
which leads us to define the following unbounded operators
\begin{align}\label{dom3}
\mathcal{L}^{a}_{\pm}:=L^{1/4}\Lambda^a_{\pm}L^{1/4}\quad\text{ with }\quad\mathrm{dom}(\mathcal{L}^{a}_{\pm})=\left\{ g\in\mathit{H}^{1/2}(\S)^4 : \Lambda_{\pm}^aL^{1/4}[g]\in\mathit{H}^{1/2}(\S)^4  \right\}.
\end{align}
In the following lemma, we study the self-adjointness character of $\mathcal{L}^{a}_{\pm}$, which will clarify the relationship between $\mathcal{H}_{\k}$ and $\Lambda^a_{\pm}$.
\begin{lemma}\label{lemme ope} Let $\k\in\rr^3$ be such that $\mathrm{sgn}(\k)\neq0$,  and let $\mathcal{L}^{a}_{\pm}$ be as above. The following hold:
\begin{itemize}
  \item  [(i)] If $\mathrm{sgn}(\k)\neq4$, then  $\mathcal{L}^{a}_{\pm}$ is self-adjoint with $\mathrm{dom}(\mathcal{L}^{a}_{\pm})=\mathit{H}^{1}(\S)^4$.
  \item  [(ii)]  If $\mathrm{sgn}(\k)=4$, then  $\mathcal{L}^{a}_{\pm}$ is essentially self-adjoint and we have 
  \begin{align}\label{dom5}
\mathrm{dom}(\overline{\mathcal{L}^{a}_{\pm}})=\left\{ g\in\mathit{L}^{2}(\S)^4 :\tilde{\Lambda}^{a}_{\pm}L^{1/4}[g]\in\mathit{H}^{1/2}(\S)^4  \right\}.
\end{align}
\end{itemize}
\end{lemma}
\textbf{Proof.} 
Since $L^{1/4}$ and $\mathit{C}^a_\S$ are self-adjoint operators on $\mathit{L}^{2}(\S)^4$, it follows that $\mathcal{L}^{a}_{\pm}$ is symmetric. Moreover, we have $\mathit{C}^{\infty}(\S)^4\subset\mathrm{dom}(\mathcal{L}^{a}_{\pm})\subset\mathit{L}^{2}(\S)^4$, which yields  that $\mathrm{dom}(\mathcal{L}^{a}_{\pm})$ is a dense subspace of $\mathit{L}^{2}(\S)^4$, therefore  $\mathcal{L}^{a}_{\pm}$ is closable. Let $h\in\mathrm{dom}(\mathcal{L}^{a\ast}_{\pm}) $ and let $g\in\mathit{C}^{\infty}(\S)^4 $.  By Proposition \ref{extension}-$(\mathrm{ii})$ we have 
\begin{align*}
 \langle h, \mathcal{L}^{a}_{\pm}[g] \rangle_{\mathit{L}^{2}(\S)^4}=  \langle L^{1/4}h, \Lambda^{a}_{\pm}L^{1/4}[g]\rangle_{\mathit{H}^{-1/2},\mathit{H}^{1/2}}=  \langle  \tilde{\Lambda}^{a}_{\pm}L^{1/4}h, L^{1/4}[g]\rangle_{\mathit{H}^{-1/2},\mathit{H}^{1/2}}.
  \end{align*}
  As $h\in\mathrm{dom}(\mathcal{L}^{a\ast}_{\pm}) $,  there is $f\in\mathit{L}^{2}(\S)^4$ such that
\begin{align*}
 \langle f,g \rangle_{\mathit{L}^{2}(\S)^4}=  \langle h, \mathcal{L}^{a}_{\pm}[g] \rangle_{\mathit{L}^{2}(\S)^4}=  \langle  \tilde{\Lambda}^{a}_{\pm}L^{1/4}h, L^{1/4}[g]\rangle_{\mathit{H}^{-1/2},\mathit{H}^{1/2}}.
  \end{align*}
 Hence, for all $g\in\mathit{C}^{\infty}(\S)^4 $, we get 
 \begin{align*}
 \langle L^{-1/4}[f],L^{1/4}g \rangle_{\mathit{H}^{-1/2},\mathit{H}^{1/2}}= \langle  \tilde{\Lambda}^{a}_{\pm}L^{1/4}h, L^{1/4}[g]\rangle_{\mathit{H}^{-1/2},\mathit{H}^{1/2}},
  \end{align*}
 which implies that $ \tilde{\Lambda}^{a}_{\pm}L^{1/4}[h]=L^{-1/4}[f]$ holds in $\mathit{H}^{-1/2}(\S)^4$ and then in $\mathit{H}^{1/2}(\S)^4$. Therefore  $ \tilde{\Lambda}^{a}_{\pm}L^{1/4}[h]\in  \mathit{H}^{1/2}(\S)^4$, and we have the inclusion 
 \begin{align*}
\mathrm{dom}(\mathcal{L}^{a\ast}_{\pm})\subset \left\{ g\in\mathit{L}^{2}(\S)^4 : \tilde{\Lambda^{a}_{\pm}}L^{1/4}[g]\in\mathit{H}^{1/2}(\S)^4  \right\}.
\end{align*}
Now, one can easily check  the other inclusion and we thus get the equality. Hence, item $(\mathrm{i})$ is an immediate consequence of Lemma \ref{commutator} and \eqref{multi}. To prove the second item, it is sufficient to show that $\mathcal{L}^{a\ast}_{\pm}\subset\overline{\mathcal{L}^{a}_{\pm}}$. For this, one can take the sequence of functions defined by \eqref{suite} (just switch the roles of $ \tilde{\Lambda}^{a}_{\pm}$ and $ \tilde{\Lambda}^{a}_{\mp}$) and use the fact that $\tilde{\Lambda}^{a}_{\pm}\tilde{\Lambda}^{a}_{\mp}$ are continuous from $\mathit{H}^{-1/2}(\S)^4$ to $\mathit{H}^{1/2}(\S)^4$, we omit the details. This finishes the proof of the lemma. \qed
\newline

Note that,  for any $ \psi=u+\Phi[g]\in\mathrm{dom}(\mathcal{H}_{\k})$ and $\varphi=v+\Phi[h]\in\mathrm{dom}(\mathcal{H}^{\ast}_{\k})$, it holds that
\begin{align}\label{similarformula} 
\langle  \mathcal{H}^{\ast}_{\k}\varphi ,\psi\rangle_{\mathit{L}^2(\rr^3)^4}- \langle \varphi, \mathcal{H}_{\k}\psi\rangle_{\mathit{L}^2(\rr^3)^4} = \langle -\tilde{\Lambda}_+[h], g\rangle_{\mathit{H}^{-1/2},\mathit{H}^{1/2}}- \langle h ,-\Lambda_+[g]\rangle_{\mathit{H}^{-1/2},\mathit{H}^{1/2}}.
\end{align}
Taking into account the above lemma, from \eqref{relationself}  and \eqref{similarformula} it easily follows that:
\begin{align}
  \mathcal{H}_{\k} \text { is (essentially) self-adjoint }\Longleftrightarrow   \mathcal{L}_+ \text{ is (essentially) self-adjoint.}
  \end{align}

As it was mentioned in the introduction, the operator $\mathcal{L}_{+}$ appears in this form when we study the self-adjoint extension of $\mathcal{H}_{\k}$ from the point of view of the boundary triples theory (see \cite{BH} and \cite{BHOP}; for a more general view of the theory we refer to \cite{BHS} and \cite{BGP} for example). Indeed, denote by $S:=\mathcal{H}\downharpoonright_{\mathit{H}^{1}_{0}(\rr^3\setminus\S)^4}$ and let $T:=\mathcal{H}$ with 
\begin{align*}
\mathrm{dom}(T)= \{ u+\Phi [g]: u\in\mathit{H}^1(\rr^3)^4, g\in\mathit{L}^{2}(\S)^4 \},
\end{align*}
and define the linear mappings $\Gamma_{1},\Gamma_{2}:\mathrm{dom}(T)\longrightarrow \mathit{L}^{2}(\S)^4$ by
\begin{align}\label{mapp}
\Gamma_{1}(\varphi)=g \quad\text{ and}\quad \Gamma_{2}(\varphi)= t_{\S}u+ \mathit{C}_{\S}[g].
\end{align}
Then, $\{ \mathit{L}^{2}(\S)^4, \Gamma_1,\Gamma_2\}$ is a quasi-boundary triples for $\overline{T}=S^{\ast}$ (adapt the arguments of \cite{BH} or \cite{BHOP}). Moreover, if  we denote by $ \tilde{\G}_1$ the extension of $\Gamma_{1}$, that is  $\tilde{\Gamma}_{1}:\mathrm{dom}(\overline{T})\longrightarrow \mathit{H}^{-1/2}(\S)^4$, we then get that $\{ \mathit{L}^{2}(\S)^4, L^{-1/4}\tilde{\Gamma}_1,L^{1/4}\Gamma_2\}$ is an ordinary boundary triple for $\overline{T}=S^{\ast}$. Now it is easy to check that 
\begin{align}
\mathcal{H}_{\k}= T\downharpoonright\mathrm{Kr}( (\ep\mathit{I}_4+\m\beta +\e(\alpha\cdot N))\G_2+\G_1)=T\downharpoonright\mathrm{Kr}(\G_2-\mathit{C}_{\S}\G_1+\Lambda_+\G_1).
\end{align}
Thus, after transforming the quasi-boundary triples to an ordinary boundary triples (see \cite[Theorem 4.5]{BH} for example) it follows that:  $\mathcal{H}_{\k}$ is self-adjoint (respectively essentially self-adjoint) if and only if $\mathcal{L}_+$ is self-adjoint (respectively essentially self-adjoint); see \cite[Corollary 2.8]{BH}.

 
 \section{Spectral properties} \label{sec4}
 \setcounter{equation}{0}
In this section, we examine the spectral properties of the operator  $\mathcal{H}_{\k}$. First, we give a necessary condition for the existence of the point spectrum in the gap $(-m,m)$ and a Krein-type resolvent formula.  More precisely, recall that $\mathrm{sgn}(\k)$ is defined in \eqref{sgndef}, then we have the following.
\begin{proposition}\label{BS-K} Let  $ \mathcal{H}_{\k}$ be as in the definition \ref{def1}. If $\mathrm{sgn}(\k)=4$ , then the following hold:
\begin{itemize}
 \item [(i)] Given $a\in(-m,m)$, then one has  $\mathrm{Kr}(\overline{\mathcal{H}_{\k}}-a)\neq\{0\}$ $\Longleftrightarrow$ $\mathrm{Kr}( \tilde{\Lambda}^{a}_{+})\neq\{0\}$ (Birman-Schwinger principle) and    $\mathrm{Kr}(\overline{\mathcal{H}_{\k}}-a)=\{ \Phi^a[g]:g\in\mathrm{Kr}( \tilde{\Lambda}^{a}_{+})\}$.

\item [(ii)] For all $z\in\cc\setminus\rr$ the operator $ \tilde{\Lambda}^{z}_{+}$ is bounded invertible from $\mathit{H}^{-1/2}(\S)^4$ to $\mathit{H}^{1/2}(\S)^4$ and we have 
 \begin{align}\label{resolvent}
    (\overline{\mathcal{H}_{\k}} -z)^{-1}= (\mathcal{H} -z)^{-1} - \Phi^{z}(\tilde{\Lambda}^{z}_{+} )^{-1} (\Phi^{\overline{z}})^\ast.
\end{align}
\end{itemize}
If $\mathrm{sgn}(\k)\neq0,4$ and $z\in\cc\setminus\rr$, then $ \Lambda^{z}_{+}$  is bounded invertible from $\mathit{H}^{1/2}(\S)^4$ to $\mathit{H}^{1/2}(\S)^4$ and  the above statements hold true with $ \Lambda^{\bullet}_{+}$  instead of $\tilde{\Lambda}^{\bullet}_{+} $. 
\end{proposition}

\textbf{Proof.} We prove the statements for  $\mathrm{sgn}(\k)=4$, the case $\mathrm{sgn}(\k)\neq0,4$  follows  the same lines. 

$(\mathrm{i})$ Let us prove the implication $(\Rightarrow)$ and the inclusion $\mathrm{Kr}(\overline{\mathcal{H}_{\k}}-a)\subset\{ \Phi^a[g]:g\in\mathrm{Kr}( \tilde{\Lambda}^{a}_{+})\}$.  Let $a\in(-m,m)$ and assume that there is a nonzero $\varphi=u+\Phi[g]\in\mathrm{dom}(\overline{\mathcal{H}_{\k}})$ such that $\overline{\mathcal{H}_{\k}}\varphi=a\varphi$. First observe that $\tilde{\Lambda}^{a}_{+}-\tilde{\Lambda}_{+}=\tilde{\mathit{C}^{a}_{\S}}-\tilde{\mathit{C}_{\S}}$. Now, using the definition of  $\overline{\mathcal{H}_{\k}}$, we then get
 \begin{align}\label{a1}
  \mathcal{H}u=a\varphi=a(u+\Phi[g]).
\end{align}  
From this we deduce that $(\mathcal{H}-a)\mathcal{H}u=ag\delta_{\S}$ holds in $\mathcal{D}^{\prime}(\rr^3)^4$, and therefore 
\begin{align}\label{a2}
\mathcal{H}u=a\Phi^a[g].
\end{align}
From this, it is clear that if $a=0$ then $u=0$. Therefore, $\varphi=\Phi[g]\neq0$  (with $g\neq0$, as otherwise $\varphi$ would be zero)  and  $g\in\mathrm{Kr}( \tilde{\Lambda}_{+})$,  which yields that  $\mathrm{Kr}(\overline{\mathcal{H}_{\k}})\subset\{ \Phi[g]:g\in\mathrm{Kr}( \tilde{\Lambda}_{+})\}$. Now assume that $a\neq0$, then from \eqref{a1} and  \eqref{a2} it follows that  $u= (\Phi^a-\Phi)[g]$. Since $\varphi=u+\Phi[g]\in\mathrm{dom}(\overline{\mathcal{H}_{\k}})$, it holds that $ t_{\S}u=-\tilde{\Lambda}_{+}[g]$, and by Proposition \ref{extension}$(\mathrm{iii})$ we also get  that  $ t_{\S}u=(\tilde{\mathit{C}^{a}_{\S}}-\tilde{\mathit{C}_{\S}})[g]=-\tilde{\Lambda}_{+}[g]$. Hence, we obtain that   $0\neq g\in\mathrm{Kr}( \tilde{\Lambda}^{a}_{+})$ and $\varphi=\Phi^a[g]$, therefore $\mathrm{Kr}(\overline{\mathcal{H}_{\k}}-a)\subset\{ \Phi^a[g]:g\in\mathrm{Kr}( \tilde{\Lambda}^{a}_{+})\}$. 

Conversely, let $a\in(-m,m)$ be such that $\tilde{\Lambda}^{a}_{+}[g]=0$,  for a nonzero $g\in\mathit{H}^{-1/2}(\S)^4$.  Then, it is clear that  $\varphi=\Phi[g]\in\mathrm{dom}(\overline{\mathcal{H}_{\k}})$ and we have  $0\neq\varphi\in\mathrm{Kr}(\overline{\mathcal{H}_{\k}})$ when $a=0$, which gives the result in this case. Now suppose that $a\neq0$, let $u=a \mathcal{H}^{-1}\Phi^a[g]\in\mathit{H}^1(\rr^3)^4$ and set $\varphi=u+\Phi[g]$. Then $\mathcal{H}u=a\Phi^a[g]$ and $( \mathcal{H}-a)u=a\Phi[g]$ in $\mathcal{D}^{\prime}(\rr^3)^4$, this amounts to saying that $\overline{\mathcal{H}_{\k}}\varphi=\mathcal{H}u=a(u+\Phi[g])=a\varphi$ and $u= \Phi^a[g]-\Phi[g]$. Furthermore, we can easily see that $ t_{\S}u=(\tilde{\mathit{C}^{a}_{\S}}-\tilde{\mathit{C}_{\S}})[g]=-\tilde{\Lambda}_{+}[g]$. Summing up, we have proved that  $\varphi=\Phi^{a}[g]\in\mathrm{dom}(\overline{\mathcal{H}_{\k}})$  and $\overline{\mathcal{H}_{\k}}\varphi=a\varphi$, which yields that $\varphi\in\mathrm{Kr}(\overline{\mathcal{H}_{\k}}-a)$. This ends the proof of $(\mathrm{i})$.

$(\mathrm{ii})$ Fix $z\in\cc\setminus\rr$. Since $\overline{\mathcal{H}_{\k}}$ is self-adjoint it follows that $ (\overline{\mathcal{H}_{\k}} -z)^{-1}$ is well defined and  bounded. Moreover, following the same arguments of the proof of $(\mathrm{i})$ one can see that $\mathrm{Kr}(\tilde{\Lambda}^{z}_{+})=\{ 0\}$, as otherwise  $z$ would be a non-real eigenvalue of $\overline{\mathcal{H}_{\k}}$. Let $u\in\mathit{L}^2(\rr^3)^4$ and set $\varphi:= (\overline{\mathcal{H}_{\k}} -z)^{-1}u\in \mathrm{dom}(\overline{\mathcal{H}_{\k}})$. Thanks to Remark \ref{appelrem}, we know that there are unique functions $v\in\mathit{H}^1(\rr^3)^4$ and $g \in\mathit{H}^{-1/2}(\S)^4$ such that $\varphi= v+ \Phi^{z}[g]$. Moreover one has  $(\overline{\mathcal{H}_{\k}} -z)\varphi=(\mathcal{H} -z)v $, and thus $v=(\mathcal{H} -z)^{-1}u$, which means actually that $ \varphi=(\mathcal{H} -z)^{-1}u+ \Phi^{z}[g]$. Next, observe that 
\begin{align}\label{}
   i\alpha\cdot\mathit{N}(t_{\S}\varphi_+ -t_{\S}\varphi_-) = g \quad \text{ and }\quad \frac{1}{2}(t_{\S}\varphi_+ +t_{\S}\varphi_-) =(\mathcal{H} -z)^{-1}u\downharpoonright_{\S} +\tilde{\mathit{C}}^{z}_{\S}[g].  
\end{align}
Using that $(\mathcal{H} -z)^{-1}u\downharpoonright_{\S} =(\Phi^{\overline{z}})^\ast u$ and the transmission condition  \eqref{trans}, we obtain that $\tilde{\Lambda}^{z}_{+}[g]=-(\Phi^{\overline{z}})^\ast u\in \mathit{H}^{1/2}(\S)^4$. Since this is true for all $u\in\mathit{L}^2(\rr^3)^4$,  we then get that $\mathrm{Rn}(\tilde{\Lambda}^{z}_{+})=\mathit{H}^{1/2}(\S)^4$. Hence $ \tilde{\Lambda}^{z}_{+}:\mathit{H}^{-1/2}(\S)^4 \longrightarrow\mathit{H}^{1/2}(\S)^4$ is a bounded bijective operator.  Summing up, we have proved that $ (\overline{\mathcal{H}_{\k}} -z)^{-1}u=(\mathcal{H} -z)^{-1}u- \Phi^{z}[(\tilde{\Lambda}^{z}_{+} )^{-1} (\Phi^{\overline{z}})^\ast u]$ holds for all $u\in\mathit{L}^2(\rr^3)^4$, which proves the identity \eqref{resolvent}.\qed

\begin{remark}\label{rr}  A careful inspection of the argument used above reveals that: $a\in(-m,m)$ is an isolated point of $\mathrm{Sp}(\overline{\mathcal{H}_{\k}})$ if and only if $0$ is an isolated point of $\mathrm{Sp}(\Lambda^{a}_{+})$  (respectively $\mathrm{Sp}(\tilde{\Lambda}^{a}_{+})$ when $\mathrm{sgn}(\k)=4$), and as $\Phi^z$ is injective we have that $\mathrm{dim}\mathrm{Kr}(\overline{\mathcal{H}_{\k}}-a)$ is equal to $\mathrm{dim}\mathrm{Kr}( \Lambda^{a}_{+})$  (respectively  $\mathrm{dim}\mathrm{Kr}( \tilde{\Lambda}^{a}_{+})$ when $\mathrm{sgn}(\k)=4$). Furthermore, item $(\mathrm{ii})$ holds true for all $z\in\rho(\overline{\mathcal{H}_{\k}})\cap\rho(\mathcal{H})$.
\end{remark}

As an immediate consequence of Lemma \ref{lemme ope}, Proposition \ref{BS-K} and Remark \ref{rr}, we have the following result.
\begin{corollary}\label{corollaire inter} Let  $ \mathcal{H}_{\k}$ be as in the definition \ref{def1}. The following holds:
\begin{itemize}
 \item [(i)] For all $ a\in(-m,m)$, one has  
 \begin{align*}
 a\in\mathrm{Sp}_{\mathrm{p}}(\overline{\mathcal{H}_{\k}})&\Longleftrightarrow 0\in\mathrm{Sp}_{\mathrm{p}}(\overline{\mathcal{L}^{a}_{+}}),\\
  a\in\mathrm{Sp}_{\mathrm{ess}}(\overline{\mathcal{H}_{\k}})&\Longleftrightarrow 0\in\mathrm{Sp}_{\mathrm{ess}}( \overline{\mathcal{L}^{a}_{+}}).
  \end{align*}
\item [(ii)] If $\mathrm{sgn}(\k)=4$, then for all $z\in\rho(\overline{\mathcal{H}_{\k}})\cap\rho(\mathcal{H})$ the operator $\overline{\mathcal{L}^{z}_{+}}$ is bounded invertible from $\mathit{L}^{2}(\S)^4$ to $\mathit{L}^{2}(\S)^4$ and we have 
 \begin{align}\label{resolvent2}
    (\overline{\mathcal{H}_{\k}} -z)^{-1}= (\mathcal{H} -z)^{-1} - \Phi^{z}L^{\frac{1}{4}}\left(\overline{\mathcal{L}^{z}_{+}}\right)^{-1}L^{\frac{1}{4}} (\Phi^{\overline{z}})^\ast.
\end{align}
\end{itemize}  
\end{corollary}

  In the remainder of this section, we focus namely on the spectral properties of $\mathcal{H}_{\k}$ when $\S$ satisfies the assumption $(\mathrm{H2})$. In order to avoid ambiguities we use the following notations:
\begin{notation} For all $\nu\geqslant0$, we denote by $\overline{\mathcal{H}^{\nu}_{\k}}$ (respectively  $\Phi^{z}_\nu$, $ \tilde{\Lambda}^{z}_{+,\nu}$ and $(\Phi^{\overline{z}}_{\nu})^\ast$) the operator  $\overline{\mathcal{H}_{\k}}$ (respectively  $\Phi^{z}$, $ \tilde{\Lambda}^{z}_{+}$ and $(\Phi^{\overline{z}})^\ast $) whenever $\S=\S_\nu$, i.e $\S$ satisfies $(\mathrm{H2})$, and we write $\overline{\mathcal{H}_{k}}$ (respectively  $\Phi^{z}$, $ \tilde{\Lambda}^{z}_{+}$ and $(\Phi^{\overline{z}})^\ast$) instead of $\overline{\mathcal{H}^0_{\k}}$ (respectively  $\Phi^{z}_0$, $ \tilde{\Lambda}^{z}_{+,0}$ and $(\Phi^{\overline{z}}_{0})^\ast$), i.e  when $\nu=0$. 
\end{notation}
\subsection{Non-critical case}
This part deals with  the basic spectral properties of $\mathcal{H}_{\k}$ when $\kappa=(\ep,\m,0)$ (i.e $\e=0$) and $\mathrm{sgn}(\k)\neq0,4$.  The following theorem gives us a complete description of the essential spectrum of  $\mathcal{H}_{\k}$ when $\S$ satisfies $(\mathrm{H2})$ and $\e=0$.  We would like  to thank  the authors of \cite{BHT} for revealing an error in the original version about the fact that the essential spectrum can emerge in the gap $(-m,m)$.


\begin{theorem}\label{Prop non} Let $\k\in\rr^2\times\{0\}$  be such that $\mathrm{sgn}(\k)=\ep^2-\m^2\neq0,4$, and suppose that $\S$ satisfies $(\mathrm{H2})$ with $\nu\geqslant0$. Set 
\begin{align}\label{defa_+}
\begin{split}
    a_\pm=&m \frac{-16\ep\m \pm (\mathrm{sgn}(\k)-4)^2\sqrt{\frac{(\mathrm{sgn}(\k)+4)^2}{(\mathrm{sgn}(\k)-4)^2}}}{(\mathrm{sgn}(\k)-4)^2+16\ep^2},   \\
    a^{\ast}=&-m \frac{-16\ep\m}{(\mathrm{sgn}(\k)-4)^2+16\ep^2}.
    \end{split}
\end{align}
 The following hold:
  \begin{itemize}
  \item [(i)] If $\ep^2-\mu^2>4$,  then 
  \begin{equation*}\label{} 
  \mathrm{Sp}_{\mathrm{ess}}(\mathcal{H}^{\nu}_{\k})= \left\{
  \begin{aligned}
(-\infty,-m]\cup [a_+,+\infty)&, \, \text{ for } \ep>0 \text{ and } \m\in\rr,\\
(-\infty,a_-]\cup [m,+\infty)&,\, \text{ for } \ep<0 \text{ and } \m\in\rr.
\end{aligned}
  \right.
\end{equation*} 
  \item[(ii)] If $0<\ep^2-\mu^2<4$,  then 
  \begin{equation*}\label{} 
  \mathrm{Sp}_{\mathrm{ess}}(\mathcal{H}^{\nu}_{\k})= \left\{
  \begin{aligned}
(-\infty,a_-]\cup [m,+\infty)&, \, \text{ for } \ep>0 \text{ and } \m\in\rr,\\
(-\infty,-m]\cup [a_+,+\infty)&,\, \text{ for } \ep<0 \text{ and } \m\in\rr.
\end{aligned}
  \right.
\end{equation*} 
    \item[(iii)] If $-4<\ep^2-\mu^2<0$,  then 
  \begin{equation*}\label{} 
  \mathrm{Sp}_{\mathrm{ess}}(\mathcal{H}^{\nu}_{\k})= \left\{
  \begin{aligned}
(-\infty,-m]\cup [m,+\infty)&, \, \text{ for } \m>0 \text{ and } \ep\in\rr,\\
(-\infty,a_-]\cup [a_+,+\infty)&,\, \text{ for } \m<0 \text{ and } \ep\in\rr.
\end{aligned}
  \right.
\end{equation*} 

      \item[(iv)]  If $\ep^2-\mu^2=-4$,  then 
  \begin{equation*}\label{} 
  \mathrm{Sp}_{\mathrm{ess}}(\mathcal{H}^{\nu}_{\k})= \left\{
  \begin{aligned}
  (-\infty,-m]\cup [m,+\infty)&, \, \text{ for } \m>0 \text{ and } \ep\in\rr,\\
(-\infty,a^{\ast}]\cup [m,+\infty)&, \,  \text{ for } \m<0 \text{ and } \ep>0,\\
(-\infty,-m]\cup [a^{\ast},+\infty)&,\, \text{ for } \m<0 \text{ and } \ep<0,\\
\rr,\quad&\, \text{ for } \m=-2 \text{ and } \ep=0.
\end{aligned}
  \right.
\end{equation*} 
       \item[(v)] If $\ep^2-\mu^2<-4$,  then 
  \begin{equation*}\label{} 
  \mathrm{Sp}_{\mathrm{ess}}(\mathcal{H}^{\nu}_{\k})= \left\{
  \begin{aligned}
(-\infty,-m]\cup [m,+\infty)&, \, \text{ for }  \m>0 \text{ and } \ep\in\rr,\\
(-\infty,a_-]\cup [a_+,+\infty)&,\, \text{ for } \m<0 \text{ and } \ep\in\rr.
\end{aligned}
  \right.
\end{equation*} 
\end{itemize}
Furthermore,  we have $\mathrm{Sp}(\mathcal{H}_{\k})=\mathrm{Sp}_{\mathrm{ess}}(\mathcal{H}_{\k})$.
\end{theorem}
\begin{remark} Note that  a similar statement can be formulated for $\eta\neq0$, in that case we have  
\begin{align}\label{defa_+bis}
    a_\pm=&m \frac{-16\ep\m \pm (\mathrm{sgn}(\k)-4)^2\sqrt{\frac{(\mathrm{sgn}(\k)+4)^2+16\e^2}{(\mathrm{sgn}(\k)-4)^2}}}{(\mathrm{sgn}(\k)-4)^2+16\ep^2}. 
\end{align}
The statements $(\mathrm{i})$ and $(\mathrm{ii})$ still hold true, and for $\mathrm{sgn}(\k)<0$ several cases should be taken into account,  so for the sake of readability we chose not to write it here.
\end{remark}

\textbf{Proof.}  We first prove assertions $(\mathrm{i})-(\mathrm{v})$  when $\nu=0$, we then use compactness arguments and  Proposition \ref{BS-K}$(\mathrm{ii})$ to get the result when $\nu>0$. To this end and for the convenience of the reader we divide the proof in three steps.

\underline{\it{Step 1.}} We analyze the spectrum of $\mathcal{H}_{\k}$ in the gap $(-m,m)$. For that, let  $a\in (-m,m)$ and set  
$$\Gamma_{\pm m,\pm a}(\xi)=\left[ \alpha\cdot(\xi_1,\xi_2,0)\pm m\beta \pm a \right].$$
 Since the $\alpha_j$'s anticommute with $\beta$, a simple computation shows that
\begin{align}\label{pm}
\begin{split}
\left(\Gamma_{ m, a}(\xi)\right)^2&= |\xi|^2+m^2-a^2 +2a \Gamma_{ m, a}(\xi),\\
\Gamma_{- m,- a}(\xi)\Gamma_{ m, a}(\xi)&= |\xi|^2+m^2-a^2 -2m\beta\Gamma_{ m, a}(\xi),\\
\Gamma_{ m, -a}(\xi) \Gamma_{ m, a}(\xi)&= |\xi|^2+m^2-a^2.
\end{split}
\end{align}
 Using  the Fourier-Plancherel operator, it is not hard to prove that $\Lambda^{a}_{+}$ is unitary equivalent to the following multiplication operator:
\begin{align}
\Pi^{a}_{+}:= \frac{1}{\mathrm{sgn}(\k)}(\ep I_4 -\mu\beta)+ \frac{1}{2\sqrt{|\xi|^2+m^2-a^2}}\Gamma_{ m, a}(\xi).
\end{align}
Moreover, taking  into account the properties \eqref{pm}, a simple computation shows that $\Pi^{a}_{+}$ is invertible and its inverse is given explicitly by
\begin{align} \label{inverse}
(\Pi^{a}_{+})^{-1}= C^{-1}\left( 1 + \frac{\ep a+\m m}{\sqrt{|\xi|^2+m^2-a^2}} -\frac{(\ep+\m\beta)}{2\sqrt{|\xi|^2+m^2-a^2}}\Gamma_{ m, a}(\xi)\right)(\ep \mathit{I}_4+\m\beta),
\end{align}
if and only if $C\neq0$ for all $\xi\in\rr^2$, where $C$ is given by
 $$C= \frac{4-\mathrm{sgn}(\k)}{4}+\frac{\ep a+\m m}{\sqrt{|\xi|^2+m^2-a^2}}.$$
 Since $\mathrm{sgn}(\k)\neq4$, it follows that $-m\m/\ep\notin\mathrm{Sp}(\mathcal{H}_{\k})$, for all $\ep\neq0$. In the following, we always  assume that $a\neq-m\m/\ep$ when $\ep\neq0$, and we look for the values of $a$ for which we have  $C=0$. Note that  
\begin{align}\label{condi1}
 C= 0\Longleftrightarrow  \sqrt{|\xi|^2+m^2-a^2}=\frac{4(\ep a+\m m)}{\mathrm{sgn}(\k)-4}.
 \end{align}
 Thus, $C=0$ for some $|\xi|\in\rr_+$, only if  
\begin{align}\label{condi2}
 \frac{4(\ep a+\m m)}{\mathrm{sgn}(\k)-4}>0.
 \end{align}
 Assume that \eqref{condi2} holds true, then $C= 0$ if and only if $ |\xi|^2=P(a)$, where the polynomial $P(a)$ is given by
 \begin{align}\label{condi3}
 P(a)= \frac{(\mathrm{sgn}(\k)-4)^2+16\ep^2}{(\mathrm{sgn}(\k)-4)^2}a^2+ \frac{32\ep\m m}{(\mathrm{sgn}(\k)-4)^2}a -  \frac{(\mathrm{sgn}(\k)-4)^2-16\m^2}{(\mathrm{sgn}(\k)-4)^2}m^2.
 \end{align}
Recall $a_+,a_-$ and $a^{\ast}$ from \eqref{defa_+}, then  $a_+$ and $a_-$ are the zeros of $P(a)$ when $\ep^2-\mu^2\neq-4$, and $a^{\ast}$ is a double root of $P(a)$ when $\ep^2-\mu^2=-4$. Thus $P(a)\geqslant0$ if and only if $a\geqslant a_+$ or $a\leqslant a_-$.  In the remainder of the proof we deal with assertion $(\mathrm{i})$, the other assertions follow in the same way. Assume that $\ep^2-\mu^2>4$, then 
 \begin{align*}
    a_\pm=&m \frac{-16\ep\m \pm (\mathrm{sgn}(\k)-4)(\mathrm{sgn}(\k)+4)}{(\mathrm{sgn}(\k)-4)^2+16\ep^2}\quad\text{ and } -m<a_-<a_+<m.
\end{align*}
 As $\mathrm{sgn}(\k)>4$, it follows that  the condition \eqref{condi2}  is equivalent to
  \begin{align*}
 a>-\frac{\m m}{ \ep}\, \text{ if } \ep>0\, \quad\text{ or }\quad \, a<-\frac{\m m}{ \ep} \,\text{ if } \ep<0.
\end{align*}
 Now using the fact that $\ep^2>\m^2$, a simple computation yields
  \begin{align*}
 a_+>-\frac{\m m}{ \ep}\, \text{ and }\, a_-<-\frac{\m m}{ \ep}. 
\end{align*}
 Hence,  if  $\ep>0$ (respectively $\ep<0$)  then for all $a\geqslant a_+$ (respectively $a\leqslant a_-$) we have   $P(a)\geqslant0$ and the condition  \eqref{condi2} holds true. Consequently, the set of $\xi$ for which  $C=0$ is given by the circle  $\{ \xi: |\xi|= \sqrt{P(a)}\}$, and in that case $0$ is in the essential spectrum of $\Lambda^{a}_+$. Therefore we conclude by Proposition \ref{BS-K} that 
 \begin{align}\label{concluA}
 \begin{split}
    (a_+,m)\subset\mathrm{Sp}_{\mathrm{ess}}(\mathcal{H}_{\k})& \text{ and } (-m,a_+) \subset \rho(\mathcal{H}_{\k}), \text{ for } \ep>0 \\
  (-m,a_-) \subset\mathrm{Sp}_{\mathrm{ess}}(\mathcal{H}_{\k}) &  \text{ and }(a_-,m) \subset \rho(\mathcal{H}_{\k}), \text{ for } \ep<0.
  \end{split}
\end{align}

\underline{\it{Step 2.}} Now we prove the inclusion $(-\infty,-m)\cup (m,+\infty)\subset\mathrm{Sp}_{\mathrm{ess}}(\mathcal{H}_{\k})$, for that  we construct a singular sequence  for $\mathcal{H}_{\k}$ and $a$. Fix $a\in (-\infty,-m)\cup(m,\infty)$ and define 
\begin{equation*}\label{} 
  \varphi: \left\{
  \begin{aligned}
\rr^3 &\longrightarrow \cc^4\\
(\overline{x},x_3)&\longmapsto \left(\frac{\xi_1-i\xi_2}{a-m},0,0,1\right)^t e^{i\overline{x}\cdot\xi},
\end{aligned}
  \right.
\end{equation*}
here $\xi= (\xi_1,\xi_2)$ and $|\xi|^2=a^2-m^2$. Observe that we have $ (-i\alpha \cdot \nabla + m\beta-a)\varphi=0$. Let $R>0$, $\chi\in \mathit{C}^{\infty}_{0}(\rr^2)$ and $\theta\in  \mathit{C}^{\infty}_{0}([0,\infty[,\rr)$ such that
\begin{equation*}\label{} 
  \theta(r)= \left\{
  \begin{aligned}
&1 \quad\text{for } r\in [2R, 3R],\\
&0 \quad \text{for }  r\in [0, R].
\end{aligned}
  \right.
\end{equation*}
For $n\in\nn^\star$, we define the sequences of functions
\begin{align}\label{suitee}
\begin{split}
  \displaystyle \varphi_{+,n}(\overline{x},x_3)&=  n^{-\frac{3}{2}} \varphi(\overline{x},x_3)\chi(\overline{x}/n)\theta(x_3/n) \,\quad \text{for }x_3>0,\\
  \displaystyle   \varphi_{-,n}(\overline{x},x_3)&=  n^{-\frac{3}{2}} \varphi(\overline{x},x_3)\chi(\overline{x}/n)\theta(-x_3/n) \, \quad\text{for }x_3<0.
  \end{split}
\end{align}
It's clear that $ \varphi_{\pm,n}\in\mathit{H}^1(\O_\pm)$ and $t_\Sigma \varphi_{\pm,n}=0$, thus $ \varphi_{n}:=(\varphi_{+,n},\varphi_{-,n})\in\mathrm{dom}(\mathcal{H}_{\k})$. Moreover,  $ (\varphi_{n})_{n\in\mathbb{N}^\star}$  converges weakly to zero and  we have
\begin{align*}\label{}
    \Vert \varphi_{n}\Vert_{\mathit{L}^2(\rr^3)^4}^2=  \Vert \varphi_{+,n}\Vert_{\mathit{L}^2(\O_+)^4}^2+  \Vert \varphi_{-,n}\Vert_{\mathit{L}^2(\O_-)^4}^2=\frac{2a}{a-m}\Vert \chi\Vert^{2}_{\mathit{L}^2(\rr^2)}\Vert \theta\Vert^{2}_{\mathit{L}^2(\rr_+)}>0,
\end{align*}
and
\begin{align*}\label{}
      \Vert\left(-i\alpha \cdot \nabla + m\beta-a\right) \varphi_{n}\Vert_{\mathit{L}^2(\rr^3)^4}^2 =&   \Vert\left(-i\alpha \cdot \nabla + m\beta-a\right) \varphi_{+,n}\Vert_{\mathit{L}^2(\O_+)^4}^2 \\
      &+ \Vert\left(-i\alpha \cdot \nabla + m\beta-a\right) \varphi_{-,n}\Vert_{\mathit{L}^2(\O_-)^4}^2\\
    &\leqslant  \frac{4a}{n^2(a-m)}\bigg( \Vert\nabla \eta\Vert_{\mathit{L}^2(\rr^2)}^2 \Vert \theta\Vert_{\mathit{L}^2(\rr_+)}^2 +\Vert \chi\Vert_{\mathit{L}^2(\rr^2)}^2 \Vert \theta^{\prime}\Vert_{\mathit{L}^2(\rr_+)}^2\bigg).
\end{align*}
Thus, we get 
\begin{align*}
\frac{  \Vert\left(-i\alpha \cdot \nabla + m\beta-a\right) \varphi_{n}\Vert_{\mathit{L}^2(\rr^3)^4}}{ \Vert \varphi_{n}\Vert_{\mathit{L}^2(\rr^3)^4}}\xrightarrow[n\to\infty]{}0.
\end{align*}
From this and  \textit{Step 1},  we deduce that 
\begin{align*}
(-\infty,-m)&\cup(a_+,m)\cup(m,\infty)\subset \mathrm{Sp}_{\mathrm{ess}}(\mathcal{H}_{\k})\subset \mathrm{Sp}(\mathcal{H}_{\k})\subset(-\infty,-m]\cup[a_+,\infty),  \text{ for } \ep>0, \\
  (-\infty,-m)&\cup(-m,a_-)\cup(m,\infty)\subset \mathrm{Sp}_{\mathrm{ess}}(\mathcal{H}_{\k})\subset \mathrm{Sp}(\mathcal{H}_{\k})\subset(-\infty,a_-]\cup[m,\infty), \text{ for } \ep<0.
\end{align*}
Since the spectrum of a self-adjoint operator is closed, the end-points  also belong to the spectrum, and hence for $\ep^2-\mu^2>4$, we get  
  \begin{equation*}\label{} 
  \mathrm{Sp}(\mathcal{H}_{\k})= \mathrm{Sp}_{\mathrm{ess}}(\mathcal{H}_{\k})= \left\{
  \begin{aligned}
(-\infty,-m]\cup [a_+,+\infty)&, \, \text{ for } \ep>0 \text{ and } \m\in\rr,\\
(-\infty,a_-]\cup [m,+\infty)&,\, \text{ for } \ep<0 \text{ and } \m\in\rr,
\end{aligned}
  \right.
\end{equation*} 
 which proves the result when  $\nu=0$.

\underline{\it{Step 3.}} Assume that $\nu>0$,  let us show the equality $\mathrm{Sp}_\mathrm{ess}(\mathcal{H}^{\nu}_{\k})=\mathrm{Sp}_{\mathrm{ess}}(\mathcal{H}_{\k})$. Recall that $F$ is the flat part of $\S_\nu$ given by  \eqref{plat}; fix $z\in\cc\setminus\rr$ and  let  $\mathcal{T}: \mathit{L}^2(\rr^3)^4\rightarrow\mathit{L}^2(\rr^3)^4$ be the bounded operator defined by
\begin{align}\label{T}
  \mathcal{T}=   \Phi^{z}_\nu(\Lambda^{z}_{+,\nu} )^{-1} (\Phi^{\overline{z}}_{\nu})^\ast - \Phi^{z}(\Lambda^{z}_{+} )^{-1} (\Phi^{\overline{z}})^\ast.
\end{align}
Then  $\mathcal{T}$ is a compact operator in $\mathit{L}^2(\rr^3)^4$. Indeed, observe that $\mathcal{T}$ can be written as follows:
\begin{align*}\label{}
    \mathcal{T}=&\left( \Phi^{z}_\nu(\Lambda^{z}_{+,\nu} )^{-1} - \Phi^{z}(\Lambda^{z}_{+} )^{-1} \right)  (\mathcal{H} -z)^{-1}\downharpoonright_{F} + \Phi^{z}_\nu(\Lambda^{z}_{+,\nu} )^{-1}  (\mathcal{H} -z)^{-1}\downharpoonright_{\S_\nu\setminus F}\\
    &- \Phi^{z} (\Lambda^{z}_{+} )^{-1}  (\mathcal{H} -z)^{-1} \downharpoonright_{\S_0 \setminus F} := \mathcal{T}_1+\mathcal{T}_2+\mathcal{T}_3.
\end{align*}
Since $ \S_\nu\setminus F$ is compact for all $\nu\geqslant0$, it follows that the injection  $\mathit{H}^{1/2}(\S_\nu\setminus F)^4\hookrightarrow\mathit{L}^{2}(\S_\nu)^4$ is compact. Using this and the fact that  $(\mathcal{H} -z)^{-1}\downharpoonright_{\S_\nu\setminus F}$ is bounded from $\mathit{L}^2(\rr^3)^4$ to $\mathit{H}^{1/2}(\S_\nu\setminus F)^4$,  it holds that $(\mathcal{H} -z)^{-1}\downharpoonright_{\S_\nu\setminus F}$ is a compact operator from $\mathit{L}^2(\rr^3)^4$ to $\mathit{L}^{2}(\S_\nu)^4$. As $\Phi^{z}_\nu(\Lambda^{z}_{+,\nu} )^{-1}$ is bounded from $\mathit{L}^{2}(\S_\nu)^4$ to $\mathit{L}^2(\rr^3)^4$, we get therefore that $\mathcal{T}_2$ and $\mathcal{T}_3$ are compact operators on $\mathit{L}^2(\rr^3)^4$. Now, let  $\chi:\S\rightarrow \rr$ be a $\mathit{C}^{\infty}$-smooth  and compactly  supported function on $\S$ and such that  $\S_\nu\setminus F\subsetneq \mathrm{supp}(\chi)$. Then we write  $ \mathcal{T}_1$ as follows 
\begin{align*}\label{}
    \mathcal{T}_1=&\left( \Phi^{z}_\nu \chi(\Lambda^{z}_{+,\nu} )^{-1} - \Phi^{z}\chi(\Lambda^{z}_{+} )^{-1} \right)  (\mathcal{H} -z)^{-1}\downharpoonright_{F}\\
    &+\chi\left( \Phi^{z}_\nu(1- \chi)(\Lambda^{z}_{+,\nu} )^{-1} - \Phi^{z}(1- \chi)(\Lambda^{z}_{+} )^{-1} \right)  (\mathcal{H} -z)^{-1}\downharpoonright_{F}\\
    &+ (1- \chi)\left( \Phi^{z}_\nu(1- \chi)(\Lambda^{z}_{+,\nu} )^{-1} - \Phi^{z}(1- \chi)(\Lambda^{z}_{+} )^{-1} \right)  (\mathcal{H} -z)^{-1}\downharpoonright_{F}:=  \mathcal{T}_4+ \mathcal{T}_5+\mathcal{T}_6.
\end{align*}
From the definition of $\Phi^{z}_\nu $, it is clear that $\mathcal{T}_6=0 $. Moreover, as $(\Lambda^{z}_{+,\nu} )^{-1}(\mathcal{H} -z)^{-1}\downharpoonright_{F} $ is bounded from $\mathit{L}^2(\rr^3)^4$ to $\mathit{H}^{1/2}(\S_\nu)^4$, using again the compactness of the Sobolev embedding, we obtain that $ \mathcal{T}_4$ and $ \mathcal{T}_5$ are also  compact operators on $\mathit{L}^2(\rr^3)^4$.   Therefore $\mathcal{T}$ is a compact operator in $\mathit{L}^2(\rr^3)^4$. Now, using Proposition \ref{BS-K}$(\mathrm{ii})$  it follows  that $\mathcal{T}= (\mathcal{H}^{\nu}_{\k} -z)^{-1}-(\mathcal{H}_{\k} -z)^{-1}$ is a compact operator in $\mathit{L}^2(\rr^3)^4$. Therefore,  by Weyl's theorem  we conclude that $\mathcal{H}^{\nu}_{\k}$ has the same essential spectrum as $\mathcal{H}_{\k}$. 
This finishes the proof of the theorem. \qed 
\\

As it was mentioned in the introduction,  in \cite{EL} the Schr\"{o}dinger operator with $\delta$-interactions (i.e the coupling $\Delta +\ep\delta_\S$ in $\rr^3$) have been considered, for a surface $\S$ satisfying the assumption $(\mathrm{H2})$.  In there,  the authors showed that for a fixed $\ep$ (such that $\mathrm{Sp}_{\mathrm{disc}}( \Delta +\ep\delta_\S)\neq\emptyset$) the discrete spectrum of $\Delta +\ep\delta_\S$ consists of exactly one simple eigenvalue for all sufficiently small $\nu>0$. Moreover, an asymptotic of this eigenvalue has been proved in terms of $\ep$, $\nu$ and $\phi$. Thus, it would be interesting to  investigate such a problem for the couplings  $\mathcal{H}+\ep\delta_\S$ and $\mathcal{H}+\beta\delta_\S$, and see if results of this type are valid.  We plan to study this problem on the future.

\subsection{Critical case}\label{Subs4}
From now, we assume that $\mathrm{sgn}(\k)=4$. The goal of this subsection is to prove the following result.  
\begin{theorem}\label{cas deforme} Let $\k=(\ep,\m,\e)\in\rr^3$ be such that $\mathrm{sgn}(\k)=4$ and let $\overline{\mathcal{H}_{\k}}$ be as in Theorem \ref{main1}.  If  $\S$ satisfies $(\mathrm{H2})$, then for all $\nu\geqslant0$ it holds that 
\begin{align}\label{EQ}
\mathrm{Sp}_{\mathrm{ess}}(\overline{\mathcal{H}^{\nu}_{\k}})=\big(-\infty,-m\big]\cup\left\{-\frac{m\m}{\ep}\right\}\cup \big[m,+\infty\big),
\end{align}
and the equality $\mathrm{Sp}(\overline{\mathcal{H}^0_{\k}})=\mathrm{Sp}_{\mathrm{ess}}(\overline{\mathcal{H}^0_{\k}})$ holds true (i.e. when $\nu=0$).
 
\end{theorem}
A few comments are in order.  Note that $\ep^2>\m^2$, thus the point $-m\m/\ep$ belongs to the gap $(-m,m)$.  Furthermore, one can imagine that the operator $\overline{\mathcal{H}_{\k}}$ is unitary equivalent to $\overline{\mathcal{H}_{\ep_1,\m_1}}$, for some $\ep_1,\m_1\in\rr$, such that $\ep^2_1-\m^2_1=4$ and $\ep_1/\ep=\m_1/\m$. Indeed, in \cite{Albert} and \cite{CLMT} it has been shown that  the potential $\eta(\alpha\cdot N)\delta_{\S}$ can always be absorbed as a change of gauge. So the existence of such a unitary transformation is not excluded. Another way to understand Theorem \ref{cas deforme} comes from the way in which we have presented the operator $\mathcal{H}_{\k}$. In fact, in this paper  we introduced the operator $\mathcal{H}_{\k}$ as the perturbation of the coupling $ \mathcal{H} +(\ep I_4 +\mu\beta )\delta_{\S}$ with the singular potential $ \eta(\alpha\cdot N)\delta_{\S}$. However, the right way is to say that $\mathcal{H}_{\k}$ is the perturbation of $ \mathcal{H} +\eta(\alpha\cdot N))\delta_{\S}$ with the singular potential $(\ep I_4 +\mu\beta )\delta_{\S}$, since  for all $\e\in\rr$, the operator $ \mathcal{H} +\eta(\alpha\cdot \mathit{N})\delta_{\S}$ is  self-adjoint (even if $\S$ is Lipschitz) and $\mathrm{Sp}( \mathcal{H} +\eta(\alpha\cdot \mathit{N})\delta_{\S})=\big(-\infty,-m\big]\cup\big[m,+\infty\big)$, cf.  \cite{BB}.
\newline

From Theorem \ref{cas deforme} we get a simple way to describe functions belonging to the domain of $\overline{\mathcal{H}_{\k}} $ when  $\S=\S_0$, i.e $\nu=0$. Indeed, we have the following result.
\begin{corollary} Assume that $\S :=\S_0$ and let $\overline{\mathcal{H}_{\k}}$ be as above. The following hold:
\begin{itemize}
 \item[(i)] If $\mu\neq0$, then 
  \begin{align}
  \mathrm{dom}(\overline{\mathcal{H}_{\k}})= \left\{ u+\Phi[-\tilde{\Lambda}^{-1}_+[t_{\S}u]]: u\in\mathit{H}^1(\rr^3)^4 \right\}.
  \end{align}
  \item[(ii)]   If $\mu=0$, then $\mathrm{dom}(\overline{\mathcal{H}_{\k}})=  \mathrm{dom}(\mathcal{H}_{\k})+\Phi[\mathrm{Kr}(\tilde{\Lambda}_+)]$.
    \end{itemize}
\end{corollary}
\textbf{Proof.} Assertion $(\mathrm{i})$ is a direct consequence of Theorem \ref{cas deforme} and Proposition \ref{BS-K}.  Assertion  $(\mathrm{ii})$ follows using the same arguments as those in \cite[Proposition 3.10]{AMV1}.\qed
 \newline
 
The main properties of the operators  $\overline{\mathcal{L}^{a}_{\pm}}$  which are relevant for us to prove Theorem \ref{cas deforme}  are collected in the following proposition.
\begin{proposition}\label{des equivalences} Let  $\k=(\ep,\m,\e)\in\rr^3$ be such that $\mathrm{sgn}(\k)\neq0$ and let $\mathcal{L}^{a}_{\pm,\k}:=\mathcal{L}^{a}_{\pm}$ be as in Lemma \ref{lemme ope}. Then, for all $a\in(-m,m)$, it holds that 
\begin{align*}
0\in\mathrm{Sp}_{\bullet}(\overline{\mathcal{L}^{a}_{+,\k}}) \Longleftrightarrow  0\in\mathrm{Sp}_{\bullet}(\overline{\mathcal{L}^{-a}_{+,\tilde{\k}}})  \Longleftrightarrow  0\in\mathrm{Sp}_{\bullet}(\overline{\mathcal{L}^{-a}_{-,\k}}),
\end{align*}
where $\tilde{\k}=(-\ep,\m,-\e)$ and $\bullet\in\{ \mathrm{ess},\mathrm{disc}\}$. In particular, $a\in\mathrm{Sp}(\overline{\mathcal{H}_{\k}})$ if and only if $-a\in\mathrm{Sp}(\overline{\mathcal{H}_{\tilde{\k}}})$.
\end{proposition}

\textbf{Proof.} Fix $\k=(\ep,\m,\e)\in\rr^3$ such that $\mathrm{sgn}(\k)\neq0$.  Following \cite[Proposition 4.2]{BEHL2}, for $f\in\mathit{L}^{2}(\S)^4$ we define 
\begin{align}\label{matt}
\mathcal{C}(f)=i\beta\alpha_2 \overline{f}^{\mathrm{c}}, \quad T(f)=\gamma_5\beta f, \quad \gamma_5 :=-i\alpha_1\alpha_2\alpha_3=\begin{pmatrix}
0& \mathit{I}_2 \\
\mathit{I}_2  & 0
\end{pmatrix}, 
\end{align}
where $\overline{f}^{\mathrm{c}}$ is the complex conjugate of $f$. Remark that  $\overline{\alpha_2}^{\mathrm{c}}=-\alpha_2$, $\gamma_5\beta=-\beta\gamma_5$ and $ \gamma_5(\alpha\cdot x)=  (\alpha\cdot x)\gamma_5$,  for all $x\in\rr^3$. Using this, it easily follows  that  $\mathcal{C}^2(f)=f$ and $T^2(f)=-f$. Moreover, a simple computation using the anticommutation relations of Dirac matrices yields that 
\begin{align}\label{des relation}
\Lambda^{\pm a}_{\pm,\k}[T(f)]= T(\Lambda^{\mp a}_{\mp,\k}[f]),\quad \quad \Lambda^{ a}_{+,\k}[\mathcal{C}(f)]= -\mathcal{C}(\Lambda^{- a}_{+,\tilde{k}}[f]),\quad \Lambda^{- a}_{+,\tilde{k}}[\mathcal{C}(f)]=-\mathcal{C}(\Lambda^{ a}_{+,\k}[f]).
\end{align}
Fix $a\in(-m,m)$ and assume that  $0\in\mathrm{Sp}_{\mathrm{ess}}(\overline{\mathcal{L}^{a}_{+}})$. Then, there exists a sequence of functions $(g_j)_{j\in\mathbb{N}}\subset\mathrm{dom}(\overline{\mathcal{L}^{a}_{+}})\subset\mathit{L}^{2}(\S)^4$, such that $ \left|\left| g_j \right|\right|_{\mathrm{L}^2(\S)^4}=1$, $(g_j)_{j\in\mathbb{N}}$ converges weakly to $0$ and $ \left|\left| \overline{\mathcal{L}^{a}_{+,\k}}g_j\right|\right|_{\mathrm{L}^2(\S)^4}\xrightarrow[j\to\infty]{}0$.  Hence, if we set $f_j=\mathcal{C}(g_j)$ and $h_j=T(g_j)$, then  it is clear that $(f_j)_{j\in\mathbb{N}}$ and $(h_j)_{j\in\mathbb{N}}$ converge weakly to zero and we have 
\begin{align*}
 \left|\left| h_j\right|\right|_{\mathrm{L}^2(\S)^4}=\left|\left| f_j\right|\right|_{\mathrm{L}^2(\S)^4}=1,\quad f_j\in\mathrm{dom}(\overline{\mathcal{L}^{-a}_{+,\tilde{\k}}}) \text{ and } h_j\in\mathrm{dom}(\overline{\mathcal{L}^{-a}_{-,\k}}),\quad  \forall j\in\mathbb{N}.
 \end{align*}
  Now using \eqref{des relation} it follows that 
\begin{align*}
\left|\left| \overline{\mathcal{L}^{-a}_{+,\tilde{\k}}}[f_j]\right|\right|_{\mathrm{L}^2(\S)^4} = \left|\left| \overline{\mathcal{L}^{-a}_{-,\k}}[h_j]\right|\right|_{\mathrm{L}^2(\S)^4}= \left|\left| \overline{\mathcal{L}^{a}_{+,\k}}[g_j]\right|\right|_{\mathrm{L}^2(\S)^4}. 
\end{align*}
 Therefore $0\in\mathrm{Sp}(\overline{\mathcal{L}^{-a}_{+,\tilde{\k}}})$ and  $0\in\mathrm{Sp}(\overline{\mathcal{L}^{-a}_{-,\k}})$. The reverse implications follow in the same way. Now that $0\in\mathrm{Sp}_{\mathrm{disc}}(\overline{\mathcal{L}^{a}_{+}}) \Longleftrightarrow  0\in\mathrm{Sp}_{\mathrm{disc}}(\overline{\mathcal{L}^{-a}_{+,\tilde{\k}}})  \Longleftrightarrow  0\in\mathrm{Sp}_{\mathrm{disc}}(\overline{\mathcal{L}^{-a}_{-,\k}})$ is a direct consequence of \eqref{des relation},  and this finishes the proof of the first statement. The last statement is a direct consequence of the first one and Corollary \ref{corollaire inter}. This completes the proof.\qed
\newline
\begin{proposition}\label{essspectra} Let $a\in(-m,m)$  and let $\overline{\mathcal{L}^{a}_{\pm}}$ be as in Lemma \ref{lemme ope}.  Assume that $\nu=0$, then it holds that 
\begin{align}
0\in\mathrm{Sp}(\overline{\mathcal{L}^{a}_{+}}) \Longleftrightarrow a=-\frac{m\m}{ \ep}\quad\text{ and }\quad 0\in\mathrm{Sp}(\overline{\mathcal{L}^{a}_{-}}) \Longleftrightarrow a=\frac{m\m}{ \ep}.
\end{align}
Moreover, $0$ is an isolated eigenvalue of $\overline{\mathcal{L}^{- m\m/\ep}_{+}}$ and $\overline{\mathcal{L}^{ m\m/\ep}_{-}}$ with infinite multiplicity.\end{proposition}
\textbf{Proof.} Given $a\in(-m,m)$, once the claimed statement is shown for $\overline{\mathcal{L}^{a}_{+}}$, by Proposition \ref{des equivalences} we get the result for  $\overline{\mathcal{L}^{a}_{-}}$.  As in Theorem \ref{Prop non}, on the Fourier side, if we let $\langle \xi\rangle:=(1+|\xi|^2)^{1/2}$ then one can check that $\overline{\mathcal{L}^{a}_{+}}$ is unitary equivalent to the following multiplication operator:
 \begin{align}
 \widetilde{\Pi^{a}_{+}}:= \langle \xi\rangle\left( \frac{1}{\mathrm{sgn}(\k)}(\ep I_4 -(\mu\beta+\e(\alpha\cdot N)) )+ \frac{1}{2\sqrt{|\xi|^2+m^2-a^2}}\Gamma_{ m, a}(\xi)\right).
\end{align}
Since $\mathrm{sgn}(\k)=4$, from \eqref{inverse} it follows that $\widetilde{\Pi^{a}_{+}}$ is invertible for all $a\neq -m\m/ \ep$, and we have 
\begin{align*} \label{inverse2}
(\widetilde{\Pi^{a}_{+}})^{-1}&=  \frac{1}{\langle \xi\rangle}\left( 1 + \frac{\sqrt{|\xi|^2+m^2-a^2}}{\ep a+\m m} -\frac{(\ep+(\m\beta+\e(\alpha\cdot N)))}{2(\ep a+ \m m)}\Gamma_{ m, a}(\xi)\right)(\ep+(\m\beta+\e(\alpha\cdot N))).
\end{align*}
Furthermore it holds that 
\begin{align*}
 \frac{1}{\langle \xi\rangle}\widetilde{\Pi^{a}_{+}}\left( 1 -\frac{(\ep+\m\beta+\e(\alpha\cdot N))}{2\sqrt{|\xi|^2+m^2-a^2}}\Gamma_{ m, a}(\xi)\right)&=0, \text{ for } a=-\frac{m\m}{ \ep}.
\end{align*}
From this, it follows that $0$ is an eigenvalue of the operators $\overline{\mathcal{L}^{- m\m/\ep}_{+}}$   with infinite multiplicity, and thereby  $0\in\mathrm{Sp}_{\mathrm{ess}}(\overline{\mathcal{L}^{- m\m/\ep}_{+}})$.
Thus, we conclude that $0\in\mathrm{Sp}(\overline{\mathcal{L}^{a}_{+}})$ if and only if $a\ep=-m\m$. Now we turn to prove the last statement for the operator $\overline{\mathcal{L}^{- m\m/\ep}_{+}}$, similar arguments give the result for  $\overline{\mathcal{L}^{ m\m/\ep}_{-}}$.  A simple computation yields
 \begin{align}
\mathrm{det}(\widetilde{\Pi^{a}_{+}}-\theta)= \bigg[\theta\bigg(\theta-\underbrace{\langle \xi\rangle\bigg( \frac{a}{\sqrt{|\xi|^2+m^2-a^2}} +\frac{\ep}{2}\bigg)}_{\theta_1}\bigg)\bigg]^2,
\end{align}
where $\mathrm{det}(\widetilde{\Pi^{a}_{+}}-\theta)$ is the determinant of $(\widetilde{\Pi^{a}_{+}}-\theta)$. By studying the variations of the non-trivial root $\theta_1$, we obtain that
\begin{align*}\label{}
\mathrm{Sp}(\overline{\mathcal{L}^{- m\m/\ep}_{+}})&=\{0\} \cup \theta_1([0,\infty))=\left\{0\right\} \cup \left[\frac{\ep}{2} - \frac{\m}{\sqrt{\ep^2-\m^2}},\infty\right] \quad \text{ if } \ep>0,\\
 \mathrm{Sp}(\overline{\mathcal{L}^{- m\m/\ep}_{+}})&=\theta_1([0,\infty))\cup \{0\} = \left[-\infty,\frac{\ep}{2} +\frac{\m}{\sqrt{\ep^2-\m^2}}\right] \cup\{0\} \quad \text{ if } \ep<0.
\end{align*}
Take into account that $\mathrm{sgn}(\k)=4$, we then  get  that $0$  is an isolated eigenvalue of  $\overline{\mathcal{L}^{- m\m/\ep}_{+}}$ with infinite multiplicity, which completes the proof of the proposition.\qed
\newline

\begin{remark} The reader should not confuse the unbounded  operator  $\mathcal{L}^{-m\m/\ep}_{+}$ with the original operator $\Lambda^{-m\m/\ep}_{+}$, which is indeed a bounded operator on $\mathit{L}^2(\S)^4$ with closed range.
\end{remark}
We are now in a position to complete the proof of our main result in this subsection.

\textbf{Proof of Theorem \ref{cas deforme}.} 
  Assume that $\S$ satisfies $(\mathrm{H2})$ and fix $\nu\geqslant0$. The result will follow from the following statements:
\begin{itemize}
  \item[(a)]  $\big(-\infty,-m\big)\cup \big(m,+\infty\big)\subset\mathrm{Sp}_{\mathrm{ess}}(\overline{\mathcal{H}^{\nu}_{\k}})$.
  \item[(b)] $\left\{-m\m/\ep\right\}\in\mathrm{Sp}_{\mathrm{ess}}(\overline{\mathcal{H}^{\nu}_{\k}})$ and $\left\{m\m/\ep \right\}\notin\mathrm{Sp}_{\mathrm{ess}}(\overline{\mathcal{H}^{\nu}_{\k}})$.
  \item [(c)]  $\mathrm{Sp}_{\mathrm{ess}}(\overline{\mathcal{H}^{\nu}_{\k}})\cap\left[(-m,m)\setminus  \{-m\mu / \ep, m\mu / \ep \}\right]=\emptyset$.
\end{itemize}

  We are going to show $(\mathrm{a})$. For that,  given $a\in(-\infty,-m)\cup(m,\infty)$ and let $(\varphi_{n})_{n\in\mathbb{N}}$ be the sequence of functions defined by \eqref{suitee} with $R=2\sup\{|x|: x\in\S_\nu\setminus \overline{F}\}$.  By construction, it is clear that $(\varphi_{n})_{n\in\mathbb{N}}$ is a singular sequence for $\overline{\mathcal{H}^{\nu}_{\k}}$ and $a$. Therefore we get the inclusion $(-\infty,-m)\cup (m,+\infty)\subset \mathrm{Sp}_{\mathrm{ess}}(\overline{\mathcal{H}^{\nu}_{\k}})$, which yields $(\mathrm{a})$.

   Now,  we turn  to the proof of $(\mathrm{b})$. Actually from Proposition \ref{essspectra} and Corollary \ref{corollaire inter}, we know that item $(\mathrm{b})$ holds true for $\nu=0$.   Next, assume that $\nu>0$, we are going to prove that $\left\{-m\m/\ep\right\}\in\mathrm{Sp}_{\mathrm{ess}}(\overline{\mathcal{H}^{\nu}_{\k}})$ and the same arguments yield that  $\left\{m\m/\ep \right\}\notin\mathrm{Sp}_{\mathrm{ess}}(\overline{\mathcal{H}^{\nu}_{\k}})$. Assume that   $\left\{-m\m/\ep\right\}\notin\mathrm{Sp}_{\mathrm{ess}}(\overline{\mathcal{H}^{\nu}_{\k}})$, then by Corollary \ref{corollaire inter} and Proposition \ref{des equivalences}  it follows that  $0\notin\mathrm{Sp}_{\mathrm{ess}}(\overline{\mathcal{L}^{-m\m / \ep}_{+,\nu}}) $ and $0\notin\mathrm{Sp}_{\mathrm{ess}}(\overline{\mathcal{L}^{m\m / \ep}_{-,\nu}}) $. 
   
   Now we set $B_{\nu}:=\tilde{\Lambda}^{-m\m/\ep}_{+,\nu}\tilde{\Lambda}^{m\m/\ep}_{-,\nu}$, $D_\nu:= \Lambda^{m\m/\ep}_{-,\nu}\Lambda^{-m\m/\ep}_{+,\nu}$ and we consider the operator  $\Upsilon^{-m\m/\ep}_\nu:\mathit{L}^2(\S_\nu)^4\longrightarrow \mathit{L}^2(\S_\nu)^4$ defined by: 
  \begin{align}
  \Upsilon^{-m\m/\ep}_\nu:=L^{\frac{1}{4}}_{\nu}D_{\nu} B_{\nu}L^{\frac{1}{4}}_{\nu}=L^{\frac{1}{4}}_{\nu}(\Lambda^{m\m/\ep}_{-,\nu}\Lambda^{-m\m/\ep}_{+,\nu})(\tilde{\Lambda}^{-m\m/\ep}_{+,\nu}\tilde{\Lambda}^{m\m/\ep}_{-,\nu}) L^{\frac{1}{4}}_{\nu}.
  \end{align}
   Observe that 
  \begin{align} 
  \begin{split}
 B_{\nu}=&\mathit{C}^{-m\m/\ep}_{\S} (\alpha\cdot \mathit{N})\lbrace \alpha\cdot \mathit{N}, \tilde{\mathit{C}}^{-m\m/\ep}_{\S}\rbrace -2\frac{m\m}{\ep}\mathit{C}^{-m\m/\ep}_{\S}S^{-m\m/\ep}\\
  &+\frac{m\m\e}{2\ep}(\alpha\cdot \mathit{N})S^{-m\m/\ep} + \frac{\e}{4} \lbrace \alpha\cdot \mathit{N}, \tilde{\mathit{C}}^{-m\m/\ep}_{\S}\rbrace,
\end{split}
  \end{align}
  and 
   \begin{align} 
    \begin{split}
 D_{\nu}=&\mathit{C}^{-m\m/\ep}_{\S} (\alpha\cdot \mathit{N})\lbrace \alpha\cdot \mathit{N}, \mathit{C}^{-m\m/\ep}_{\S}\rbrace -2\frac{m\m}{\ep}S^{-m\m/\ep}\mathit{C}^{-m\m/\ep}_{\S}\\
  &+\frac{m\m\e}{2\ep}(\alpha\cdot \mathit{N})S^{-m\m/\ep}+ \frac{\e}{4} \lbrace \alpha\cdot \mathit{N}, \mathit{C}^{-m\m/\ep}_{\S}\rbrace.
\end{split}
  \end{align}
 As $L^{\frac{1}{4}}_{\nu}$ is an isomorphism, using Lemma \ref{commutator} it easily follows  that  $\Upsilon^{-m\m/\ep}_\nu$ is a bounded,  self-adjoint operator on $\mathit{L}^2(\S_\nu)^4$.  Note that by hypothesis, it holds that   $0\notin\mathrm{Sp}_{\mathrm{ess}}(\Upsilon^{-m\m/\ep}_\nu)\big)$.  Next, we introduce the unitary transformation  $U:\mathit{L}^2(\S_\nu)^4\longrightarrow \mathit{L}^2(\rr^2)^4$ defined by $Ug(x)= J^{1/2}_\nu(\overline{x})g(\tau(\overline{x}))$,  where $J^{1/2}_\nu$ is given by \eqref{leJ}. We claim that $U\Upsilon^{-m\m/\ep}_\nu U^{-1}-\Upsilon^{-m\m/\ep}_0$ is a compact operator on $\mathit{L}^2(\S_\nu)^4$. Indeed, let $\chi:\S\rightarrow \rr$ be a $\mathit{C}^{\infty}$-smooth  and compactly  supported function on $\S$ and such that  $\S_\nu\setminus F\subsetneq \mathrm{supp}(\chi)$,  we then get 
 \begin{align*}
  U\Upsilon^{-m\m/\ep}_\nu U^{-1}-\Upsilon^{-m\m/\ep}_0:=  ( U\Upsilon^{-m\m/\ep}_\nu U^{-1}-\Upsilon^{-m\m/\ep}_0)\chi +( U\Upsilon^{-m\m/\ep}_\nu U^{-1}-\Upsilon^{-m\m/\ep}_0)(1-\chi).
  \end{align*}
   Since the embedding $\chi\mathit{L}^{2}(\rr^2)^4\hookrightarrow\mathit{H}^{-1/2}(\rr^2)^4$ is compact, and $ U\Upsilon^{-m\m/\ep}_\nu U^{-1}$ is bounded from $\mathit{H}^{-1/2}(\rr^2)^4$ to $\mathit{L}^{2}(\rr^2)^4$, for all $\nu\geqslant0$, it follows that $( U\Upsilon^{-m\m/\ep}_\nu U^{-1}-\Upsilon^{-m\m/\ep}_0)\chi $ is a compact operator on $\mathit{L}^{2}(\rr^2)^4$. Next, observe  that $ ( U\Upsilon^{-m\m/\ep}_\nu U^{-1}-\Upsilon^{-m\m/\ep}_0)(1-\chi) :=   T_1+T_2$, where
 \begin{align*} 
  T_1&= UL^{\frac{1}{4}}_{\nu}D_{\nu}B_{\nu}\chi L^{\frac{1}{4}}_{\nu}U^{-1}(1-\chi)- L^{\frac{1}{4}}_{0}D_{0}B_{0}\chi L^{\frac{1}{4}}_{0}(1-\chi),\\
  T_2&= \left(UL^{\frac{1}{4}}_{\nu}D_{\nu}B_{\nu}- L^{\frac{1}{4}}_{0}D_{0}B_{0}\right)(1-\chi) L^{\frac{1}{4}}_{0}(1-\chi).
   \end{align*}
Recall that  $L^{\frac{1}{4}}_{\nu}U^{-1}$ is bounded from $\mathit{L}^2(\rr^2)^4$ to $\mathit{H}^{-1/2}(\S_\nu)^4$, for all $v\geqslant0$, and the embedding  $\chi\mathit{H}^{-1/2}(\S_\nu)^4\hookrightarrow\mathit{H}^{-1}(\S_\nu)^4$ is compact. Since $B_{\nu}$ is bounded from $\mathit{H}^{-1}(\S_\nu)^4$ to $\mathit{L}^{2}(\S_\nu)^4$ (see Remark \ref{ext remark}) it follows that  $B_{\nu}\chi L^{\frac{1}{4}}_{\nu}U^{-1}$ is a compact operator from $\mathit{L}^2(\rr^2)^4$  to $\mathit{L}^2(\S_\nu)^4$. Now, as $UL^{\frac{1}{4}}_{\nu}D_{\nu}$ is bounded from $\mathit{L}^2(\S_\nu)^4$ to $\mathit{L}^2(\rr^2)^4$,  we then get that $T_1$ is a compact operator on $\mathit{L}^2(\rr^2)^4$.
We now apply this argument again to the operator $T_2$, and we get
\begin{align*} 
  T_2=& \left(UL^{\frac{1}{4}}_{\nu}D_{\nu}\chi B_{\nu}- L^{\frac{1}{4}}_{0}D_{0}\chi B_0\right)(1-\chi) L^{\frac{1}{4}}_{0}(1-\chi)\\
  &+ \bigg(UL^{\frac{1}{4}}_{\nu}D_{\nu}(1-\chi)B_{\nu}- L^{\frac{1}{4}}_{0}D_0(1-\chi)B_0\bigg)(1-\chi) L^{\frac{1}{4}}_{0}(1-\chi):=T_3+T_4.
   \end{align*}
Then in the same manner, it is easy to check that $T_3$ is a compact operator on $\mathit{L}^2(\rr^2)^4$. Now set 
\begin{align*} 
\tilde{B_\nu}:=(1-\chi)B_\nu(1-\chi)=(1-\chi)(\tilde{\Lambda}^{-m\m/\ep}_{+,\nu}\tilde{\Lambda}^{m\m/\ep}_{-,\nu})(1-\chi),
 \end{align*}
Observe that 
\begin{align*} 
\tilde{B_\nu}=&\bigg(\mathit{C}^{-m\m/\ep}_{\S} \chi(\alpha\cdot \mathit{N})\left\{ \alpha\cdot \mathit{N}, \tilde{\mathit{C}}^{-m\m/\ep}_{\S}\right\} -2\frac{m\m}{\ep}\mathit{C}^{-m\m/\ep}_{\S}\chi S^{-m\m/\ep}\\
&+\frac{m\m\e}{2\ep}\chi(\alpha\cdot \mathit{N})S^{-m\m/\ep} 
+\frac{\e}{4} \chi \left\{ \alpha\cdot \mathit{N}, \tilde{\mathit{C}}^{-m\m/\ep}_{\S}\right\} \\
   &-2\frac{m\m}{\ep}\chi\mathit{C}^{-m\m/\ep}_{\S}(1-\chi) S^{-m\m/\ep} -2\chi\frac{m\m}{\ep}\mathit{C}^{-m\m/\ep}_{\S}(1-\chi) S^{-m\m/\ep}\bigg)(1-\chi)\\
  &+(1-\chi)\bigg( -2\frac{m\m}{\ep}\mathit{C}^{-m\m/\ep}_{\S}(1-\chi) S^{-m\m/\ep}+\frac{m\m\e}{2\ep}(\alpha\cdot \mathit{N})S^{-m\m/\ep}\bigg)(1-\chi):=B_{\nu,1}+B_{2}.
   \end{align*}
 Here we used the fact that $(1+ \chi)(\alpha\cdot \mathit{N})\lbrace \alpha\cdot \mathit{N}, \tilde{\mathit{C}}^{-m\m/\ep}_{\S}\rbrace (1+ \chi)$ vanishes identically. Therefore we obtain that 
 \begin{align*}
 T_4=& \bigg(UL^{\frac{1}{4}}_{\nu}D_{\nu}B_{\nu,1}- L^{\frac{1}{4}}_{0}D_0B_{0,1}\bigg) L^{\frac{1}{4}}_{0}(1-\chi)+\bigg(UL^{\frac{1}{4}}_{\nu}D_{\nu}- L^{\frac{1}{4}}_{0}D_{0}\bigg)B_2 L^{\frac{1}{4}}_{0}(1-\chi):=T_5+T_6.
   \end{align*}
   Again,  using the compactness of the Sobolev injection, one can show that $T_5$ is a compact operator on $\mathit{L}^2(\rr^2)^4$. Next, remark that 
\begin{align*} 
D_\nu(1-\chi)=&\bigg(\mathit{C}^{-m\m/\ep}_{\S}\chi (\alpha\cdot \mathit{N})\left\{ \alpha\cdot \mathit{N}, \mathit{C}^{-m\m/\ep}_{\S}\right\} -2\frac{m\m}{\ep}\chi S^{-m\m/\ep}\mathit{C}^{-m\m/\ep}_{\S}\\
&+\frac{m\m\e}{2\ep}(\alpha\cdot \mathit{N})S^{-m\m/\ep} + \frac{\e}{4} \chi\left\{ \alpha\cdot \mathit{N}, \mathit{C}^{-m\m/\ep}_{\S}\right\}\bigg)(1-\chi) \\
  &+(1-\chi) \bigg( -2\frac{m\m}{\ep}S^{-m\m/\ep}(1-\chi) \mathit{C}^{-m\m/\ep}_{\S}+\frac{m\m\e}{2\ep}(\alpha\cdot \mathit{N})S^{-m\m/\ep}\bigg)(1-\chi)\\
  &:= D_{\nu,1}(1-\chi)+D_2(1-\chi).
   \end{align*}
Note that  $(U(1-\chi) L^{\frac{1}{4}}_{\nu}-(1-\chi) L^{\frac{1}{4}}_{0})D_2B_2 L^{\frac{1}{4}}_{0}(1-\chi)=0$. Therefore, we obtain that  
\begin{align*}
 T_6=\bigg(UL^{\frac{1}{4}}_{\nu}D_{\nu,1}- L^{\frac{1}{4}}_{0}D_{0,1}\bigg)B_2 L^{\frac{1}{4}}_{0}(1-\chi)+ \bigg(U\chi L^{\frac{1}{4}}_{\nu}- \chi L^{\frac{1}{4}}_{0}\bigg)D_2 B_2 L^{\frac{1}{4}}_{0}(1-\chi).
   \end{align*}
Then, we conclude as above that $T_6$  is a compact operator on $\mathit{L}^2(\rr^2)^4$. Therefore,  $ U\Upsilon^{-m\m/\ep}_\nu U^{-1}-\Upsilon^{-m\m/\ep}_0$ is compact in $\mathit{L}^2(\rr^2)^4$. As $0\in\mathrm{Sp}_{\mathrm{ess}}\left(\Upsilon^{-m\m/\ep}_0\right)$ because $0\in\mathrm{Sp}_{\mathrm{ess}}\left(\overline{\mathcal{L}^{ m\m/\ep}_{-}}\right)$ by  Proposition \ref{essspectra},  then by Weyl's theorem we get that $0\in\mathrm{Sp}_{\mathrm{ess}}\left(\Upsilon^{-m\m/\ep}_\nu\right)$. This contradicts the fact that $0\notin\mathrm{Sp}_{\mathrm{ess}}(\Upsilon^{-m\m/\ep}_\nu)\big)$, which proves $(\mathrm{b})$.

 We now  show $(\mathrm{c})$, so assume that $a\in(-m,m)\setminus\{  -m\m/ \ep,m\m/ \ep\}$. We introduce the operator $G^a_\nu:\mathit{L}^2(\S_\nu)^4\longrightarrow \mathit{L}^2(\S_\nu)^4$ defined by: 
$$G^a_\nu:=L^{\frac{1}{4}}_{\nu}(\Lambda^{a}_{-,\nu}\Lambda^{a}_{+,\nu})(\tilde{\Lambda}^{a}_{-,\nu}\tilde{\Lambda}^{a}_{+,\nu}) L^{\frac{1}{4}}_{\nu}.$$
 Clearly, $G^a_\nu$ is bounded self-adjoint in $\mathit{L}^2(\S_\nu)^4$, since $\tilde{\Lambda}^{a}_{-,\nu}\tilde{\Lambda}^{a}_{+,\nu}=\tilde{\Lambda}^{a}_{+,\nu}\tilde{\Lambda}^{a}_{-,\nu}$. Moreover, by definition we have
 \begin{align}\label{comparaison} 
 0\in\mathrm{Sp}_{\mathrm{ess}}(\overline{\mathcal{L}^{a}_{\pm,\nu}}) \Longrightarrow  0\in\mathrm{Sp}_{\mathrm{ess}}(G^a_\nu).
 \end{align}
As  $\tilde{\Lambda}^{a}_{+,0}$ and $\tilde{\Lambda}^{a}_{-,0}$ are bounded, invertible operators for all $a\in(-m,m)\setminus\{ -m\m/ \ep,m\m/ \ep\}$, from Proposition \ref{essspectra} it follows that  $0\in\mathrm{Sp}_{\mathrm{ess}}(G^a_0) $ if and only if $a=\mp m\m/ \ep$. Next, we claim that: if $a\neq\mp m\m/ \ep$, then  $0\notin\mathrm{Sp}_{\mathrm{ess}}(G^a_\nu) $. To prove this, let    $U$ be the unitary transformation defined in the proof of $(\mathrm{b})$, and set $T= UG^a_\nu U^{-1}- G^a_0$. Then,  $T$ is a compact operator in $\mathit{L}^{2}(\rr^2)^4$. Indeed, this may be handled in much the same way as in the proof of the  previous statement, we omit the details. Therefore  $\mathrm{Sp}_{\mathrm{ess}}(G^a_\nu)=\mathrm{Sp}_{\mathrm{ess}}(G^a_0) $ holds by Weyl's theorem. This proves the claim because  $0\in\mathrm{Sp}_{\mathrm{ess}}(G^a_0) $ if and only if $a=\mp m\m/ \ep$. Using this, from \eqref{comparaison} it follows that, if $a\neq\mp m\m/ \ep$ then $ 0\notin\mathrm{Sp}_{\mathrm{ess}}(\overline{\mathcal{L}^{a}_{\pm,\nu}})$. Therefore,  Corollary \ref{corollaire inter} yields that $\mathrm{Sp}_{\mathrm{ess}}(\overline{\mathcal{H}^{\nu}_{\k}})\cap\left[(-m,m)\setminus  \{-m\mu / \ep, m\mu / \ep \}\right]=\emptyset$, which proves $(\mathrm{c})$.

Summing up, from $(\mathrm{a})$ and $(\mathrm{b})$  we obtain that $\big(-\infty,-m\big)\cup\left\{-m\m/\ep\right\}\cup \big(m,+\infty\big)\subset\mathrm{Sp}_{\mathrm{ess}}(\overline{\mathcal{H}^{\nu}_{\k}})$. From  $(\mathrm{b})$ and $(\mathrm{c})$ we get the inclusion $\mathrm{Sp}_{\mathrm{ess}}(\overline{\mathcal{H}^{\nu}_{\k}}) \subset\big(-\infty,-m\big]\cup\left\{-m\m/\ep \right\}\cup \big[m,+\infty\big)$. Since the essential spectrum of a self-adjoint operator is closed, we get then the equality \eqref{EQ}. This completes the proof of the theorem. \qed
\newline

 Actually in the case $\S=\rr^2\times\{ 0\}$, one can check directly using the separation of variables that $a=-m\m/\ep$ is an eigenvalue of $ \overline{\mathcal{H}_{\k}}$  with infinite multiplicity. Indeed, let $a=-m\m/\ep$ and $\varphi \in \mathrm{dom}(\overline{\mathcal{H}_{\k}})$ such that: 
\begin{align}\label{Ev}
  ( \overline{\mathcal{H}_{\k}} -a)\varphi=0, \quad\text{ in }  \mathit{L}^{2}(\rr^3)^4.
  \end{align}
A simple computation yields  the following relations
 \begin{align}\label{a4}
 \begin{split}
  \left[\frac{1}{2}(\ep I_4-\mu\beta+ \eta\alpha_3) +i \alpha_3\right] \left[\frac{1}{2}(\ep I_4+\mu\beta- \eta\alpha_3) -i \alpha_3\right]&=(2-i\e)I_4,\\
   \left[\frac{1}{2}(\ep I_4-\mu\beta+ \eta\alpha_3) +i \alpha_3\right] \left[\frac{1}{2}(\ep I_4+\mu\beta- \eta\alpha_3) -i \alpha_3\right]&=i\alpha_3(\ep +\m\beta).
   \end{split}
 \end{align}
Hence, using this relation and the Definition \ref{def2},  another way of stating \eqref{Ev} is to say:
\begin{eqnarray}\label{3}
 \left\{
    \begin{split}\
     (\mathcal{H}-a)\varphi&= 0 \quad\text{for all } x_3\neq 0,\\
 (2-i\e)t_{\Sigma}\varphi_+&= -i\alpha_3(\ep +\m\beta)t_{\Sigma}\varphi_- \quad\text{ for }x_3=0.   
   \end{split}
  \right.  
\end{eqnarray}
  Since $(\mathcal{H}+a)(\mathcal{H}-a) = (-\D +m^2 -a^2)I_4$, one get that $\varphi$ is also solution of the following equation 
$$ (-\D +m^2 -a^2)I_4\varphi=0,\quad \text{  for all }  x_3\neq0$$
Thus, applying Fourier-Plancherel operator on $\overline{x}=(x_1,x_2)$, we get that 
\begin{equation}\label{11} 
   \mathcal{F}_{\overline{x}}\left[\varphi\right](\xi,x_3)=\left\{
    \begin{aligned}
      e^{-x_3\sqrt{| \xi|^2 +m^2-a^2}}\mathcal{F}_{\overline{x}}\left[\psi_+\right](\xi) \quad\text{for  } x_3>& 0,\\ 
      e^{x_3\sqrt{| \xi|^2 +m^2-a^2}}\mathcal{F}_{\overline{x}}\left[\psi_-\right](\xi)  \quad\text{for  } x_3< &0,   
      \end{aligned}
  \right.
\end{equation}
 for some $\psi_{\pm}\in \mathit{H}^{\frac{1}{2}}(\rr^2)$. Since $(\mathcal{H}+a)\varphi = 2a\varphi$, by applying the inverse Fourier-Plancherel operator,  we obtain that
 \begin{equation}\label{12} 
 \begin{split}
  \varphi(\overline{x},x_3)=\left\{
    \begin{aligned}
      \frac{1}{2\pi}&\int_{\rr^2}e^{i\overline{x}\cdot\xi} e^{-x_3\sqrt{| \xi|^2 +m^2-a^2}}\Gamma_{ i}(\xi)\mathcal{F}_{\overline{x}}\left[\psi_+\right](\xi)\mathrm{d}\xi \quad\text{for  } x_3> 0,\\ 
       \frac{1}{2\pi}&\int_{\rr^2}e^{i\overline{x}\cdot\xi} e^{x_3\sqrt{| \xi|^2 +m^2-a^2}}\Gamma_{ -i}\mathcal{F}_{\overline{x}}\left[\psi_-\right](\xi)\mathrm{d}\xi \quad\text{for  } x_3< 0,   
      \end{aligned}
  \right.
  \end{split}
\end{equation}
where  $\Gamma_{\pm i}(\xi)=\left[ \alpha\cdot(\xi_1,\xi_2,\pm i\sqrt{| \xi|^2 +m^2-a^2})+m\beta + a \right]$. From this,  it is clear that 
$$\varphi_{\pm},(\alpha\cdot\nabla)\varphi_{\pm}\in\mathit{L}^2(\O_{\pm}).$$ 
Now, if we set $$\psi_-= -\frac{\e-i2}{\ep^2-\m^2}(\ep+\m\beta)\alpha_3\psi_+,$$
 then we get 
 \begin{align}\label{4}
 \begin{split}
(2-i\e)\Gamma_{ +i}(\xi) \mathcal{F}_{\overline{x}}\left[\psi_+\right](\xi)=-i\alpha_3(\ep I_4 -\mu\beta)\Gamma_{ -i}(\xi) \mathcal{F}_{\overline{x}}\left[\psi_-\right](\xi).
\end{split}
\end{align}
Thus $(\varphi_+,\varphi_-)$ satisfies the transmission condition. Therefore, for all $\psi \in \mathit{H}^{\frac{1}{2}}(\rr^2)$ the function 
 \begin{equation}\label{12} 
 \begin{split}
  \varphi(\overline{x},x_3)=\left\{
    \begin{aligned}
      \frac{1}{2\pi}&\int_{\rr^2}e^{i\overline{x}\cdot\xi} e^{-x_3\sqrt{| \xi|^2 +m^2-a^2}}\Gamma_{ i}(\xi)\mathcal{F}_{\overline{x}}\left[\psi_+\right](\xi) \quad\text{for  } x_3> 0,\\ 
       \frac{-i}{4\pi}&\int_{\rr^2}e^{i\overline{x}\cdot\xi} e^{x_3\sqrt{| \xi|^2 +m^2-a^2}}\Gamma_{ -i}(\l-\m\beta)\alpha_3\mathcal{F}_{\overline{x}}\left[\psi_+\right](\xi) \quad\text{for  } x_3< 0,   
      \end{aligned}
  \right.
  \end{split}
\end{equation}
is an eigenvector associated to the eigenvalue $a=-m\m/\ep$.

\section{Quantum Confinement induced by Dirac operators with  anomalous magnetic $\delta$-shell interactions}\label{sec6}

\setcounter{equation}{0}
 The main goal of this section is to derive a new model of Dirac operators with $\delta$-shell interactions which generate confinement.  Let us explain how to derive this model. Using the unit $ c=\hbar=1$, where $c$ is the speed of light and $\hbar$ is the Planck's constant,  the Dirac operator with an electromagnetic field is given by (see \cite{Tha}):
 \begin{align}
\tilde{\mathcal{H}} =  \alpha \cdot(- i \nabla  -eA(x))+ m\beta+ e\phi_{el}(x)\mathit{I}_4,
\end{align}  
where $e$ is the charge of the particle, $\phi_{el}(x)$ is the electric field, and $A(x)$ is the magnetic vector potential.  Here the electric and magnetic field strengths are
 \begin{align*}
E(x)=  -\nabla\phi_{el}(x)  -\frac{\partial A(x)}{\partial t},\quad B(x)=\nabla \times A(x),
\end{align*}  
where $\partial/\partial t$ denotes the partial derivative with respect to time $t\in\rr$. In this setting, the anomalous magnetic potential is given by:
 \begin{align}
V(x)=  \u\left(i\beta(\alpha\cdot E(x))  - \frac{1}{4}\beta ((\alpha\times \alpha)\cdot B(x))\right).
\end{align}  
here the coupling constant $\u$ is the magnitude of the anomalous potential. Now, we put $ \phi_{el}(x)= |x|$ and $A(x)= 0$, we then obtain
 \begin{align*}
V(x)=  i\u\beta\left(\alpha\cdot \frac{x}{|x|}\right).
\end{align*} 
 Now, given $R>0$, if $x\in\mathbb{S}^{2}_{R}=\{ x\in\rr^3: |x|=R\}$, then $x/|x|$ coincide with the normal vector field $\mathit{N}(x)$. Thus we get
 \begin{align}
V_\u(x):=V(x)=  i\u\beta(\alpha\cdot \mathit{N}(x)).
\end{align} 
 Now,  given a surface $\S\subset\rr^3$  satisfying the assumption $(\mathrm{H1})$, we can consider the following Dirac operator
 \begin{align}
  \mathcal{H} + V_{\u}=  \mathcal{H}+ V_\u\delta_{\S}, \quad \u\in\rr.
\end{align} 
 and called it  Dirac operator with  anomalous magnetic $\delta$-shell interactions of strength $\u$. As we already mentioned in the introduction,  when finalizing the current manuscript, it turns out that the authors of the article \cite{CLMT} considered this problem in dimension two. However, instead of deriving the potential $V_\u$ as we had done here, they rigorously proved that the two-dimensional analog of $V_\u$ can be approximated  by regular shrinking potentials of magnetic type, and hence they justified the fact that $V_\u$ is a "magnetic" $\delta$-shell interactions. 
 \begin{remark} If one choose the magnetic field $A(x)$ so that $B(x)=x/|x|$ and put $ \phi_{el}(x)=0$, we then get the $\delta$-potential  $V_{\tilde{\u}}(x)=\tilde{\u}\beta ((\alpha\times \alpha)\cdot \mathit{N}(x)\delta_{\S}$.  Note that in dimension 2, $V_{\tilde{\u}}$  coincides with the electrostatic $\delta$-potential, and in dimension 3 it gives rise to a different $\delta$-potential, see \cite{BB} for more details 
 on the spectral properties of $(\mathcal{H}+V_{\tilde{\u}})$.
 \end{remark}
 
Recall  the matrix  $\g_5$ defined in \eqref{matt}, we set $V_\z=\z\g_5$, for all  $\z\in\rr$. To our knowledge,   the potential $V_\z$  does not seem to have a physical interpretation, but mathematically, it has the same characteristics as the electrostatic potential when $\z=\pm2$; cf. Remark \ref{like the electro}.

Unless otherwise specified, throughout this section we assume that $\S$ satisfies  the assumption $(\mathrm{H1})$, and we consider the Dirac operator $\mathcal{H}_{\z,\u}$ defined formally by
  \begin{align}
\mathcal{H}_{\z,\u} =  \mathcal{H} + V_{\z,\u}= \mathcal{H} +\left( \z\g_5 + i\u\beta(\alpha\cdot \mathit{N})\right)\delta_{\S}, \quad \z,\u\in\rr.
\end{align} 
 
 Comparing with the operators studied before, the operator $\mathcal{H}_{\z,\u}$ is very different. Indeed, because of the presence of anomalous magnetic  potential, several commutativity properties  are no longer true in this case. In addition, $\mathcal{H}_{0,\u}$ (i.e $\z=0$) has the particularity of combining two important phenomena that we have seen before. In fact, as it was  indicated in the introduction,   in the critical case,  $\mathcal{H}_{0,\pm2}$ is essentially self-adjoint and $\S$ becomes impenetrable; see Theorem \ref{main7} below. 
\newline 
  
Now, given $z\in\cc\setminus\left((-\infty,-m]\cup[m,\infty)\right)$, we define the operators $\Lambda^{z}_{\pm}$ as follows:
\begin{align}
\Lambda^{z}_{\pm}=\frac{1}{\z^2+\u^2}(\z\g_5+ i\u\beta(\alpha\cdot \mathit{N}) )\pm\mathit{C}^{z}_{\S}.
\end{align}
Since $i\beta(\alpha\cdot N)$ is $\mathit{C}^{1}$-smooth and symmetric, it follows that $\Lambda^{z}_{\pm}$ are bounded   from  $\mathit{L}^{2}(\S)^4$ onto itself (respectively from $\mathit{H}^{1/2}(\S)^4$ onto itself). Moreover,  $\Lambda^{z}_{\pm}$ are self-adjoint on $\mathit{L}^{2}(\S)^4$,  for all $z\in(-m,m)$.  

 Now, using the same notations as in Section \ref{sec1},  the Dirac operator $\mathcal{H}_{\z,\u}$ (acting in $ \mathit{L}^2(\rr^3)^4$) is defined on the domain
 \begin{align}\label{dom6}
\mathrm{dom}(\mathcal{H}_{\z,\u})=\left\{ u+\Phi[g]: u\in\mathit{H}^1(\rr^3)^4, g\in\mathit{L}^2(\S)^4, t_{\S}u=-\Lambda_{+}[g]\right\},
\end{align}
and the potential $V_{\z,\u}$ is defined by:
\begin{align}\label{}
V_{\z,\u}(\varphi)=\frac{1}{2}(\z\g_5+ i\u\beta(\alpha\cdot \mathit{N}))(\varphi_+ +\varphi_-)\delta_{\S},
\end{align}
with $\varphi_\pm= t_\S u + \mathit{C}^{}_\pm [g]$. Here  $\mathcal{H}_{\z,\u}$ acts in the sense of distributions as $ \mathcal{H}_{\z,\u}(\varphi)= \mathcal{H}u$, for all $\varphi=u+\Phi^{}[g]\in\mathrm{dom}(\mathcal{H}_{\z,\u})$. 

We remind the reader that $\tilde{\Lambda}^{z}_{\pm}$ denotes the continuous extension of $\Lambda^{z}_{\pm}$ defined from $\mathit{H}^{-1/2}(\S)^4$ onto itself. Using the same method as in Section \ref{sec3}, one can show that $ \mathcal{H}_{\z,\u}$ is closable and the domain of the adjoint is given by
 \begin{align}\label{dom6}
\mathrm{dom}(\mathcal{H}^{\ast}_{\z,\u})=\left\{ u+\Phi[g]: u\in\mathit{H}^1(\rr^3)^4, g\in\mathit{H}^{-1/2}(\S)^4, t_{\S}u=-\tilde{\Lambda}_{+}[g]\right\}.
\end{align}
In the following, we briefly discuss the basic spectral properties of $ \mathcal{H}_{\z,\u}$ in the non-critical case, i.e $\z^2+\u^2\neq4$.   The following theorem gathers the most important properties of $ \mathcal{H}_{\z,\u}$. 
\begin{theorem}\label{main6} Let $(\z,\u)\in\rr^2$ be such that $\z^2+\u^2\neq0,4$. Then $ \mathcal{H}_{\z,\u}$ is self adjoint and we have 
 \begin{align}\label{dodom}
\mathrm{dom}(\mathcal{H}_{\z,\u})=\left\{ u+\Phi(g): u\in\mathit{H}^1(\rr^3)^4, g\in\mathit{H}^{1/2}(\S)^4, t_{\S}u=-\Lambda_{+}[g]\right\}\subset\mathit{H}^1(\rr^3\setminus\S)^4.
\end{align}
Moreover, the following statements hold true:
\begin{itemize}
 \item [(i)] Given $a\in(-m,m)$, then $\mathrm{Kr}(\mathcal{H}_{\z,\u}-a)\neq\{0\}$ $\Longleftrightarrow$ $\mathrm{Kr}( \Lambda^{a}_{+})\neq\{0\}$.
\item [(ii)] For all $z\in\rho(\mathcal{H}_{\z,\u})\cap\rho(\mathcal{H})$, it holds that 
 \begin{align}\label{resolvent6}
    (\mathcal{H}_{\z,\u} -z)^{-1}= (\mathcal{H} -z)^{-1} - \Phi^{z}(\Lambda^{z}_{+} )^{-1} (\Phi^{\overline{z}})^\ast.
\end{align}
\item [(iii)] $\mathrm{Sp}_{\mathrm{ess}}(\mathcal{H}_{\z,\u})=(-\infty,-m]\cup [m,+\infty)$. 
 \item[(iv)] $\mathrm{Sp}_{\mathrm{disc}}(\mathcal{H}_{\z,\u})\cap(-m,m)$ is finite.
\end{itemize}
\end{theorem}
Recall that  $[ A, B]= AB -BA$ is the usual commutator bracket. Before giving the proof of the above theorem, we need the following proposition.
\begin{proposition}\label{Pr6} Let $z\in\cc\setminus\left((-\infty,-m]\cup[m,\infty)\right)$. Then, the  commutator $[\beta(\alpha\cdot \mathit{N}), \mathit{C}^{z}_{\S}]$ gives rise to a bounded operator 
\begin{align}
[\beta(\alpha\cdot \mathit{N}), \mathit{C}^{z}_{\S}]:\mathit{H}^{-1/2}(\S)^4\rightarrow \mathit{H}^{1/2}(\S)^4.
\end{align}
In particular, $[\beta(\alpha\cdot \mathit{N}), \mathit{C}^{z}_{\S}]$ is compact in $ \mathit{L}^{2}(\S)^4$.
\end{proposition}
\textbf{Proof.} Let $x,y\in\S$, using \eqref{clifford} a trivial verification shows that 
\begin{align}\label{cliiden}
\begin{split}
\beta(\alpha\cdot\mathit{N}(x))(\alpha\cdot(x-y))-(\alpha\cdot(x-y))\beta(\alpha\cdot\mathit{N}(y))=& (\alpha\cdot(x-y))\beta(\alpha\cdot(\mathit{N}(y)-\mathit{N}(x))\\
&+ (\mathit{N}(x)\cdot(x-y))\beta.
\end{split}
\end{align}
 Now, let $g\in\mathit{L}^{2}(\S)^4$, then using the identity \eqref{cliiden}, similar arguments to those of Lemma \ref{commutator} yield 
\begin{align}
[\beta(\alpha\cdot \mathit{N}), \mathit{C}^{z}_{\S}][g](x)= z\beta[(\alpha\cdot \mathit{N}), S^z][g](x)-m\{(\alpha\cdot \mathit{N}), S^z\}[g](x)+T_z[g](x), 
\end{align}
where the integral representation of $T_z$ is given by
\begin{align}
T_z[g](x)= \int_{\S}K_z(x,y)g(y)dS(y),
    \end{align}
with
\begin{align*}
 K_z(x,y)= \beta \frac{e^{i\sqrt{z^2-m^2}|x-y|}}{4\pi |x-y|^3}(1-i\sqrt{z^2-m^2}|x-y|)\bigg(& (\alpha\cdot(x-y))(\alpha\cdot(\mathit{N}(x)-\mathit{N}(y))\\
 &- 2\left(\mathit{N}(x)\cdot(x-y) \right)\mathit{I}_4\bigg).
 \end{align*}
As $\mathit{N}$ is $\mathit{C}^{1}$-smooth and $S^z$ is bounded from $\mathit{H}^{-1/2}(\S)^4$ to $ \mathit{H}^{1/2}(\S)^4$, it follows that $\beta[(\alpha\cdot \mathit{N}), S^z]$ and $\{(\alpha\cdot \mathit{N}), S^z\}$ are bounded from $\mathit{H}^{-1/2}(\S)^4$ to $ \mathit{H}^{1/2}(\S)^4$. Now, that  $T_z$ is bounded from $\mathit{H}^{-1/2}(\S)^4$ to $ \mathit{H}^{1/2}(\S)^4$ is a direct consequence of Lemma \ref{commutator}. This completes the proof of the first statement, the second statement is a consequence of the Sobolev embedding.\qed
\newline

Now we are in position to prove Theorem \ref{main6}.

\textbf{Proof of Theorem \ref{main6}.} Fix $z\in\cc\setminus\left((-\infty,-m]\cup[m,\infty)\right)$, then a simple computation yields 
\begin{align*}
\Lambda_{\mp}^{z}\Lambda_{\pm}^{z}&=\frac{1}{\z^2+\u^2}-\frac{1}{4}-\mathit{C}^{z}_{\S} (\alpha\cdot \mathit{N})\lbrace \alpha\cdot \mathit{N}, \mathit{C}^{z}_{\S}\rbrace \pm\frac{\z}{\z^2+\u^2}[\g_5,\mathit{C}^{z}_{\S}] \pm\frac{i\u}{\z^2+\u^2} [\beta(\alpha\cdot \mathit{N}), \mathit{C}^{z}_{\S}].
\end{align*}
Now, it is easy to that $[\g_5,\mathit{C}^{z}_{\S}] =2m\g_5\beta S^z$. Using this, it follows that
\begin{align}
\begin{split}
\Lambda_{\mp}^{z}\Lambda_{\pm}^{z}=&\frac{1}{\z^2+\u^2}-\frac{1}{4}-\mathit{C}^{z}_{\S} (\alpha\cdot \mathit{N})\lbrace \alpha\cdot \mathit{N}, \mathit{C}^{z}_{\S}\rbrace\\
& \pm\frac{2m\z}{\z^2+\u^2}\g_5\beta S^z\pm \frac{i\u}{\z^2+\u^2} [\beta(\alpha\cdot \mathit{N}), \mathit{C}^{z}_{\S}].\label{M1}
\end{split}
\end{align}
As  $\z^2+\u^2\neq0,4$, using Lemma \ref{commutator}, Proposition \ref{Pr6}  and applying the same method as in the proof of Theorem \ref{main1} we obtain that 
 \begin{align*}
\mathrm{dom}(\mathcal{H}^{\ast}_{\z,\u})=\mathrm{dom}(\mathcal{H}_{\z,\u})=\left\{ u+\Phi[g]: u\in\mathit{H}^1(\rr^3)^4, g\in\mathit{H}^{1/2}(\S)^4, t_{\S}u=-\Lambda_{+}[g]\right\}.
\end{align*}
Thus, $\mathcal{H}_{\z,\u}$ is self-adjoint  and $\mathrm{dom}(\mathcal{H}_{\z,\u})\subset \mathit{H}^1(\rr^3\setminus\S)^4$. 
Assertions $(\mathrm{i})$, $(\mathrm{ii})$  and $(\mathrm{iii})$ can be proved as in Proposition \ref{BS-K}. Assertion $(\mathrm{iv})$ is a consequence of the Sobolev injection. Indeed, one can easily adapt  the proof of \cite[Theorem 4.1 (ii)]{BEHL2}  and show that $\mathrm{Sp}_{\mathrm{disc}}(\mathcal{H}_{\z,\u})\cap(-m,m)$ is finite. We omit the details.\qed
\newline


 In the following theorem, we discuss the self-adjointness of $\mathcal{H}_{\z,\u}$ in the critical case, i.e $\z^2+\u^2=4$. We mention that assertions $(\mathrm{a})$ and $(\mathrm{c})$ have already been proved in \cite{MH}, where the author studied the inner part of $ \mathcal{H}_{0,\pm2}$ which acts on $\O_+$, known as the Dirac operator with zig-zag boundary conditions, we refer to  \cite{CLMT} for the two-dimensional case. 
\begin{theorem}\label{main7} Let $(\z,\u)\in\rr^2$ be such that $\z^2+\u^2=4$ , then $ \mathcal{H}_{\z,\u}$ is essentially  self -adjoint and we have
 \begin{align*}
\mathrm{dom}(\overline{\mathcal{H}_{\z,\u}})=\left\{ u+\Phi[g]: u\in\mathit{H}^1(\rr^3)^4, g\in\mathit{H}^{-1/2}(\S)^4, t_{\S}u=-\tilde{\Lambda}_{+}[g]\right\}.
\end{align*}
Moreover, the following assertions hold true:
\begin{itemize}
 \item [(i)]  $a\in\mathrm{Sp}(\overline{\mathcal{H}_{\z,\u}})\Longleftrightarrow -a\in\mathrm{Sp}(\overline{\mathcal{H}_{\z,\u}})$.
 \item [(ii)] For all $z\in\cc\setminus\rr$, the operator $ \tilde{\Lambda}^{z}_{+}$ is bounded invertible from $\mathit{H}^{-1/2}(\S)^4$ to $\mathit{H}^{1/2}(\S)^4$ and we have 
 \begin{align}\label{resolventsect6}
    (\overline{\mathcal{H}_{\z,\u}} -z)^{-1}= (\mathcal{H} -z)^{-1} - \Phi^{z}(\tilde{\Lambda}^{z}_{+} )^{-1} (\Phi^{\overline{z}})^\ast.
\end{align}
\item [(iii)]  If $\z=0$, then $\S$ becomes impenetrable and it holds that 
\begin{align}\label{sommedir2}
\overline{\mathcal{H}_{0,\u}}= \mathcal{H}^{\Omega_+}_{\u}\oplus\mathcal{H}^{\Omega_-}_{\u}=\left(-i\alpha\cdot \nabla+m\beta\right)\oplus\left(-i\alpha\cdot \nabla+m\beta\right),
\end{align}
 where $\mathcal{H}^{\Omega_\pm}_{\u}$ are the self-adjoint Dirac operators defined on
\begin{align*}
\mathrm{dom}(\mathcal{H}^{\Omega_\pm}_{\u})&=\left\{ \varphi_{\pm}\in\mathit{L}^2(\Omega_\pm)^4:  (\alpha\cdot\nabla) \varphi_\pm\in\mathit{L}^2(\Omega_{\pm})^4 \text{ and } P_{\mp,\u} t_{\S}\varphi_{\pm}=0 \right\},
\end{align*} 
where the boundary condition has to be understood as an equality in $\mathit{H}^{-1/2}(\S)^4$, and $P_{\pm,\u}$ are the projectors defined by
\begin{align}
P_{\pm,\u}= \frac{1}{2}\left(\mathit{I}_4 \pm\frac{\u}{2}\beta\right).
\end{align}
Furthermore, we have 
\begin{itemize}
\item [(a)] $-m$ and $m$ are eigenvalues of $\overline{\mathcal{H}_{0,\u}}$ with infinite multiplicities. 
\item [(b)] $\mathrm{Sp}(\overline{\mathcal{H}_{0,\u}})=(-\infty,-m]\cup [m,+\infty)$. 
 \item[(c)] There is a sequence $(\l_j(m))_{j\in\mathbb{N}}\subset \mathrm{Sp}(\overline{\mathcal{H}_{0,\u}})$ such that each  $\l_j(m)$ is an eigenvalue with finite multiplicity, with $\l_j(m)^2>m^2$ for all $j\in\mathbb{N}$, and $\l_j(m)^2\longrightarrow \infty$ as $j\longrightarrow \infty$.
\end{itemize}
\end{itemize}
\end{theorem}

\textbf{Proof.} Let us show the first statement. The proof is a relatively straightforward modification of the technique used in the proof of Theorem \ref{main1}. Indeed, as $\mathcal{H}_{\z,\u}$ is closable the only thing left to prove is the inclusion $\mathcal{H}^{\ast}_{\z,\u}\subset\overline{\mathcal{H}_{\z,\u}}$. For this, let  $\varphi=u+\Phi[g]\in  \mathrm{dom}(\mathcal{H}^{\ast}_{\z,\u}) $ and let $(h_j)_{j\in\mathbb{N}}\subset \mathit{H}^{1/2}(\S)^4$ be a sequence of functions that converges to $g$ in $\mathit{H}^{-1/2}(\S)^4$. Set 
   \begin{align}\label{suite6}
   g_j:= \frac{1}{2}\left(\z\g_5+i\u\beta(\alpha\cdot\mathit{N})\right) \left( \tilde{\Lambda}_{+}[g]+\Lambda_{-}[h_j ]\right) ,\quad\forall j\in\mathbb{N}.
   \end{align} 
   Clearly,  $(g_j)_{j\in\mathbb{N}}\subset \mathit{H}^{1/2}(\S)^4$. Since $\Lambda_+$ is bounded from $\mathit{H}^{1/2}(\S)^4$ onto itself, we then get that $ (\Lambda_{+}[g_j])_{j\in\mathbb{N}}\subset \mathit{H}^{1/2}(\S)^4$. Now, remark that 
   \begin{align}
   -\frac{1}{2}\left(\z\g_5+i\u\beta(\alpha\cdot\mathit{N})\right)\tilde{\Lambda}_{-}[g]= -g  +\frac{1}{2}\left(\z\g_5+i\u\beta(\alpha\cdot\mathit{N})\right)\tilde{\Lambda}_{+}[g].
   \end{align}
  Using this, it follows that 
   \begin{align}\label{suite7}
   g_j:= g -\frac{1}{2}\left(\z\g_5+i\u\beta(\alpha\cdot\mathit{N})\right)  \tilde{\Lambda}_{-}[h_j -g],\quad\forall j\in\mathbb{N}.
   \end{align}
  As $\tilde{\Lambda}_-$ is bounded from $\mathit{H}^{-1/2}(\S)^4$ onto itself, it follows that  $g_j \xrightarrow[j\to\infty]{} g $  in $ \mathit{H}^{-1/2}(\S)^4$. Moreover, we have 
   \begin{align}\label{suite8}
    \tilde{\Lambda}_{+}[g_j- g]&= -\frac{1}{2}\tilde{\Lambda}_{+}\left(\z\g_5+i\u\beta(\alpha\cdot\mathit{N})\right)  \tilde{\Lambda}_{-}[g-h_j ]= \left( \tilde{\Lambda}_{+}  \tilde{\Lambda}_{-} \tilde{\Lambda}_{-} +\Lambda_{+}  \tilde{\Lambda}_{+} \tilde{\Lambda}_{-}  \right) [h_j-g ].
   \end{align}
   From Lemma \ref{commutator}, Proposition \ref{Pr6} and \eqref{M1} it follows that  $\tilde{\Lambda}_{\pm}\tilde{\Lambda}_{\mp}$ are bounded from $\mathit{H}^{-1/2}(\S)^4$ to $\mathit{H}^{1/2}(\S)^4$. Therefore,  \eqref{suite8} yields that 
 \begin{align}\label{cv2}
\Lambda_{+}[g_j] \xrightarrow[j\to\infty]{} \tilde{\Lambda}_{+}[g] ,  \text{ in } \mathit{H}^{1/2}(\S)^4.
\end{align}
Let 
$$v_j= E\left(\frac{1}{2} \tilde{\Lambda}_+\left(\z\g_5+i\u\beta(\alpha\cdot\mathit{N})\right)  \tilde{\Lambda}_{-}[h_j-g] \right) \in \mathit{H}^{1}(\rr^3)^4, \quad \forall j\in\mathbb{N},$$
and define  $\varphi_j:= u_j +\Phi[g_j]$, where $u_j=u-v_j$.  It is clear that $u_j\in\mathit{H}^1(\rr^3)^4 $ and $t_\S u_j= -\Lambda_+[g_j]\in \mathit{H}^{1/2}(\S)^4$, hence $(\varphi_j)_{j\in\mathbb{N}}\subset \mathrm{dom}(\mathcal{H}_{\z,\u})$. Moreover, since $(h_j)_{j\in\mathbb{N}}$ (resp $(g_j)_{j\in\mathbb{N}}$)  converges to $g$ in $\mathit{H}^{-1/2}(\S)^4$ as $j\longrightarrow \infty$, using the continuity of $\tilde{\Lambda}_{\pm}\tilde{\Lambda}_{\mp}$ it follows that  $(\varphi_j,\mathcal{H}_{\z,\u}\varphi_j) \xrightarrow[j\to\infty]{}(\varphi, \mathcal{H}^{\ast}_{\z,\u}\varphi)$ in $\mathit{L}^2(\rr^3)^4$. Therefore $ \mathcal{H}^{\ast}_{\z,\u}\subset\overline{\mathcal{H}_{\z,\m}}$ and hence $ \overline{\mathcal{H}_{\z,\u}}$ is self-adjoint with   
  \begin{align}
\mathrm{dom}(\overline{\mathcal{H}_{\z,\u}})=\left\{ u+\Phi[g]: u\in\mathit{H}^1(\rr^3)^4, g\in\mathit{H}^{-1/2}(\S)^4, t_{\S}u=-\tilde{\Lambda}_{+}[g]\right\}.
\end{align}
 This finishes the proof of the first statement. In order to continue  the proof of the theorem we use the definition of $\mathrm{dom}(\overline{\mathcal{H}_{\z,\u}})$ with transmission condition.   As in Definition \ref{def2},  using the Plemelj-Sokhotski formula, one can show that $\overline{\mathcal{H}_{\z,\u}}$ acts in the sense of distributions as 
 
 \begin{align}
\overline{\mathcal{H}_{\z,\u}}\varphi= \left(-i\nabla\cdot\alpha+m\beta\right)\varphi_+\oplus\left(-i\nabla\cdot\alpha+m\beta\right)\varphi_-,
 \end{align}
 for $\varphi= (\varphi_+,\varphi_-)\in\mathit{L}^2(\rr^3)^4$ such that $ (\alpha\cdot\nabla) \varphi_\pm\in\mathit{L}^2(\Omega_{\pm})^4$ and satisfies the following transmission condition in $\mathit{H}^{-1/2}(\S)^4$:
   \begin{align*}
    \left(\frac{1}{2}(\z\g_5+ i\u\beta(\alpha\cdot \mathit{N}))+ i(\alpha\cdot \mathit{N})\right)t_{\S}\varphi_+ =-\left(\frac{1}{2}(\z\g_5+ i\u\beta(\alpha\cdot \mathit{N}))- i(\alpha\cdot \mathit{N})\right)t_{\S}\varphi_-.
    \end{align*}
 Now, let us show item $(\mathrm{i})$, for that recall the operator $\mathcal{C}$ defined in \eqref{matt}. Then, a trivial computation yields that 
 \begin{align*}
 \varphi\in\mathrm{dom}(\overline{\mathcal{H}_{\z,\u}})&\Longleftrightarrow  \mathcal{C}[\varphi]\in\mathrm{dom}(\overline{\mathcal{H}_{\z,\u}}).
   \end{align*}
   Since for all $u\in\mathit{L}^2(\rr^3)^4$, we have 
    \begin{align*}
 \mathcal{C}[(- i \alpha \cdot\nabla + m\beta)u]&=-(- i \alpha \cdot\nabla + m\beta)\mathcal{C}[u], 
  \end{align*}
 it follows that $a$ belongs to $\mathrm{Sp}(\overline{\mathcal{H}_{\z,\u}})$ if and only if $-a$ belongs to $\mathrm{Sp}(\overline{\mathcal{H}_{\z,\u}})$, which yields  $(\mathrm{i})$. Item $(\mathrm{ii})$ follows in the same way as Proposition \ref{BS-K}. To prove item $(\mathrm{iii})$, observe that 
  \begin{align*}
 \begin{split}
\mathrm{dom}(\overline{\mathcal{H}_{0,\u}})=\bigg\{ \varphi=(\varphi_+,\varphi_-)\in\mathit{L}^2(\Omega_+)^4\oplus&\mathit{L}^2(\Omega_-)^4 : (\alpha\cdot\nabla) \varphi_\pm\in\mathit{L}^2(\Omega_{\pm})^4 \text{ and }\\
& i(\alpha\cdot N)P_{-,\u}t_{\S}\varphi_+ =i(\alpha\cdot N)P_{+,\u}t_{\S}\varphi_- \bigg\}.
\end{split}
\end{align*}  
Since $P_{\pm,\u}$ are projectors, we deduce that a function $\varphi=(\varphi_+,\varphi_-)\in\mathit{L}^2(\Omega_+)^4\oplus\mathit{L}^2(\Omega_-)^4$ with $ (\alpha\cdot\nabla) \varphi_\pm\in\mathit{L}^2(\Omega_{\pm})^4$ belongs to $\mathrm{dom}(\overline{\mathcal{H}_{0,\u}})$ if and only if $ P_{\mp,\u} t_{\S}\varphi_{\pm}=0$ holds in $\mathit{H}^{-1/2}(\S)^4$. Therefore, $\S$ becomes impenetrable and  the decomposition \eqref{sommedir2} holds true. 

In the rest of the proof we assume that $\u=2$, the case $\u=-2$ can be recovered with the same arguments.
Let us show assertion $(\mathrm{a})$. For that, we first show that $-m$ is an eigenvalue of $\mathcal{H}^{\Omega_+}_{2}$ with infinite multiplicity, and hence of $\overline{\mathcal{H}_{0,2}}$. Observe that, for any $\varphi=(\varphi_1,\varphi_2)^{\top}\in\mathrm{dom}(\mathcal{H}^{\Omega_+}_{2})$, we have $\varphi_1\in\mathit{H}^1_0(\O_+)^2$ and $(\sigma\cdot\nabla)\varphi_2\in\mathit{L}^2(\O_+)^2$. Let $\psi\in\mathit{C}^{2}(\O_+)^2$ be a harmonic function with respect to $\sigma\cdot\nabla$, i.e $(\sigma\cdot\nabla)\psi=0$ in $\O_+$, and set $\varphi=(0,\psi)^{\top}$.Then, it is clear that $\varphi\in\mathrm{dom}(\mathcal{H}^{\Omega_+}_{2})$ and we have
\begin{align*}
(\mathcal{H}^{\Omega_+}_{2}+m)\varphi=  \begin{pmatrix}
-i(\sigma\cdot\nabla)\psi\\
0
\end{pmatrix}  +m\begin{pmatrix}
2&0\\
0&0
\end{pmatrix} \begin{pmatrix}
0\\
\psi
\end{pmatrix}=0.
\end{align*}
As the set of harmonic functions with respect to $(\sigma\cdot\nabla)$ is infinite dimensional,  we get that $-m$ is an eigenvalue of  $\overline{\mathcal{H}_{0,2}}$ with infinite multiplicity. By $(\mathrm{i})$ we also have that   $m$ is an eigenvalue of  $\overline{\mathcal{H}_{0,2}}$ with infinite multiplicity, which proves assertion $(\mathrm{a})$.

 Now we are going to prove $(\mathrm{b})$ and $(\mathrm{c})$, for that we consider the following Dirac operators
\begin{align}\label{dautreope}
D^{\Omega_\pm}_{2}\psi=\left(-i\alpha\cdot \nabla+m\beta\right)\psi,\quad \psi\in \mathrm{dom}(D^{\Omega_\pm}_{2})=\left\{ \varphi_{\pm}\in\mathit{H}^1(\Omega_\pm)^4: P_{\mp,2} t_{\S}\varphi_{\pm}=0 \right\}.
\end{align} 
Then, one can easily verify that $D^{\Omega_\pm}_{2}$ are symmetric and closable operators. Moreover, it holds that  $\overline{D^{\Omega_\pm}_{2}}=\mathcal{H}^{\Omega_\pm}_{2}$. 
Indeed, denote by $\mathcal{Q}^{\Omega_\pm}_{2}$ the quadratic form associated to $(D^{\Omega_\pm}_{2})^2$. Given $\varphi\in \mathrm{dom}(D^{\Omega_\pm}_{2})$, using the Green's formula and the boundary conditions, it easily follows that:
\begin{align}\label{QuadraticQ}
\mathcal{Q}^{\Omega_\pm}_{2}[\varphi]=&\Vert (\alpha\cdot\nabla)\varphi\Vert^{2}_{\mathit{L}^2(\Omega_\pm)^4} +m^2\Vert \varphi\Vert^{2}_{\mathit{L}^2(\O_\pm)^4}.
\end{align}
Hence, we get $\mathcal{Q}^{\Omega_\pm}_{2}[\varphi]\geqslant m^2\Vert \varphi\Vert^{2}_{\mathit{L}^2(\O_\pm)^4}$.
Thus $(D^{\Omega_\pm}_{2})^2$ is lower semi-bounded. Therefore, by \cite[Theorem 6.3.2]{EE} it follows that $(\mathcal{H}^{\Omega_\pm}_{2})^2$ is the Friedrichs extension of $(D^{\Omega_\pm}_{2})^2$  and it holds that 
 \begin{align*}
 \mathrm{Sp}(\mathcal{H}^{\Omega_\pm}_{2})\subset (-\infty,-m]\cup [m,+\infty).
 \end{align*}
Now, let $(-\Delta^{\Omega_\pm})$ be the Dirichlet realization of $(-\Delta)$ in $\O_\pm$, with domain $\mathit{H}^2(\O_\pm)\cap\mathit{H}^1_0(\O_\pm)$. Using the Weyl's theorem and the fact that $\mathit{H}^1_0(\O_+)$ is compactly embedded in  $\mathit{L}^2(\O_+)$, it is not hard to show that 
 \begin{align}\label{lesspec}
 \begin{split}
 \mathrm{Sp}(-\Delta^{\Omega_-}+m^2)&=  [m^2,+\infty),\\
 \mathrm{Sp}(-\Delta^{\Omega_+}+m^2)&=\mathrm{Sp}_{\mathrm{disc}}(-\Delta^{\Omega_+}+m^2)=\{ m^2+\l_j, j\in\mathbb{N}\},
 \end{split}
 \end{align}
 with $\l_j>0$ for all $j\in\mathbb{N}$, and $\l_j\longrightarrow \infty$ as $j\longrightarrow \infty$.  Using the boundary condition it follows that  
 \begin{align}\label{H0}
\varphi= \begin{pmatrix}
 \varphi_1\\
\varphi_2
\end{pmatrix}\in  \mathrm{dom}(D^{\Omega_+}_{2})\Longrightarrow \varphi_2\in \mathit{H}^1_0(\O_+)^2, \quad \varphi= \begin{pmatrix}
 \varphi_1\\
\varphi_2
\end{pmatrix}\in  \mathrm{dom}(D^{\Omega_-}_{2})\Longrightarrow \varphi_1\in \mathit{H}^1_0(\O_-)^2.
  \end{align}
 Denote by $\tilde{\mathcal{Q}}^{\Omega_\pm}$ the quadratic form associated to $(-\Delta^{\Omega_\pm}+m^2)\mathit{I}_2$. Using \eqref{H0} and the Green's formula, from \eqref{QuadraticQ} it follows that 
\begin{align}\label{QuadraticQ1}
\begin{split}
\mathcal{Q}^{\Omega_+}_{2}[\varphi]&=\Vert (\sigma\cdot\nabla)\varphi_1\Vert^{2}_{\mathit{L}^2(\Omega_+)^2} +m^2\Vert \varphi_1\Vert^{2}_{\mathit{L}^2(\O_+)^2}+\tilde{\mathcal{Q}}^{\Omega_+}[\varphi_2], \quad \forall\varphi\in\mathrm{dom}(D^{\Omega_+}_{2}),\\
\mathcal{Q}^{\Omega_-}_{2}[\varphi]&=\Vert (\sigma\cdot\nabla)\varphi_2\Vert^{2}_{\mathit{L}^2(\Omega_-)^2} +m^2\Vert \varphi_2\Vert^{2}_{\mathit{L}^2(\O_-)^2}+\tilde{\mathcal{Q}}^{\Omega_-}[\varphi_1],\quad \forall\varphi\in\mathrm{dom}(D^{\Omega_-}_{2}).
\end{split}
\end{align}
Thus, \eqref{lesspec} together with assertion $(\mathrm{i})$ yield that  $\l_j(m)=\pm \sqrt{m^2+\l_j}$ is an eigenvalue of $\overline{\mathcal{H}_{0,2}}$  with finite multiplicity, and we have
 \begin{align}
 \mathrm{Sp}(\overline{\mathcal{H}_{0,2})}= (-\infty,-m]\cup [m,+\infty),
 \end{align}
which yields  $(\mathrm{b})$ and $(\mathrm{c})$, and achieves the proof of theorem. \qed
\newline

We finish this paper by pointing out  the following remarks.
\begin{remark}\label{like the electro} Let $\z=\pm2$ and let $\overline{\mathcal{H}_{\z,0}}$ be as in Theorem \ref{main7}. Given $(\varphi_+,\varphi_-)\in\mathrm{dom}(\overline{\mathcal{H}_{\z,0}})$, we write $\varphi_\pm=(\varphi_{\pm,1},\varphi_{\pm,2})^{\top}$.  Then, one can write the transmission condition as follows:
\begin{align}
t_{\S}\varphi_{+,1}=\frac{i\z}{2}(\sigma\cdot\mathit{N})t_{\S}\varphi_{-,2},\quad t_{\S}\varphi_{+,2}=\frac{i\z}{2}(\sigma\cdot\mathit{N})t_{\S}\varphi_{-,1}.
\end{align}
Thus, we deduce that $\overline{\mathcal{H}_{\z,0}}$ coincide with the Dirac operator coupled with the electrostatic $\delta$-interactions of strength $-\z$. Thus, in this sense, one can consider the potential $V_{\z}$ as an electrostatic potential for $\z=\pm2$.
\end{remark}
\begin{remark}  If one assume that $\S$ satisfies the assumption $(\mathrm{H2})$, then $\mathcal{H}_{\z,\u}$ is essentially self-adjoint, when $\z^2+\u^2=4$. In particular, if $\u=0$, then Remark \ref{like the electro} and Theorem \ref{cas deforme} yield that 
\begin{align}
\mathrm{Sp}_{\mathrm{ess}}(\overline{\mathcal{H}_{\pm2,0}})=\big(-\infty,-m\big]\cup\{0\}\cup \big[m,+\infty\big).
\end{align}
 However, if $\z=0$, then there is no embedded eigenvalues in the essential spectrum of $\overline{\mathcal{H}_{0,\pm2}}$, and we have 
 \begin{align}
\mathrm{Sp}(\overline{\mathcal{H}_{0,\pm2}})=\mathrm{Sp}_{\mathrm{ess}}(\overline{\mathcal{H}_{0,\pm2}})=\big(-\infty,-m\big]\cup \big[m,+\infty\big).
\end{align}
   \end{remark}

\section*{Acknowledgement} I would like to thank my PhD supervisors Vincent Bruneau and Luis Vega for their encouragement, and for their precious discussions and advices during the preparation of this paper.  This project is based upon work supported by the Government of the Basque Country under Grant PIFG$18/06$ "ERC Grant: Harmonic Analysis and Differential Equations", researcher in charge: Luis Vega.


\end{document}